%% file: ArxivShahidi.tex
\input amstex
\documentstyle{amams} 
\document
\annalsline{155}{2002}
\input amssym.def
\input amssym.tex
 
\startingpage{837}
\font\tenrm=cmr10

\catcode`\@=11
\font\twelvemsb=msbm10 scaled 1100

\font\ninemsb=msbm10 scaled 800
\newfam\msbfam
\textfont\msbfam=\twelvemsb  \scriptfont\msbfam=\ninemsb
  \scriptscriptfont\msbfam=\ninemsb
\def\msb@{\hexnumber@\msbfam}
\def\Bbb{\relax\ifmmode\let\next\Bbb@\else
 \def\next{\errmessage{Use \string\Bbb\space only in math
mode}}\fi\next}
\def\Bbb@#1{{\Bbb@@{#1}}}
\def\Bbb@@#1{\fam\msbfam#1}
\catcode`\@=12

 \catcode`\@=11
\font\twelveeuf=eufm10 scaled 1100
\font\teneuf=eufm10
\font\nineeuf=eufm7 scaled 1100
\newfam\euffam
\textfont\euffam=\twelveeuf  \scriptfont\euffam=\teneuf
  \scriptscriptfont\euffam=\nineeuf
\def\euf@{\hexnumber@\euffam}
\def\frak{\relax\ifmmode\let\next\frak@\else
 \def\next{\errmessage{Use \string\frak\space only in math
mode}}\fi\next}
\def\frak@#1{{\frak@@{#1}}}
\def\frak@@#1{\fam\euffam#1}
\catcode`\@=12

\def\bL{{\Bbb L}}
\def\tilde{\widetilde}
\def\la{\langle}
\def\ra{\rangle}

\def\frakh{{\frak h}}
\def\var{\varepsilon}
\def\bG{{\bf G}}
\def\bM{{\bf M}}
\def\bN{{\bf N}}
\def\bP{{\bf P}}
\def\bB{{\bf B}}
\def\bT{{\bf T}}
\def\bU{{\bf U}}

\def\bH{{\bf H}}

\def\bR{{\Bbb R}}
\def\bC{{\Bbb C}}
\def\bQ{{\Bbb Q}}
\def\ds{\displaystyle}
\def\bA{\bold A}
\def\bAB{\Bbb A}
\def\bAA_F{\Bbb A_F}

\def\bZ{\Bbb Z}
\def\sP{\Cal P}

\def\sym{\text{\rm Sym}}

\def\Ind{{\rm Ind}}
\def\Spin{\text{\rm Spin}}

\font\eti=cmr10 scaled \magstep 3
 3 

 3
\font\tits=cmr10

\title{Functorial products for {\eti GL}\lower2pt\hbox{\tits 2}$\times$ {\eti GL}\lower2pt\hbox{\tits 3}\\ and the symmetric cube for {\eti
GL}\lower2pt\hbox{\tits 2}}
\shorttitle{The symmetric cube for GL$_2$} 
\acknowledgements{The first author was partially supported by NSF grant
DMS-9988672, NSF grant DMS-9729992 (at IAS), and by 
Clay Mathematics Institute.  He is a Clay Mathematics Institute Prize Fellow. The second author was partially
supported by NSF grant DMS-9970156 and by Clay 
Mathematics Institute.  He is a Clay Prize Fellow as well as a Guggenheim Fellow.}
 \twoauthors{Henry H. Kim}{Freydoon Shahidi}
   \institutions{University of Toronto, Toronto, Ontario, Canada\\
{\eightpoint {\it E-mail address\/}: henrykim@math.toronto.edu}\\ \vglue6pt
Purdue University, West Lafayette, IN\\
{\eightpoint {\it E-mail address\/}: shahidi@math.purdue.edu}}

 \centerline{\it Dedicated to  Robert P.\ Langlands}

\bigbreak\centerline{\bf Introduction}
\vglue12pt

In this paper we prove two new cases of Langlands functoriality.
The first is a functorial product for cusp forms on ${\rm GL}_2\times {\rm GL}_3$ as automorphic forms on ${\rm GL}_6$, from which we obtain our second case, the long awaited functorial symmetric cube map for cusp forms on ${\rm GL}_2$.
We prove these by applying a recent version of converse theorems of Cogdell and Piatetski-Shapiro to analytic properties of certain $L$-functions obtained from the method of Eisenstein series (Langlands-Shahidi method).
As a consequence we prove the bound 5/34 for Hecke eigenvalues of Maass forms over any number field and at every place, finite or infinite, breaking the crucial bound 1/6 (see below and Sections 7 and 8) towards Ramanujan-Petersson and Selberg conjectures for ${\rm GL}_2$. As noted below, many other applications follow.

To be precise, let $\pi_1$ and $\pi_2$ be two automorphic cuspidal
representations of ${\rm GL}_2(\bAA_F)$ and ${\rm GL}_3(\bAA_F)$, respectively,
where $\bAA_F$ is the ring of ad\`eles of a number field $F$. Write $\pi_1=\hskip-6pt  \mathbold{\otimes}_v\pi_{1v}$ and $\pi_2=\hskip-6pt  \mathbold{\otimes}_v\pi_{2v}$.
For each $v$, finite or otherwise, let $\pi_{1v}\boxtimes\pi_{2v}$ be the irreducible admissible representation of ${\rm GL}_6(F_v)$, attached to $(\pi_{1v},\pi_{2v})$ through the local Langlands correspondence by Harris-Taylor [HT], Henniart [He], and Langlands [La4].
We point out that, if $\varphi_{iv},\ i=1,2$, are the two- and the three-dimensional representations of Deligne-Weil group, parametrizing $\pi_{iv}$,
respectively, then $\pi_{1v}\boxtimes\pi_{2v}$ is attached to the
 six-dimensional representation $\varphi_{1v}\otimes\varphi_{2v}$.
Let $\pi_1\boxtimes\pi_2=\hskip-6pt  \mathbold{\otimes}_v (\pi_{1v}\boxtimes\pi_{2v})$.
Our first result\break (Theorem 5.1) is:
\pagebreak

\nonumproclaim{Theorem A} The representation $\pi_1\boxtimes\pi_2$ of ${\rm GL}_6(\bAA_F)$ is automorphic{\rm ,} i.e.{\rm ,}
functorial products {\rm [La3]} for ${\rm GL}_2\times {\rm GL}_3$ exist{\rm .} It is isobaric and more specifically{\rm ,} is irreducibly
induced from unitary cuspidal representations{\rm ,} i.e{\rm .,} $\pi_1\boxtimes\pi_2=
{\rm Ind}\, \sigma_1\otimes\cdots\otimes\sigma_k${\rm ,} where $\sigma_i$\/{\rm '}\/s are
 unitary cuspidal representations of ${\rm GL}_{n_i}(\bAA_F),\ n_i > 1${\rm .}
Moreover{\rm ,} $k=3,\ n_1=n_2=n_3=2${\rm ,}
 or $k=2,\ n_1=2$, $n_2=4$ occur if and only if $\pi_2$ is a twist of ${\rm Ad}(\pi_1)${\rm ,} the adjoint of $\pi_1$ {\rm (}\/see below\/{\rm ),} by
a gr{\rm \"{\it o}}ssencharacter{\rm .}
\endproclaim

We   remark that at the moment, we are unable to characterize the image of this functorial product or completely determine all the pairs for which
the image is cuspidal.

Next, let $\pi$ be a cuspidal representation of ${\rm GL}_2(\bAB_F)$.
Write $\pi=\hskip-6pt  \mathbold{\otimes}_v\pi_v$.
Let ${\rm Ad}\colon {\rm GL}_2(\bC)\to {\rm GL}_3(\bC)$ be the adjoint representation of ${\rm GL}_2(\bC)$.
Then ${\rm Ad}({\rm diag}((a,b))$ $={\rm diag}((ab^{-1},1,a^{-1} b)$.
Let $\varphi_v$ be the two-dimensional representation of the  Deligne-Weil group attached to $\pi_v$ ([Ku], [La4]).
Then ${\rm Ad}(\varphi_v)=Ad\cdot\varphi_v$ is a three-dimensional representation.
Let ${\rm Ad}(\pi_v)$ be the irreducible admissible representation of ${\rm GL}_3(F_v)$ attached to ${\rm Ad}(\varphi_v)$ ([He2], [La4]).
Set ${\rm Ad}(\pi)=\hskip-6pt  \mathbold{\otimes}_v {\rm Ad}(\pi_v)$.
Then in [GeJ], Gelbart and Jacquet proved that ${\rm Ad}(\pi)$ is an automorphic representation of ${\rm GL}_3(\bAB_F)$ which is cuspidal unless $\pi$ is monomial, i.e., $\pi\simeq\pi\otimes\eta$, where $\eta\not= 1$ is a gr\"ossencharacter of $F$.
Observe that, if ${\rm Sym}^2$ is the three-dimensional irreducible representation of ${\rm GL}_2(\bC)$ on symmetric tensors of rank 2, this
implies the same facts about ${\rm Sym}^2(\pi)=\hskip-6pt  \mathbold{\otimes}_v {\rm Sym}^2 (\pi_v)$, where ${\rm Sym}^2(\pi_v)$ is the irreducible admissible
representation of ${\rm GL}_3(F_v)$ attached to ${\rm Sym}^2(\varphi_v)={\rm Sym}^2\cdot\varphi_v$. Moreover ${\rm
Sym}^2(\pi)={\rm Ad}(\pi)\otimes\omega_\pi$, where $\omega_\pi$ is the central character of $\pi$.

We now proceed to the next case and, as before, let $\pi$ be a cuspidal representation of ${\rm GL}_2(\bAA_F)$.
Write $\pi=\hskip-6pt  \mathbold{\otimes}_v\pi_v$ and denote by ${\rm Sym}^3\colon {\rm GL}_2(\bC)\to {\rm GL}_4(\bC)$, the four-dimensional irreducible representation of
${\rm GL}_2(\bC)$ on symmetric tensors of rank 3. Simply put, for each $g\in {\rm GL}_2(\bC)$, ${\rm Sym}^3(g)\in {\rm GL}_4(\bC)$ can be taken to be
the matrix that gives the change in coefficients of an arbitrary homogeneous cubic polynomial in two variables, under the change of variables by $g$.
Again, as before, let $\varphi_v$ be the two-dimensional representation of the Deligne-Weil group attached to $\pi_v$. Then
${\rm Sym}^3(\varphi_v)=\text{ Sym}^3\cdot\varphi_v$ is a four-dimensional representation. Let ${\rm Sym}^3(\pi_v)$ be the irreducible admissible
representation of ${\rm GL}_4(F_v)$ attached to ${\rm Sym}^3(\varphi_v)$ by the local Langlands correspondence (see references above). Set
${\rm Sym}^3(\pi)=\hskip-6pt  \mathbold{\otimes}_v\text{ Sym}^3 (\pi_v)$ which we call the symmetric cube of $\pi$. Applying Theorem A to $\pi_1=\pi$ and $\pi_2=\text{
Ad}(\pi)$, and using the classification theorem [JS1] for ${\rm GL}_6(\bAB_F)$, we obtain (Theorem~6.1):

\nonumproclaim{Theorem B} The representation ${\rm Sym}^3(\pi)$ is an automorphic representation of ${\rm GL}_4(\bAA_F)${\rm .}
It is cuspidal{\rm ,} unless $\pi$ is either of dihedral or  of tetrahedral type{\rm ,} i.e.{\rm ,}
 those attached to dihedral and tetrahedral representations of the Galois group of $\overline F/F$ {\rm (}\/cf.~{\rm [Ge], [La2]).}
In particular{\rm ,} if $F=\bQ$ and $\pi$ is the automorphic representation attached to a nondihedral holomorphic form of weight $\geq 2${\rm ,}
 then ${\rm Sym}^3(\pi)$ is cuspidal{\rm .}
\endproclaim

Due to the large number of applications and   resilience to previously available methods, the existence of the symmetric cube has fascinated a good
number of experts in the field (cf.~[KS1, 2, 3] and references therein, also see [BK]) ever since symmetric squares were established by Gelbart and
Jacquet [GeJ]. Here the map is symplectic (cf.~Section 9) and a much richer geometry is involved, which seems to be the pattern for odd symmetric
powers, making them harder to get than even ones. Moreover, neither of the trace formulas, or for that matter, any other approach, can be used to
prove either of these theorems.

Our present proof of Theorem B is quite surprising.
In fact, we were originally planning to prove the existence of the symmetric cube for ${\rm GL}_2$ by further twisting by forms on ${\rm
GL}_2(\bAB_F)$ which at the time was totally out of reach [KS1]. On the other hand, our present proof of the existence of the  symmetric cube is an
indirect consequence of the functorial product for ${\rm GL}_2\times {\rm GL}_3$, from which the properties of $L$-functions attached to the above
twisting immediately follow. We should point out that one does not need the full
functorial product for ${\rm GL}_2\times {\rm GL}_3$ to prove the existence of the functorial symmetric cube for cusp forms on ${\rm GL}_2(\bAB_F)$ (Remark 6.7).

Theorem B is a very pleasant conclusion to a project   the second author has pursued since 1978. In fact, following Langlands [La6] and his theory of
Eisenstein series [La1], this has led him to develop a machinery [Sh1, 2, 4, 5, 7], whose full power and subtlety (for example cf.~Sections 3-7
of [Sh1] and Sections 2--5 here) are necessary to prove Theorem A from which existence of the symmetric cube, i.e., Theorem B, follows. In fact, as
explained in the next paragraph, this is not possible, unless we bring in two new and crucial results [Ki1], [GeS], but still using the same method.

Theorem A is proved by applying a recent and ingenious version [CP-S3] of converse theorems of Cogdell and Piatetski-Shapiro to analytic properties
of
$L$-functions obtained from our method (cf.~Theorem 3.2 here). While functional equations and their subtle consequences were proved in full
generality earlier in [Sh4], [Sh1], two new and important ingredients are needed. The first one is a crucial observation of Kim (Proposition 2.1 of [Ki1]
or Lemma 7.5 of [La1]), allowing us to utilize holomorphy of highly ramified twists of certain $L$-functions obtained from our method.  In fact, it
was generally believed that it would not be possible to obtain the holomorphy of the $L$-functions using the Langlands-Shahidi method.  But the fact
that cuspidal representations which are not invariant under certain Weyl group elements do not contribute to the residual spectrum, combined
with the Langlands-Shahidi method, gives the holomorphy of corresponding $L$-functions. In view of recent powerful converse theorems [CP-S1],
[CP-S3], this is sufficient. The second is a recent paper of Gelbart-Shahidi [GeS], where boundedness in finite vertical strips for every
$L$-function obtained from our method, is proved. In the present situation, we need to apply the main theorem of [GeS] to four different cases in
the lists in [La6] and [Sh2] (case (iii) of [La6] and the  triple product $L$-function cases $D_5-2,\ E_6-1$, and $E_7-1$ of [Sh2]), all of which but one,
have more than one
$L$-function in their constant terms, requiring us to use the full power and subtlety of Theorem 4.1 of [GeS]  (the number of $L$-functions
$m=1,2,3$, and 4, in these four cases, respectively).

As striking as this result is, it only allows us to prove the existence of a weak lift (Theorem 3.8).
To prove the existence of a lift whose
local components are everywhere those attached by Harris-Taylor [HT],
Henniart [He], and Langlands [La4], through the local Langlands correspondence,
a lot more work is needed (Theorem 5.1).

What we need to do is to prove equalities (3.2.1) and (3.2.2) of Section 3, i.e.,
\vglue-1pt
\centerline{$
L(s,\pi_{1v}\times \pi_{2v}\times\sigma_v)=L(s,(\pi_{1v}\boxtimes \pi_{2v})\times \sigma_v)
$}

\noindent 
and

\centerline{$
\varepsilon(s,\pi_{1v}\times \pi_{2v}\times \sigma_v,\psi_v)=\varepsilon(s,(\pi_{1v}
\boxtimes \pi_{2v})\times \sigma_v,\psi_v),
$}
\vglue4pt
\noindent 
for all irreducible admissible generic representations $\sigma_v$ of
${\rm GL}_n(F_v)$, $1\le n\le 4$.  These equalities are not obvious at all since the
two factors on the left and right are defined using completely different techniques.  The factors on the left are those defined by the triple $L$-functions in our 
four cases
 [Sh1, 7], [La4], while the ones on the right are those of
Rankin-Selberg for ${\rm GL}_6(\Bbb A_F)\times {\rm GL}_n(\Bbb A_F)$, $1\le n\le 4$
[JP-SS], [Sh1].  By the local Langlands correspondence, they are Artin factors.

Using [Sh5] on the equality of $L$-functions defined from our method and those of Artin for
${\rm GL}_k\times {\rm GL}_l$, we can show the above equalities when $\pi_{1v}$, $\pi_{2v}$
are not both supercuspidal representations.  If $v\nmid 2$, then any supercuspidal representation of ${\rm GL}_2(F_v)$ is attached to an induced representation of the corresponding Weil group.  Hence, using quadratic
base change [AC], [La2], we can reduce to the case when $\pi_{1v}$ is not supercuspidal.  We do the same when $v|2$ and $\pi_{1v}$ is not an
extraordinary supercuspidal representation.  Now, suppose $v|2$ and $\pi_{1v}$ is an
extraordinary supercuspidal representation.  Then $\pi_{2v}$ is attached to an
induced representation from a cubic extension (normal or otherwise).  If $n\le 3$, then using either a normal cubic base change [AC], [La2], or a
nonnormal one [JP-SS2, 3], we can reduce to the case when $\pi_{2v}$ is not supercuspidal.  However, due to the fact that the theory of
nonnormal cubic base change for ${\rm GL}_4$ is not available at present, this does not work for $n=4$.  Namely, our argument cannot proceed when
$v|2$,
$\pi_{1v}$ is an extraordinary supercuspidal representation of ${\rm GL}_2(F_v)$, and $\pi_{2v}$ is a supercuspidal representation of ${\rm
GL}_3(F_v)$, attached to an induced representation from a character of a nonnormal cubic extension.  By using Proposition 5.1 of [Sh1], which allows
us to produce cusp forms with prescribed supercuspidal components, we can obtain an irreducible admissible representation $\Pi_v$ of ${\rm
GL}_6(F_v)$, which satisfies the above equalities for $n\le3$, and which differs from $\pi_{1v}\boxtimes \pi_{2v}$ by at most a quadratic
character. Then the self-contained appendix by Bushnell and Henniart [BH], which uses certain subtle results from the theory of types and conductor
for pairs [BHK], [Sh5], proves that in fact $\Pi_v\simeq \pi_{1v}\boxtimes \pi_{2v}$.  To complete the theorem, we need to apply the converse
theorem twice again. (See the argument at the end of Section 5.)

On the other hand, we should point out that it is possible to show that
${\rm Sym}^3(\pi)$ is functorial even at $v|2$ without resorting to the
appendix [BH] (cf.~Remark 6.7).

Section 4 is devoted to a verification of Assumption 1.1 and Conjecture 1.2 (cf.~Section 1 here) in our four cases. They ought to be verified, if our machine is going to work, as they are needed both for our Proposition 2.1, as well as Theorem 4.1 of
[GeS]. The proof relies heavily on repeated application of multiplicativity (Theorem 3.5.3 of [Sh1]) and certain results and ideas from [Ki1, 2,  5],
[CS], [Z].

Our paper concludes with a large number of applications (Sections 7--10).
Many more, some yet to be formulated, are expected (e.g., in dimension  four generalization of Wiles' program [W]).

Sections 7 and 8 are devoted to analytic number theory. The first result establishes a new bound towards Ramanujan-Petersson and Selberg conjectures for ${\rm GL}_2$. More precisely, our Theorem 7.1, which is a consequence of applying
Theorem B to main estimates in [LRS1], states:

\nonumproclaim{Theorem C} Let $\pi$ be a cuspidal representation of ${\rm GL}_2(\bAA_F)${\rm .}
Let $\pi_v$ be a local {\rm (}\/finite or infinite\/{\rm )} spherical component of $\pi${\rm ,} i.e.{\rm ,}
 $\pi_v=\Ind(|\ \ |_v^{s_{1v}}, |\ \ |_v^{s_{2v}}),\break
s_{jv}\in\bC${\rm .} Then $|\text{\rm Re}(s_{jv})|\leq 5/34,\ j=1,2${\rm .}
\endproclaim

When $F=\bQ$ and $\lambda_1(\Gamma)$ is the smallest positive eigenvalue for Laplacian on $L^2(\Gamma\backslash\frak h)$, where $\frak h$ denotes the upper half plane and $\Gamma$ is a congruence subgroup, then Theorem C implies:\ $\lambda_1(\Gamma)\geq 66/289\cong 0.22837$. The earlier bound of 0.21 was due to Luo-Rudnick-Sarnak 
[LRS2]. Observe that $1/7 < 5/34 < 1/7+0.004$.

The fact that the estimate $5/34$ is sharper than $1/6$ is crucial and allows us to prove some fundamental results in analytic number theory.
For example, our Propositions 8.1 and 8.2 on shifted sums [Go] and hyperbolic circle \pagebreak
 problem [I], yield as sharp a result as if one assumes the full
Selberg  Conjecture:\ $\lambda_1(\Gamma)\geq 0.25$. The next crucial bound is 1/12. Section 8 was suggested to us by Peter Sarnak. Further
applications to theory of automorphic forms (e.g., Jacquet's conjecture), using the fact that $5/34 < 1/6$, are also expected.

In Section 9, following a suggestion of Joseph Shalika, we prove a
conditional theorem (Theorem 9.2) on the existence of Siegel modular
cusp forms of weight $3$ for ${\rm GSp}_4(\bAB_{\bQ})$. Here we need the validity of Arthur's multiplicity formula for ${\rm GSp}_4(\bAB_{\bQ})$.
As we remark there (Remark 9.3), this must follow from stabilization of the trace formula for ${\rm GSp}_4 (\bAB_\bQ)$ and the most general form
of the twisted trace formula for ${\rm GL}_4(\bAB_\bQ)$.

Section 10 is devoted to new examples of Artin's conjecture and global Langlands correspondence. In view of the recent progress in Taylor's program
[Tay] on Artin's conjecture for icosahedral representations of global Weil group $W_F$, our Theorem B immediately proves new cases of Artin's
conjecture for  four-dimensional irreducible {\it primitive} representations of $W_F$. Together with certain other new cases, this is recorded as
our Theorem 10.1. We refer to Theorem 10.2 for new cases of Artin's conjecture coming from Theorem A, when the representations are of
icosahedral type. The solvable cases of our theorem are already established by Ramakrishnan in [R3].

Finally, we   refer to [Ki4] for other cases of functoriality obtained using this method, namely, the exterior square lift from ${\rm GL}_4$ to ${\rm
GL}_6$ and the symmetric fourth powers of cusp forms on ${\rm GL}_2$. In fact, when the exterior square lift from ${\rm GL}_4$ to ${\rm GL}_6$
is combined with our symmetric cube, it leads to the existence of symmetric fourth powers of cusp forms on ${\rm GL}_2(\bAA_F)$. 
Further consequences of the existence of these two symmetric powers are collected in [KS3].

There are a number of mathematicians whose comments have influenced this paper.
We first thank Herv\'e Jacquet whose comments as a Comptes Rendus editor on our announcement [KS2] led to the present form of Section 5.  Next,
we  thank Peter Sarnak for his continued encouragement and support for this project. In particular, we  thank him for suggesting the
material in Section 8. The problem in Section 9 was suggested by Joseph Shalika for which we  thank him. The authors   thank him as
well for a discussion which led to formulation of the existence of the  symmetric cube in the form presented in Section 6. This paper owes much, as
well, to Dinakar Ramakrishnan, particularly for many conversations   the first author   had with him during his stay at the Institute for
Advanced Study during the 1999--2000 Special Year.  Thanks are also due to James Cogdell for insights provided by him on his converse theorems
with Piatetski-Shapiro. Our paper concludes with an appendix by Colin Bushnell and Guy Henniart which allows us to remove the final hurdle in
establishing our functorial product. We are grateful to them for\break providing it.

Next, we  thank the Institute for Advanced Study, and particularly the organizers of their Special Year in the Theory of Automorphic Forms and
the $L$-function:\ E.~Bombieri, H.~Iwaniec, R. P.~Langlands, and P.~Sarnak, for providing the support and the setting which led to this progress during
the first author's year visit and the  second author's two month visits.

The second author   thanks Jian-Shu Li and The Hong Kong University of Science and Technology, as well as the  Centre Emile Borel, Institut Henri
Poincar\'e, and the organizers of the Hecke Semester for their support while the material in Section 5 was being developed.  In
fact, it was at IHP where questions raised to  Bushnell and Henniart led to
preparation of their appendix.

Both authors  thank the Clay Mathematics Institute for its support during the final stages of preparation of this paper.

The authors would like to dedicate this paper to Robert Langlands with
admiration. Developing his ideas over the years   made these surprising results possible.

 \vfil
\section{Preliminaries}
 
\vfil
Let $F$ be a number field and denote by $\bAB_F$ its ring of ad\`eles.
For each place $v$ of $F$, let $F_v$ denote the corresponding completion of $F$.
When $v < \infty$, we let $O_v$ denote the ring of integers of $F_v$.
Let $\varpi_v$ be a uniformizing element for the maximal prime ideal $\sP_v$ of $O_v$.
Let $q_v$ be the cardinality of $O_v/\sP_v$ and fix an absolute value
$|\ |_v$ in the class of $v$ such that $|\varpi_v| = q^{-1}_v$.

Let $\bG$ be a quasi-split connected reductive algebraic group over $F$.
Fix a Borel subgroup $\bB$ of $\bG$ over $F$. Write $\bB=\bT\bU$ where $\bT$ is a maximal torus and $\bU$ is the unipotent radical of $\bB$. Let $\bP \supset\bB$ be a parabolic subgroup of $\bG$ over $F$. Write $\bP=\bM\bN$ for a Levi decomposition of $\bP$ with a Levi factor $\bM$ and unipotent radical $\bN$. Fix $\bM$ by assuming $\bM \supset\bT$. Observe that $\bN\subset\bU$.
We let $\bA_0$ be the maximal $F$-split subtorus of $\bT$. Let $W$ be the Weyl group of $\bA_0$ in $\bG$. The choice of $\bU$ determines a set of positive roots for $\bA_0$. Let $\Delta$ be the set of simple roots and denote by $\theta$ the subset of $\Delta$ generating $\bM$.

Throughout this paper, we shall assume $\bP$ is maximal. Then $\Delta\backslash \theta$ is  a singleton. Let $\Delta\backslash\theta= \{\alpha\}$.
The simple root $\alpha$ can be identified as the unique reduced root of $\bA$, the split component of the center of $\bM$, in $\bN$. If $\rho_{\bP}$
is half the sum of roots in $\bN$, we let $\tilde\alpha=\la \rho_\bP,\alpha\ra^{-1}\rho_\bP$ as in [Sh2].

There exists a unique element $\tilde w_0\in W$ such that $\tilde w_0(\theta)\subset\Delta$ while $\tilde w_0(\alpha)$ is a negative root.
We shall fix a representative $w_0$ for $\tilde w_0$ in $K \cap \bG (F)$ as in [Sh7], where $K$ is a fixed maximal compact subgroup of $G=\bG
(\bAA_F)$. Similarly we use $B, T, U, P, M, N, A$, and $A_0$ to denote the corresponding adelic points. Finally, given $v$, and a group $\bH$ over
$F_v$, we use $H_v$ to denote $\bH (F_v)$. We, therefore, have $G_v, B_v$, and so on.

Given a connected reductive algebraic group $\bH$ over $F$, let ${}^L\!H$ be its
$L$-group. We use $\hat{H}$ to denote the connected component ${}^L\!H^0$ of ${}^L\!H$. Considering $\bH$ as a group over each $F_v$, we then denote by ${}^L\!H_v$ its $L$-group over $F_v$, and we have a natural homomorphism
${}^L\!H_v \longrightarrow {}^L\!H$. Let $\eta_v\colon {}^L\!M_v \to {}^L\!M$ be this map for $\bM$.

Let ${}^L\!N$ be the $L$-group of $\bN$ defined naturally in [Bor]. Let ${}^L\!\frak{n}$ be its (complex) Lie algebra, and let $r$ denote the adjoint
action of ${}^L\!M$ on ${}^L\!\frak{n}$. Decompose $r=\bigoplus\limits^m_{i=1} r_i$ to its irreducible subrepresentations, indexed according to values
$\la \tilde\alpha,\beta\ra=i$ as $\beta$ ranges among the positive roots of $\bT$.
More precisely, $X_{\beta^\vee}\in ^L\!\frak n$ lies in the space of $r_i$, if and only if $\la \tilde\alpha,\beta\ra=i$.
Here $X_{\beta^\vee}$ is a root vector attached to the coroot $\beta^\vee$, considered as a root for the $L$-group.
Moreover, $\la\ ,\ \ra$ denotes the Killing form, i.e., for every pair of positive roots $\gamma$ and $\delta$ of $\bT$, $\la \gamma,\delta\ra=2
(\gamma,\delta)/(\delta,\delta)=(\gamma,\delta^\vee)$, where $\delta^\vee$ is the coroot $2\delta/(\delta,\delta)$ attached to $\delta$ (cf.~[Sh2]).

Now, let $\pi =\hskip-6pt  \mathbold{\otimes}_v \pi_v$ be an irreducible unitary globally generic
cuspidal representation of $M= \bM(\bAA_F)$. For each place $v$ and each $i, 1 \leq i \leq m$, let $L(s, \pi_v,r_i \cdot \eta_v)$ and $\varepsilon (s, \pi_v, r_i\cdot \eta_v, \psi_v)$ be the $L$-function and the root number
defined in [Sh1], where $\psi=\hskip-6pt  \mathbold{\otimes}_v\psi_v$ is a nontrivial additive
character of $F\backslash \bAA_F$. When $v=\infty$ or $\pi_v$ is unramified or has only an Iwahori-fixed vector, then they
are exactly those defined by Langlands [La4, 5, 6], [Sh7]. Moreover, if
$$ \align
L(s, \pi, r_i)& = \prod_v L(s, \pi_v, r_i\cdot \eta_v)
\\
\noalign{\noindent 
and}
\noalign{\vskip-4pt}
\varepsilon (s, \pi, r_i) &= \prod_v \varepsilon
(s, \pi_v, r_i\cdot \eta_v,\psi_v),
\\
\noalign{\noindent 
then} 
L(s,\pi,r_i)&=\varepsilon (s,\pi,r_i) L(1-s,\tilde\pi,r_i),\tag1.1 \endalign
$$
where $\tilde\pi$ is the contragredient of $\pi$ (cf.~[Sh1]). We set $r_{i,v}=r_i\cdot \eta_v$.

Next, with notation as in [Sh1, 2], we have the globally induced
representation $I(s, \pi) = I(s \tilde\alpha, \pi)$ from $\pi\otimes \exp
\la s\tilde\alpha,H_P(\ \ )\ra$ as well as the local ones $I(s, \pi_v)$,
induced from $\pi_v\otimes \exp\la s\tilde\alpha,H_{P_v}(\ \ )\ra$, for each $v$.
Let us point out that in our notation, $s=0$ always corresponds to a unitarily induced representation.
We finally recall the global intertwining operator $M(s, \pi)$ defined by
$$
M(s, \pi) f(g) = \int_{N'} f(w^{-1} ng) dn\qquad (g\in G),\tag1.2
$$
where $f\in I (s, \pi)$ and $\bN'$ is the unipotent radical of the standard parabolic subgroup which has the Levi subgroup $\bM^\prime \supset \bT$, generated by $\tilde w_0 (\theta)$. Then
\pagebreak
$$
M(s,\pi) = \hskip-6pt  \mathbold{\otimes}_v M(s \tilde\alpha, \pi_v, w_0),\tag1.3
$$
where the local operators are defined the same way. Finally at each $v$, let $N(s \tilde\alpha, \pi_v, w_0)$ be the normalized operator
$$ \multline
  N(s\tilde\alpha,\pi_v, w_0)\\ = \prod^m_{i=1}\varepsilon
(is, \pi_v, \tilde r_{i,v}, \psi_v)
   L(1 + is, \pi_v, \tilde r_{i,v})
L(i s, \pi_v, \tilde r_{i,v})^{-1}
   M(s \tilde\alpha, \pi_v, w_0).
\endmultline\tag1.4
$$
To proceed, we need the following assumption, originally called Assumption A
in [Ki1], and later Assumption 2.1 in [GS].

Given a place $v$, the map $f^0_v (g) \mapsto \exp \la s \tilde\alpha, H_{P_v} (g)\ra f^0_v (g),\ g \in G_v$, defines a bijection from the space of $I(0, \pi_v)$
 onto $I(s,\pi_v)$. Set $$f_v (g)= \exp \la s \tilde\alpha,
H_{P_v} (g)\ra f^0_v (g).$$

\vglue5pt
\nonumproclaim{Assumption 1.1} Fix a place $v$ and assume ${\rm Re}(s)\geq 1/2${\rm .}
Then there exists a function $f^0_v\in V (0, \pi_v)$ such that $N(s\tilde\alpha, \pi_v,w_0) f_v (g)$ 
is holomorphic and nonzero at $s$ for some $g\in G_v${\rm .}
\endproclaim

The usual duality arguments show that if Assumption 1.1 holds for $\pi$, then
it holds for $\tilde\pi$.

\demo{{R}emark}
The proof of the assumption (usually called Assumption A) in many cases is the subject matter of a work in preparation by Kim [Ki5].
As outlined there, besides the Standard Module Conjecture (cf.~[CS], [V]), one needs the following conjecture (Conjecture 7.1 of [Sh1]).
\enddemo

\nonumproclaim{Conjecture 1.2} Assume $\pi_v$ is tempered{\rm ,} then each
$L(s, \pi_v, r_{i,v})$ is holomorphic for ${\rm Re}(s) > 0${\rm .}
\endproclaim

We shall prove both the assumption and the conjecture in any of the cases which
we need in this paper. We refer to [CS],  [Ki5], and [As] for the
progress made on Conjecture 1.2.

We conclude this section by recording the following well-known
(cf.~[Sh2], [La6]) equation for the sake of completeness.

$$
M(s,\pi)f= \prod^m_{i=1} \varepsilon (is, \pi,\tilde r_i)^{-1}
L(is, \pi, \tilde r_i) L(1 + is, \pi, \tilde r_i)^{-1}
\otimes_v N(s \tilde\alpha, \pi_v, w_0) f_v,\tag1.5
$$
where $f=\hskip-6pt  \mathbold{\otimes}_v f_v \in I(s, \pi)$.
\pagebreak

\section{A general result}

In this section we will generalize an idea from [Ki1, Prop.\ 2.1] and 
[La1, Lemma 7.5]  to establish the holomorphy of each $L$-function $L(s, \pi, r_i)$ under highly ramified twists. This will be the global
analogue of [Sh6]. We shall continue to make our Assumption 1.1 for $(\bG, \bM, \pi)$ as well as each $(\bG_i, \bM_i, \pi')$ and any other cusp form
(such as
$\pi_\chi$; see below)  that appears (cf.~[GS]).

\vglue6pt
\nonumproclaim{Proposition  2.1} {\rm a)}\ Given a pair $(\bG,\bM)${\rm ,}
 where $\bG$ is a quasisplit connected reductive algebraic group over a number field $F$ and $\bM$ is an $F$\/{\rm -}\/rational 
maximal Levi subgroup of
$\bG${\rm ,} there exists a rational character $\xi$ of $\bM${\rm ,} i.e.{\rm ,}
 $\xi\in X(\bM)_F${\rm ,} with following property\/{\rm :}\/ Let $S$ be a nonempty
finite set of finite places of $F${\rm .}
 For every globally generic cuspidal representation $\pi$ of $M=\bM(\bAA_F)${\rm ,} there exist nonnegative integers
$f_v,\ v\in S${\rm ,}
 such that for every gr{\rm \"{\it o}}ssencharacter $\chi=\hskip-6pt  \mathbold{\otimes}_v\chi_v$ of $F$ for which conductor of $\chi_v,\ v\in S${\rm ,} is
larger than or equal to $f_v${\rm ,} every
$L$\/{\rm -}\/function $L(s,\pi_\chi,r_i),\ i=1,\ldots,m${\rm ,} is entire{\rm ,} where $\pi_\chi=\pi\otimes (\chi\cdot \xi)${\rm .}

{\rm b)}\ The integers $f_v$ depend only on the conductors of the local central characters of $\pi_v$ for all $v\in S${\rm .}
\endproclaim

\demo{{P}roof} We first specify $\xi$.
Let
$$
\xi(m)=\det ({\rm Ad}(m)|\frak n),
$$
$m\in\bM$, where $\frak n$ is the Lie algebra of $\bN$, as in [Sh6]. 
 \enddemo

Clearly $\xi\in X(\bM)_F$ and therefore $\xi\in X(\bM)_{F_v}$. Consequently, given a gr\"ossencharacter $\chi=\bigotimes_v\chi_v$ of $F$, each $\chi_v\cdot\xi$ becomes a character of $M_v$. Later on in the proof we will replace $\xi$ with a power of $\xi$, if need be.

Next, let $\tilde w_0$ be the longest element in $W(\bA_0)$ modulo that of the Weyl group of $\bA_0$ in $\bM$. Then $\tilde w_0$ sends $\alpha$ to a negative root, while $\tilde w_0 (\theta)\subset \Delta$.

We now apply Proposition 2.1 of [Ki1] (cf.~Lemma 7.5 of [La1]) to
equation (1.5) for $\pi_\chi$.
It states that if $\bP$ is not self-conjugate, or if $\bP$ is self-conjugate but $w_0(\pi_\chi)\not\simeq \pi_\chi$, then the constant term $M(s\tilde\alpha,\pi_\chi)$ of the corresponding Eisenstein series is holomorphic for ${\rm Re}(s)\geq 0$.
Then, together with Assumption 2.1 for $\pi_\chi$, this implies:

\nonumproclaim{Lemma 2.2} Let $\chi$ be a gr{\rm \"{\it o}}ssencharacter of $F${\rm .}
Suppose $w_0(\pi_\chi)\not\simeq \pi_\chi$ which is in particular valid if $\bP$ is not self\/{\rm -}\/conjugate{\rm .} Then
$$
\prod^m_{i=1} L(is,\pi_\chi,r_i)/L(1+is,\pi_\chi,r_i)
$$
is holomorphic for ${\rm Re}(s)\geq 1/2${\rm .}
 Here $\pi_\chi=\pi\otimes (\chi\cdot\xi)$ and $w_0$ is a representative for $\tilde w_0$ in $\bG(F)\cap K$ as in Section {\rm 1.}
\endproclaim

It is clear that from now on we may assume $\bP$ is self-conjugate and
therefore $w_0 (\bA)=\bA$. We must first show that we can find
$f_v, v \in S$, as demanded in the proposition. By Lemma 2.2, it would be enough to find a $v \in S$ and a positive integer $f_v$ for which $w_0 (\pi_v \otimes (\chi_v \cdot \xi))\not= \pi_v \otimes (\chi_v \cdot \xi)$, if\ Cond$(\chi_v) \geq f_v$. Let $\omega_v$ be the central character of $\pi_v$. Since
the center of $M_v$ contains $A_v$, it is enough to show that for a highly ramified character $\chi_v, w_0 (\omega_v (\chi_v \cdot \xi))
\neq\omega_v (\chi_v\cdot\xi)$, upon restriction to $A_v$.

Let
$$
\bA^1 = \left\{a\in\bA|\tilde w_0(a)=a^{-1}\right\}.
$$
It is a connected subgroup of $\bA$. The rational morphism $\alpha\colon \bA^1\to \bG_m$ is surjective and consequently $\alpha (A^1_v)$ is open in $F^*_v$. Let $\omega_{v,1} = \omega_v|A^1_v$. By
the conductor of $\omega_{v,1}$ we mean the smallest nonnegative integer $n_v$ such that $1 + \sP_v^{n_v}$ is contained in $\alpha (\text{Ker}
(\omega_{v,1}))$.

Next, observe that, by Lemma 2 of [Sh6], $\xi (A^1_v)$ is also open in
$F^*_v$.
Choose $\ell_v \in\bZ^+$ such that $1 + \sP^{\ell_v}_v
\subset \xi (A^1_v)$.
Now take $\chi_v$ so that the conductor of $\chi^2_v$
is larger than $\ell_v$ and that of $\omega^{-2}_{v,1}$.
This guarantees that $\chi^2_v \cdot \xi \neq \omega^{-2}_v$ upon restriction to $A^1_v$, implying $w_0 ((\chi_v \cdot \xi) \omega_v) \neq (\chi_v \cdot \xi) \omega_v$ on $A_v$ as desired. We summarize this as:

\nonumproclaim{Lemma 2.3} Suppose $\chi$ is a gr{\rm \"{\it o}}ssencharacter of $F$ which is appropriately highly ramified at one place $v\in S${\rm .}
 Then
$$
\prod^m_{i=1} L(is,\pi_\chi,r_i)/L(1+is,\pi_\chi,r_i)
$$
is holomorphic for ${\rm Re}(s)\geq 1/2${\rm .}
\endproclaim

Next we need to show that under the assumption of Proposition 2.1, each
$L(s, \pi_\chi, r_i)$ is holomorphic for ${\rm Re}(s)\geq 1/2$. The statement of the proposition then immediately 
follows from functional equation (1.1).
To do this we need to use the general induction from [Sh1, 2] (also see
Proposition 3.1 of [GS]). Simply said:\ To use the fact that for each
$i$, there exists a triple $(\bG_i,\bM_i,\pi')$ such that $$L(s,\pi,r_i)=L(s,\pi', r'_1)\hbox{ and }
\varepsilon(s, \pi,r_i)= \varepsilon (s, \pi', r'_1),
$$
where
$r'\!=\!\bigoplus\limits^{m'}_{j=1} r'_j$ is the corresponding decomposition for $(\bG_i, \bM_i)$, and $m' < m$. This is necessary since we next need
to show inductively  that each $L(s, \pi_\chi, r_i)$, $i=2, \ldots, m$, is nonzero for ${\rm Re}(s) \geq 1$, if $\chi$ is highly ramified.

For our purposes it is more convenient to choose $\bG_i$ such that $\bM_i = \bM$ and $\pi'= \pi$. We will therefore use 
Arthur's version of our induction which makes use of endoscopic groups ([A1, Prop.~5]).

In Arthur's version, groups $\bG_i$ can be chosen to be among endoscopic groups of $\bG$ which share $\bM$, and where $\pi'= \pi$. Although we do
not need this in our proof, it is worth observing that:\ {\it Let $\bA_{\bG}$ and $\bA_{\bG_i}$ be connected components of centers of $\bG$ and
$\bG_i${\rm ,} respectively{\rm ,}
 and let $\tilde\chi$ be a character of $M=\bM(\bAA_F)${\rm .} Then $\tilde\chi|(\bA_{\bG} (\bAB_F)=\tilde\chi | \bA_{\bG_i}
(\bAB_F)$ and therefore one restriction is trivial if and only if the other one is}. In our setting $\chi\cdot\xi| \bA_{\bG_i} (\bAB_F)$ is trivial since
$\chi\cdot\xi|
\bA_\bG (\bAB_F)$ is.

Next, let $\xi_i$ be the corresponding rational character defined for
$(\bG_i, \bM)$; i.e.
$$
\xi_i (m) = \det ({\rm Ad}(m)| \frak n_i).
$$
Since $X (\bA^1) \otimes_{\bZ} \bQ$ is equal to the $\bQ$-span of $\xi$, it is clear that there exist positive integers $n_i, i=1, 2, \ldots, m$, such that
$$
\xi^{n_i}_i = \xi^{n_1}.
$$
Our earlier arguments are still valid if $\xi$ is changed to $\xi^{n_1}$.
Then $\chi \cdot \xi^{n_1} = \chi \cdot \xi^{n_i}_i$ for each $i=1, \ldots,
m$, and if at a place $v \in S, \chi_v$ is ramified enough, then
$w_{i,0} (\pi_\chi) \not\simeq \pi_\chi$, where $w_{i,0}$ is the
corresponding longest element for each $i=1,2,\ldots,m$ and $w_{1,0}= w_0$.

We can now again appeal to Proposition 2.1 of [Ki1] (Lemma 7.5 of [La1]) to conclude that the corresponding Eisenstein series and consequently its Fourier coefficients are all holomorphic for ${\rm Re}(s) \geq 0$, if $\chi_v$
is highly ramified at a place $v\in S$. It then immediately follows, say for example from Proposition 5.2 of [GS], that
$$
\prod^m_{i=1} L(1 + is, \pi_\chi, r_i)\tag2.1
$$
is nonzero for ${\rm Re}(s)\geq 0$. (One needs to notice that local $L$-functions are never zero.)

Our induction hypothesis is now that:\ {\it If $\chi$ is highly ramified at a place $v\in S${\rm ,}
 then each $L(s,\pi_\chi,r_i),\ i=2,\ldots,m${\rm ,} is holomorphic for ${\rm Re}(s)\geq 1/2$ and nonzero for ${\rm Re}(s)\geq 1$}.

We now assume $\chi$ is appropriately ramified such that $w_{i,0} (\pi_\chi)
\not\simeq \pi_\chi$ for all $i=1,2,\ldots,m$. Applying our induction hypothesis for $L(s,\pi_\chi, r_i), i=2,\ldots,m$, first to (2.1) to conclude that $L(s,\pi_\chi, r_1)$ is nonzero for ${\rm Re}(s) \geq 1$, and next,
applying this  to Lemma 2.3, to get the holomorphy of $L(s,\pi_\chi, r_1)$ again for ${\rm Re}(s)\geq 1/2$, complete  the induction.  Proposition 2.1
is now proved.
\phantom{almost time}\hfill\qed

\proclaimtitle{of the proof}
\specialnumber{2.4}
\proclaim{{C}orollary} Suppose $\chi$ is highly ramified for some $v\in S${\rm .} Then every $L$-function $L(s,\pi_\chi,r_i)$
is nonzero for ${\rm Re}(s)\geq 1,\ 1\leq i\leq m${\rm .}
\endproclaim

\nonumproclaim{Corollary 2.5} Suppose $m=1${\rm ,} but $\pi$ is any cuspidal
representation of $M=\bM(\bAA_F)${\rm ,} i{\rm .}e{\rm .,} not necessarily globally
generic{\rm .} Then $L(s,\pi_\chi,r_1)$ is analytic for ${\rm Re}(s)\geq 1/2${\rm ,} if
$\chi$ is unitary and highly ramified{\rm .} Moreover{\rm ,} if it satisfies a
standard functional equation{\rm ,} then it is entire{\rm .} Similar statements are
true if $m=2$ and $r_2$ is one\/{\rm -}\/dimensional{\rm .}
\endproclaim

\vglue-16pt
\section{Functorial products for ${\rm GL}_2\times {\rm GL}_3$;  The weak lift}

\vglue-8pt

In this section we shall start proving our main results. We will first establish Theorem 3.8 from which other results follow. To establish the lifts we shall use an appropriate version of converse theorems of Cogdell and Piatetski-Shapiro 
[CP-S1], and it is in fact quite surprising that analytic properties of appropriate $L$-functions that are needed are mainly established using our
method [Sh1, 2, 4, 7], [Ki1], [GS] rather than that of Rankin-Selberg, where converse theorems have their roots.

We shall start by considering four special cases of our general machinery outlined in Section 1.

Let $\pi_1$ and $\pi_2$ be two cuspidal representations of ${\rm GL}_2(\bAA_F)$ and ${\rm GL}_3(\bAA_F)$, respectively. Here are the four cases. In each case we give a triple $(\bG,\bM,\pi)$ 
as in Section 1 and in what follows $\sigma$ always denotes a cuspidal representation of ${\rm GL}_n(\bAA_F),\ n=1,2,3$, and 4.

\vglue4pt
{\it Case} 1 ({\it Case} iii   {\it of}  [La6]). Here $\bG={\rm GL}_5,\ \bM={\rm GL}_2\times {\rm GL}_3,\
\pi=(\pi_1\otimes\sigma)\otimes\tilde\pi_2,\ n=1$ and therefore $\sigma$ is a gr\"ossencharacter of $F$. The index $m=1$. The theory of
$L$-functions of Section 1 in this case is that of Rankin-Selberg $L$-functions $L(s,(\pi_1\otimes\sigma)\times\pi_2))$ which is very well
developed. In particular the necessary analytic properties of these $L$-functions are all well-known (cf.~[JS1], [JP-SS], [Sh4, 5], [MW2]).
\vglue4pt

 {\it Case} 2 ($D_5-2$ {\it of}  [Sh2]). In this case $\bG$ is the simply connected type of $D_5$ and the pair $(\bG,\bM)$ is as in Case $D_5-2$
of [Sh2], i.e., $\bG=\Spin_{10}$ and the derived group of $\bM$ is ${\rm SL}_3\times {\rm SL}_2\times {\rm SL}_2$. The index $m=2$.
\vglue4pt
{\it Case} 3 ($E_6-1$ {\it of}  [Sh2]). Here the group $\bG$ is the simply connected $E_6$ and the derived group of $\bM$ is ${\rm
SL}_3\times {\rm SL}_2\times {\rm SL}_3$. The index $m=3$.
\vglue4pt
{\it Case} 4 ($E_7-1$ {\it of} [Sh2]). Finally we take $\bG$ to be the simply connected $E_7$ and take the derived group of $\bM$ equal to
${\rm SL}_3\times {\rm SL}_2\times {\rm SL}_4$. Here $m=4$.
\vglue4pt

We look at Case 4 in detail. The other cases are similar.
We will use Bourbaki's notation.
Let $\theta=\Delta-\{\alpha_4\}$.
Let $\bold P=\bold P_{\theta}=\bold M\bold N$ and let $\bold A$ be the connected component of the center of $\bold M$.
Then
$$
\align
\bold A&= \left(\bigcap_{\alpha\in\theta} {\rm ker}\,\alpha\right)^0\cr
&= \{a(t)=H_{\alpha_1}(t^4)H_{\alpha_3}(t^8)H_{\alpha_2}(t^6)
H_{\alpha_4}(t^{12})H_{\alpha_5}(t^9)H_{\alpha_6}(t^6)H_{\alpha_7}(t^3): t\in\overline F^*\}.
\endalign
$$ 
Since $\bold G$ is simply connected, the derived group $\bold M_D$ of $\bold M$ is simply connected, and consequently
$$\bold M_D={\rm SL}_2\times {\rm SL}_3\times {\rm SL}_4.
$$
Moreover
$$\bold A\cap \bold M_D=\{H_{\alpha_1}(t^4)H_{\alpha_3}(t^8)H_{\alpha_2}(t^6)
H_{\alpha_4}(t^{12})H_{\alpha_5}(t^9)H_{\alpha_6}(t^6)H_{\alpha_7}(t^3):\,
 t^{12}=1\}.
$$
If we identify $\bold A$ with ${\rm GL}_1$, then
$$
\bold M=({\rm GL}_1\times {\rm SL}_2\times {\rm SL}_3\times {\rm SL}_4)/(\bold A\cap \bold M_D).
$$

We define a map $\bar f: \bold A\times \bold M_D\longrightarrow {\rm GL}_1\times {\rm GL}_1\times {\rm GL}_1\times {\rm SL}_2\times {\rm SL}_3\times {\rm SL}_4$ by
$$\bar f: (a(t),x,y,z)\longmapsto (t^6,t^4,t^3,x,y,z).
$$
It induces a map $f: \bold M\longrightarrow {\rm GL}_2\times {\rm GL}_3\times {\rm GL}_4$
which is in fact an injection.
We need to determine $f(H_{\alpha_4}(t))$, $t\in\overline F^*$.
Write $H_{\alpha_4}(t)=axyz,\ a\in\bold A,\ x\in {\rm SL}_2,\ y\in {\rm SL}_3$, and $z\in {\rm SL}_4$.
Then, a glance at
$$
H_{\alpha_4}(t^{12})=a(t) H_{\alpha_2} (t^{-6})H_{\alpha_1}(t^{-4})
H_{\alpha_3} (t^{-8}) H_{\alpha_5} (t^{-9}) H_{\alpha_6} (t^{-6})
H_{\alpha_7}(t^{-3})
$$
shows that for a fixed $12^{\rm th}$ root $t^{1/12}$ of $t$, whose choice is irrelevant upon going to $f,\ x=H_{\alpha_2}(t^{-1/2})$,
$y=H_{\alpha_1} (t^{-1/3}) H_{\alpha_3}(t^{-2/3})$, and
$$
z=H_{\alpha_5} (t^{-3/4}) H_{\alpha_6} (t^{-1/2}) H_{\alpha_7} (t^{-1/4}),
$$
where $t^{1/2}=(t^{1/12})^6,\ t^{1/3}=(t^{1/12})^4$, and $t^{1/4}=(t^{1/12})^3$.
Moreover, $a=a(t^{1/12})$.

Using this, we see easily  that
$$
f(H_{\alpha_4}(t))=({\rm diag}(1,t),\, {\rm diag}(1,1,t),\, {\rm diag}(1,1,1,t)).\tag3.1
$$

Let $\pi_i,\ i=1,2$, be cuspidal representations of ${\rm GL}_{1+i}(\bAB_F)$ with central characters $\omega_{\pi_i},\ i=1,2$, respectively.
Let $\sigma$ be a cuspidal representation of ${\rm GL}_4(\bAB_F)$ whose central character is denoted by $\omega_\sigma$.

Since $f$ is $F$-rational, it induces an injection
$$
f\colon \bM(\bAB_F)\longrightarrow {\rm GL}_2(\bAB_F)\times {\rm GL}_3(\bAB_F)\times {\rm GL}_4(\bAB_F).
$$
Moreover $\bM(\bAB_F) (\bAB_F^*)^2$ is co-compact in ${\rm GL}_2(\bAB_F)\times {\rm GL}_3(\bAB_F)\times {\rm GL}_4(\bAB_F)$, where $(\bAB_F^*)^2$ is embedded as 
a center  of, say, the first two factors.
Consequently $\pi_1\otimes \pi_2\otimes\sigma|f(M),\ M=\bM(\bAB_F)$, decomposes to a direct sum of irreducible cuspidal representations of $M$.
Let $\pi$ be any irreducible constituent of this direct sum.
Then the central character of $\pi$ is
$\omega_\pi=\omega_{\pi_1}^6 \omega_{\pi_2}^4 \omega_\sigma^3$.
As we shall see, choice of $\pi$ is irrelevant.
Write $\pi=\hskip-6pt  \mathbold{\otimes}_v\pi_v$.
(The fact that $f$ is an injection is not important.
All that we need is that there is an $F$-rational map which is the identity on $\bM_D$.)

Now suppose each $\pi_{iv}$ is an unramified representation, given by
$$\pi_{1v}=\pi(\eta_1,\eta_2),\quad \pi_{2v}=\pi(\nu_1,\nu_2,\nu_3),\quad \sigma_{v}=\pi(\mu_1,\mu_2,\mu_3,\mu_4).
$$
Then $\pi_v$ is an unramified representation of $\bold M(F_v)$,
induced from the character $\chi$ of the torus of $\bold M(F_v)$ whose image under $f$ lies in the product of diagonal subgroups of corresponding 
${\rm GL}_i(F_v)$'s.
Since $f$ is the identity on $M_D$,  
$$
\alignat3
\chi\circ H_{\alpha_7}(t)&=\mu_1\mu_2^{-1}(t),\quad\chi\circ H_{\alpha_6}(t)&=\mu_2\mu_3^{-1}(t),\quad \chi\circ
H_{\alpha_5}(t)&=\mu_3\mu_4^{-1}(t), \\
\chi\circ H_{\alpha_3}(t)&=\nu_2\nu_3^{-1}(t),\quad \chi\circ H_{\alpha_1}(t)&=\nu_1\nu_2^{-1}(t),\quad\chi\circ
H_{\alpha_2}(t)&=\eta_1\eta_2^{-1}(t).
\endalignat
$$ 
Moreover, (3.1) implies
$$
\chi\circ H_{\alpha_4}(t)=\mu_4\nu_3\eta_2 (t).
$$
We conclude that $m=4$, and
$$
\align
L(s,\pi_v,r_1)&=L(s,\pi_{1v}\times\pi_{2v}\times\sigma_v),\\
L(s,\pi_v,r_2)&=L(s,\tilde\pi_{2v}\otimes\sigma_v, (\rho_3\otimes\omega^2_{\pi_{1v}}\omega_{\pi_{2v}})\otimes\wedge^2\rho_4),\\
L(s,\pi_v,r_3)&=L(s,(\pi_{1v}\otimes\omega_{\pi_{1v}}\omega_{\pi_{2v}}\omega_{\sigma_v})\times\tilde\sigma_v),\\
L(s,\pi_v,r_4)&=L(s,\pi_{2v}\otimes\omega_{\pi_{1v}}\omega_{\pi_{2v}}\omega_{\sigma_v}).
\endalign
$$
(By Proposition 2.1 we do not need to know the precise form of the last three $L$-functions.)

In Case 2, we do the same, namely, we construct a map $f: \bold M\longrightarrow {\rm GL}_2\times {\rm GL}_3\times {\rm GL}_2$.
Let $\pi_i,\ i=1,2$, be cuspidal representations of ${\rm GL}_{1+i}(\Bbb A_F)$ with central characters $\omega_{\pi_i}$, $i=1,2$, resp. Let $\sigma$ be a cuspidal representation of ${\rm GL}_2(\Bbb A_F)$ with central character $\omega_{\sigma}$.
Finally, let $\pi$ be an irreducible cuspidal representation of $\bold M(\Bbb A_F)$, induced by the map $f$ from $\pi_1,\pi_2,\sigma$ as before.
Then the central character of $\pi$ is $\omega_\pi=\omega_{\pi_1}^3\omega_{\pi_2}^2 \omega_\sigma^3$ and for an unramified place $v$ 
$$
\align
L(s,\pi_v,r_1)&=L(s,\pi_{1v}\times\pi_{2v}\times\sigma_v),\\
L(s,\pi_v,r_2)&=L(s,\tilde\pi_{2v}\otimes\omega_{\pi_{1v}}\omega_{\pi_{2v}}\omega_{\sigma_v}).
\endalign
$$

In Case 3, we again construct a map $f: \bold M\longrightarrow {\rm GL}_2\times {\rm GL}_3\times {\rm GL}_3$ and proceed as before.
Then the central character of $\pi$ is $\omega_\pi=\omega_{\pi_1}^3 \omega_{\pi_2}^2 \omega_\sigma^2$ and for an unramified place $v$
$$
\align
L(s,\pi_v,r_1)&=L(s,\pi_{1v}\times\pi_{2v}\times\sigma_v),\\
L(s,\pi_v,r_2)&=L(s,(\tilde\pi_{2v}\otimes\omega_v)\times\tilde\sigma_v),\\
L(s,\pi_v,r_3)&=L(s,\pi_{1v}\otimes\omega_v),
\endalign
$$
where $\omega_v=\omega_{\pi_{1v}}\omega_{\pi_{2v}}\omega_{\sigma_v}$.

In ramified places, we take $L(s,\pi_v,r_i)$ to be the one defined in [Sh1] for each of these cases. Observe that in particular, if $v=\infty$, then $L(s,\pi_v,r_i)$ is the corresponding Artin $L$-function
 (cf. [La4,  5], [Sh7]) in each case.

\demo{{R}emark} One can also use similitude groups to get these $L$-functions since local coefficients depend only on derived groups, while
the Levi subgroups are less complicated.
For example, the corresponding Levi in $GE_7=({\rm GL}_1\times E_7^{sc})/Z(E_7^{sc})$ is isomorphic to ${\rm GL}_3\times M_0$, where $M_0$ is the standard Levi subgroup ${\rm SL}_6\cap ({\rm GL}_2\times {\rm GL}_4)$ of ${\rm SL}_6$.
For $GE_6$, $M={\rm GL}_2\times M_0$, where $M_0={\rm SL}_6\cap ({\rm GL}_3\times {\rm GL}_3)$, while for ${\rm GSpin}(10),\ M={\rm GL}_3\times {\rm GSpin}(4)$ (cf.~[As]).
Note that in each case there is an $F$-rational injection $f$ from $\bold M$ to ${\rm GL}_2\times {\rm GL}_3\times {\rm GL}_k,\ k=2,3,4$.
\enddemo

Next, fix a nontrivial additive character $\psi=\hskip-6pt  \mathbold{\otimes}_v\psi_v$ of $F\backslash \bAA_F$ and define local root numbers
$\var(s,\pi_v,r_1,\psi_v)$ again as in [Sh1]. Similar comment applies when $v=\infty$.

To proceed, we must dispose of Assumption 1.1 in our four cases, as well as in the inductive cases $(\bG_i,\bM_i,\pi')$ attached to them.
We will prove it for components of arbitrary $\pi$ and all the corresponding $\pi'$.
With notation as in Section 1, let $N(s\tilde\alpha,\pi_v,w_0)$ denote the normalized local intertwining operator (1.4) in each of our cases. In the next section we will prove:

\nonumproclaim{Proposition  3.1} The normalized local intertwining operators 
$$N(s\tilde\alpha,\pi_v,\omega_0)$$ are holomorphic and nonzero for ${\rm Re}(s)\geq 1/2$ and for all $v${\rm .}
\endproclaim

For each $v$, let $\pi_{1v}\boxtimes\pi_{2v}$ be the irreducible admissible representation of ${\rm GL}_6(F_v)$ attached
 to $\pi_{1v}\otimes\pi_{2v}$ through the local Langlands correspondence. More precisely,  if $\delta_{1v}$ and $\delta_{2v}$ are representations of
the Deligne-Weil group, parametrizing $\pi_{1v}$ and $\pi_{2v}$ through the local Langlands correspondence for ${\rm GL}_2(F_v)$ and ${\rm
GL}_3(F_v)$, respectively, then let $\pi_{1v}\boxtimes\pi_{2v}$ be the representation of ${\rm GL}_{6}(F_v)$ attached to
$\delta_{1v}\otimes\delta_{2v}$ ([HT], [He], [La4]).

Set
$$\pi_1\boxtimes\pi_2=\hskip-6pt  \mathbold{\otimes}_v (\pi_{1v}\boxtimes\pi_{2v}).
$$
It is an irreducible admissible representation of ${\rm GL}_{6}(\Bbb A_F)$. Langlands' functoriality in our case is equivalent to the assertion that $\pi_1\boxtimes\pi_2$ is an automorphic representation.

To proceed we need to state the converse theorem that applies to our situation.

\proclaimtitle{Theorem 2 of [CP-S1]}
\specialnumber{3.2} 
\proclaim{Theorem} Suppose $\Pi=\hskip-6pt  \mathbold{\otimes}_v \Pi_v$ is an irreducible admissible representation of 
${\rm GL}_m(\bAA_F)$ whose central character is a gr{\rm \"{\it o}}ssencharacter{\rm .}
Let $S$ be a finite set of finite places of $F$ and let $\Cal T^S(n)$ be the set of cuspidal representations of 
${\rm GL}_n(\bAA_F)$ that are unramified at all places $v\in S${\rm .}
Suppose for each $n\leq m-2$ and every cuspidal representation
$\sigma\in\Cal T^S (n), L(s,\Pi\times\sigma)$ is {\rm ``}\/nice\/{\rm ''} in the sense that it satisfies the following properties\/{\rm :}
\vglue2pt
\itemitem{\rm a)}The $L$-function $L(s,\Pi\times \sigma)$ is entire{\rm ,}
\vglue2pt
\itemitem{\rm b)}it is bounded in every vertical strip of finite width{\rm ,} and
\vglue2pt
\itemitem{\rm c)}it satisfies a standard functional equation of type {\rm (1.1).}
\vglue2pt
\noindent Then there exists an automorphic representation $\Pi'$ of ${\rm GL}_m(\bAA_F)$ such that $\Pi_v\simeq \Pi'_v$ for all $v\not\in S$ 
and in
particular for every archimedean place of $F${\rm .} We will then say that   $\Pi$ is quasi\/{\rm -}\/automorphic with respect to $S${\rm .}
\endproclaim

We apply the converse theorem to $\pi_1\boxtimes\pi_2=\hskip-6pt  \mathbold{\otimes}_v (\pi_{1v}\boxtimes\pi_{2v})$. For that, we need
 to consider Rankin triple product $L$-functions  found in our cases\break 1--4. To show that they can be made ``nice,'' we will twist our
representation $\pi$ by an appropriately ramified gr\"ossencharacter of $F$ to which we can apply Proposition 2.1 to verify condition a) of Theorem
3.2. Condition b) is proved in full generality in [GS]. Condition c) is delicate and requires multiplicativity of local factors (Theorem 3.5.3 of [Sh1]).
We need to prove
$$
L(s,\pi_{1v}\times \pi_{2v}\times\sigma_v)=L(s,(\pi_{1v}\boxtimes\pi_{2v})\times\sigma_v)\tag3.2.1
$$
and
$$
\var(s,\pi_{1v}\times \pi_{2v}\times\sigma_v,\psi_v)=\var(s,(\pi_{1v}\boxtimes\pi_{2v})\times\sigma_v,\psi_v),\tag3.2.2
$$
for all irreducible admissible generic representations $\sigma_v$ of ${\rm GL}_n(F_v),\ 1\leq n\leq 4$.

These equalities are not obvious at all since the two factors on the left and right are defined using completely different techniques. The factors on the left are those defined by the triple $L$-functions of our four cases (to be recalled in Section 5), while the ones on the right are those of Rankin-Selberg for ${\rm GL}_6(\bAA_F)\times {\rm GL}_n (\bAA_F),\ 1\leq n\leq 4$
 [JP-SS], [Sh5, 7].  By the local Langlands correspondence, they are Artin factors.

In this section we will apply the converse theorem to $\pi_1\boxtimes\pi_2$
 where $S$ is a finite set of places of $F$ which contain
all the finite places $v$ where at least one of $\pi_{iv}$, $i=1,2$, is
ramified. This makes the application of the converse theorem simpler because
if $\sigma\in \Cal T^S(n)$, one of $\sigma_v,\pi_{1v}, \pi_{2v}$, is in the principal series for $v<\infty$.

In Section 5, we shall prove that the equalities (3.2.1) and (3.2.2) hold
for all $v$.
Machinery of converse theorems then implies that $\pi_1\boxtimes\pi_2$ is automorphic.
(See the proof immediately after the proof of Proposition 5.8 as well as the paragraph before Proposition 5.4.)
We start with the following definition:

\vglue2pt {\it Definition {\rm 3.3}}. Let $\pi_1$ ($\pi_2$, resp.) be irreducible cuspidal representations of ${\rm GL}_2(\bAA_F)$ $({\rm GL}_3
(\bAA_F)$, resp.). {An automorphic representation} $\Pi=\hskip-6pt  \mathbold{\otimes}_v \Pi_v$ of ${\rm GL}_6(\bAA_F)$ is a {\it strong lift} or
{\it transfer} of
$\pi_1\otimes\pi_2$, if for every $v$, $\Pi_v$ is a local lift or transfer of $\pi_{1v}\otimes\pi_{2v}$, in the sense that
$$ \align
L(s,\pi_{1v}\times \pi_{2v}\times\sigma_v)&= L(s,\Pi_v\times\sigma_v)\tag3.3.1
\cr
\noalign{\noindent  and}
\var(s,\pi_{1v}\times \pi_{2v}\times\sigma_v,\psi_v)&=\var(s,\Pi_v\times\sigma_v,\psi_v),\tag3.3.2 \endalign
$$
for all irreducible admissible generic representations $\sigma_v$ of ${\rm GL}_n(F_v),\ 1\leq n\leq 4$. If these equalities hold for all $v\not\in S$,
then   $\Pi$ is  {\it a weak lift} of $\pi_1\otimes\pi_2$ with respect to $S$.
\vglue2pt

Let $S$ be a finite set of finite places of $F$ which contain all the places
$v$ where at least one of $\pi_{iv}$, $i=1,2$, is ramified.
Fix $v_0\in S$. Take a gr\"ossencharacter $\chi=\hskip-6pt  \mathbold{\otimes}_v \chi_v$.
We will assume $\chi_{v_0}$ is appropriately highly ramified so that Proposition 2.1 can be applied to $L$-functions $L(s,\pi_\chi,r_1)$ 
in each of our four cases.

Replacing $\chi$ with an appropriate integral power of $\chi$ if necessary, one can find gr\"ossencharacters $\chi_n,n=1,2,3,4$, such that
$$
L(s,(\pi_1\otimes\chi)\times\pi_2\times\sigma)=L(s,\pi_{\chi_n},r_1)
$$
for $n=1,2,3,4$, representing each of our four cases, or equivalently as
$\sigma\in\Cal T^S (n),\ n=1,2,3,4$.

As $\sigma$ ranges in $\Cal T^S (n)$, the conductor of the central character of $\pi$ will not change and therefore in view of Proposition 3.1, Part b) of Proposition 2.1 applies, implying:

\nonumproclaim{Proposition  3.4} 
There exists a positive integer $f_0$ such that for every gr{\rm \"{\it o}}ssencharacter $\chi=\hskip-6pt  \mathbold{\otimes}_v \chi_v$ 
with {\rm Cond} $(\chi_{v_0})\geq f_0,\
L(s,(\pi_1\otimes\chi)\times \pi_2\times\sigma)$ is entire for every $\sigma\in \Cal T^S (n),\ n=1,2,3,4${\rm .}
\endproclaim

Next, taking into account our Proposition 3.1 and applying Theorem 4.1 of [GS] to our four cases, we have:

\nonumproclaim{Proposition  3.5} Fix $f_0$ as in Proposition {\rm 3.4.}
Then each $$L(s,(\pi_1\otimes\chi)\times\pi_2\times\sigma)$$ is bounded in every vertical strip of finite width{\rm .}
\endproclaim

Next we show:

\nonumproclaim{Proposition  3.6}\hglue-4pt Notation being as in Proposition {\rm 3.5,} let $\sigma\!\in\! \Cal T^S(n)\otimes\chi${\rm .} Then
for each $v${\rm ,}
$$\align
\gamma(s,\pi_{1v}\times \pi_{2v}\times\sigma_v,\psi_v)&=\gamma(s,(\pi_{1v}\boxtimes\pi_{2v})\times\sigma_v,\psi_v),\\
L(s,\pi_{1v}\times \pi_{2v}\times\sigma_v)&=L(s, (\pi_{1v}\boxtimes\pi_{2v})\times\sigma_v).
\endalign
$$
The equality of $\varepsilon$\/{\rm -}\/factors follows from the above equalities{\rm .}
\endproclaim

\demo{Proof} As we noted before, if $\sigma\in \Cal T^S(n)\otimes\chi$, then for each $v<\infty$, 
one of $\sigma_v,\pi_{1v},\pi_{2v}$, is in the principal series. When $v=\infty$, as has been shown in [Sh7], the factors on the left are Artin
factors, and hence we have the equalities. If $v<\infty$, by multiplicativity of $\gamma$-factors and $L$-functions (Part 3 of Theorem 3.5 and
Section 7 of [Sh1]; also see the discussions at the beginning of Section 5 here), the factors on the left-hand side are a product of those for ${\rm
GL}_k\times {\rm GL}_l$. Shahidi [Sh5] has shown that in the case of ${\rm GL}_k\times {\rm GL}_l$, his factors are those of Artin. Since the same
multiplicativity holds for the factors on the right-hand side, we have the equalities.  
\enddemo

\demo{{R}emark} Multiplicativity of $L$-functions is particularly transparent in this case since the principal series representation among $\pi_{iv}$ and $\sigma_v$ at each place $v$ is in fact its own standard module (cf.~Section 7 of [Sh1] and Section 5 here).
\enddemo

\nonumproclaim{Proposition  3.7} Fix $S$ and $\chi$ as in Proposition {\rm 3.4.}
Then each\break $L$\/{\rm -}\/function $L(s,((\pi_1\boxtimes\pi_2)\otimes\chi)\times\sigma)$ is {\rm ``}\/nice\/{\rm ''} 
 as $\sigma$ runs over the sets $\Cal
T^S (n),\break n=1,2,3,4${\rm .}
\endproclaim

\demo{Proof} This follows immediately from Propositions 3.4, 3.5, 3.6, and the functional equation satisfied by 
$L(s,\pi_{1}\times \pi_{2}\times\sigma)$, proved in [Sh1].  
\enddemo

We now apply the Converse Theorem 3.2 to Proposition 3.7 to conclude:

\nonumproclaim{Theorem 3.8} Let $S$ be a finite set of finite places containing all the places $v$ for which either $\pi_{1v}$ or $\pi_{2v}$ is ramified{\rm .}
Then there exists an automorphic representation $\Pi=\hskip-6pt  \mathbold{\otimes}_v \Pi_v$ 
of ${\rm GL}_6(\Bbb A_F)$ such that $\Pi_v\simeq\pi_{1v}\boxtimes \pi_{2v}$ for $v\notin S${\rm .}
\endproclaim

\demo{Proof} Using Proposition 3.7, we only need to apply Theorem 3.3 to\break
$(\pi_1\boxtimes\pi_2)\otimes\chi$ for some gr\"ossencharacter $\chi$ which is highly ramified on $S$. Thus
 $(\pi_1\boxtimes\pi_2)\otimes\chi$ is quasi-automorphic with respect to $S$ and therefore so is $\pi_1\boxtimes\pi_2$. 
\enddemo

\demo{{R}emark {\rm 3.9}} Even if one does not have the local Langlands correspondence, one can still use the converse theorem. For each $v\in S$,
take $\Pi_v$ to be arbitrary, up to its central character, which is predetermined by central characters of $\pi_{iv},\ i=1,2$. 
We would then need to use stability of Rankin-Selberg $\gamma$-functions under highly ramified twists (Proposition 4 of [JS2]) which is a deep
result from the Rankin-Selberg method. We refer to [CP-S3] for more detail, as well as [CKPSS] for a very important application.
\enddemo

Now, let $\Pi=\hskip-6pt  \mathbold{\otimes}_v \Pi_v$ denote a {\it weak lift} of $\pi_1\otimes\pi_2$, with respect to $S$, i.e.~an automorphic representation for which $\Pi_v\simeq\pi_{1v}\boxtimes
\pi_{2v}$,  for all  $v\not\in S$. Choose real numbers $r_i$ and (unitary) cuspidal representations $\sigma_i$ of ${\rm GL}_{n_i}(\bAA_F),\break
i=1,\ldots,k$, such that
$\Pi$ is equivalent to a subquotient of
$$
I=\text{Ind }\sigma_1 |\det(\ )|^{r_1}\otimes\cdots\otimes\sigma_k|
\det(\ )|^{r_k}.
$$
Since the central character $\omega_{\Pi}=\omega_{\pi_1}^3\omega_{\pi_2}^2$ is unitary,
$$
n_1 r_1+\cdots+n_k r_k=0.
$$
Observe that $n_i > 1$ since $L_S(s,\Pi\otimes \mu)=L_S(s,(\pi_1\otimes \mu)\times\pi_2)$ is entire for every gr\"ossencharacter $\mu$.
Here $L_S$ is the partial $L$-function with respect to $S$, i.e., the product of all the local $L$-functions for $v$ outside of $S$.
We shall prove:

\nonumproclaim{Proposition  3.10} The exponents $r_1=\cdots=r_k=0$ and
$I=\Ind \ \sigma_1\otimes\cdots\otimes\sigma_k$ are irreducible.
Therefore in the notation of {\rm [La3]}
$$
\Pi=\sigma_1\boxplus\cdots\boxplus\sigma_k.\tag3.10.1
$$
Moreover{\rm ,} $\Pi$ is unique and each local component $\Pi_v$ of $\Pi$ is irreducible{\rm ,} unitary{\rm ,} and generic{\rm .}
\endproclaim

\demo{Proof} We need to use the weak Ramanujan property for a (unitary) cuspidal representation $\pi=\hskip-6pt  \mathbold{\otimes}_v\pi_v$ of ${\rm GL}_n(\bAA_F)$, where $n$ is a positive integer. Let $\pi_v$ be an unramified component of $\pi$.
Write:
$$
\pi_v=\Ind \mu_1|\ \ |^{s_1}\otimes\cdots\otimes\mu_\ell| \ \ |^{s_\ell}
\otimes\nu_1\otimes\cdots\otimes\nu_u\otimes\mu_\ell|\ \ |^{-s_\ell}
\otimes\cdots\otimes \mu_1|\ \ |^{-s_1},\tag3.10.2
$$
where $\mu_i$ and $\nu_j$ are unitary characters of $F_v^*$ and
$0 < s_\ell \leq\cdots\leq s_1 < 1/2$, by classification of irreducible
unitary generic representations of ${\rm GL}_n(F_v)$ (cf.~[Tad]).
Here we have suppressed the dependence of all the factors on $v$ for simplicity of notation.
Then $\pi$ is said to satisfy the {\it weak Ramanujan property}, if given
$\var > 0$, the set of places $v$ for which $s_1\geq\var$  is of density
zero (cf.~[CP-S2] for the original idea; also see [Ki2]).

It then follows from Ramakrishnan (Lemma 3.1 of [Ra2]) that for $n=2$ and 3
every irreducible cuspidal representation satisfies the weak Ramanujan property.
Consequently, so does every weak lift $\Pi$ of $\pi_1\otimes\pi_2$ (with the same definition for an automorphic representation).

Suppose $\tau_i=\sigma_i|\det(\ \ )|^{r_i},\ 1\leq i\leq k$.
For an unramified pair $(\pi_{1v},\pi_{2v})$, let $t_v={\rm diag}((a_{iv},
\eta_v a^{-1}_{iv},\ldots)$ denote the Hecke-Frobenius (Satake) parameter of $\sigma_{iv},\ |\eta_v|=1$, by the fact that $\sigma_{iv}$ has the form
(3.10.2). We may assume $|a_{iv}|\geq 1$.

By the equality $n_1 r_1+\cdots +n_k r_k=0$, if $r_i$ are not all zero, then there is $i$ such that $r_i > 0$.
Using the fact that $|\varpi_v|_v^{r_i} t_v$ defines the
Hecke-Frobenius parameter of $\tau_{iv}$, we have
$$
|\tilde a_{iv}|=|a_{iv}|^{-1} q_v^{-r_i}\leq q_v^{-r_i} < 1,
$$
where $\tilde a_{iv}=\eta_v a_{iv}^{-1}|\varpi_v|_v^{r_i}$.
Now take $\var=r_i$ to contradict the weak Ramanujan property for $\Pi$.
Thus $r_i=0$ for all $i$.

Note that for each $v,{\rm Ind}\, \sigma_{1v}\otimes\cdots\otimes\sigma_{kv}$ is irreducible, unitary, and generic [Tad].
Hence $\Pi_v={\rm Ind}\, \sigma_{1v}\otimes\cdots\otimes\sigma_{kv}$ for all $v$.
Thus $\Pi=\sigma_1\boxplus \cdots \boxplus \sigma_k$.

{}From the classification theorem for $GL(n)$ (cf.~[JS1]) and (3.10.1) it is clear that the weak lift is unique.
The proposition is now proved. 
\enddemo

\section{Proof of Proposition 3.1}

In [Ki5], it was shown that in a fairly general setting, the normalized local intertwining operators $N(s,\sigma_v,w_0)$ are holomorphic and nonzero for ${\rm Re}\, s\geq \frac 12$ for all $v$. For the sake of completeness, we include a proof here in our cases. First we 
need:

\nonumproclaim{Lemma 4.1} Conjecture {\rm 7.1} of {\rm [Sh1]} is true in the following cases\/{\rm :}
\vglue3pt
\item{\rm (1)} $D_n-2$ $(n\geq 4)${\rm :} i.e.{\rm ,} the local $L$\/{\rm -}\/function 
$L(s,\pi_{1v}\times\pi_{2v}\times\sigma_{v})$ is holomorphic for ${\rm Re}\, s>0${\rm ,}
where
$\pi_{1v},\pi_{2v}${\rm ,} and $\sigma_{v}$ are tempered representations of ${\rm GL}_2(F_v), {\rm GL}_{n-2}(F_v), {\rm GL}_2(F_v)${\rm ,} 
\vglue3pt\item{\rm (2)} $E_6-1$ and $E_7-1${\rm :} i.e.{\rm ,} the local $L$\/{\rm -}\/functions
$L(s,\pi_{1v}\times\pi_{2v}\times\sigma_{v})$ are holomorphic for ${\rm Re}\, s>0${\rm ,} where $\pi_{1v},\pi_{2v}${\rm ,} and $\sigma_{v}$
are tempered representations of ${\rm GL}_2(F_v), {\rm GL}_3(F_v)${\rm ,} and ${\rm GL}_n(F_v)${\rm ,}
$n=3,4${\rm ,} resp{\rm .}
\vglue3pt\item{\rm (3)} $D_6-3:$ i.e.{\rm ,} the local $L$\/{\rm -}\/function
$L(s,\pi_{2v}\otimes\sigma_{v},\rho_3\otimes\wedge^2\rho_4)$ is holomorphic for ${\rm Re}\, s>0${\rm ,} where $\pi_{2v}$ and $\sigma_{v}$
are tempered representations of ${\rm GL}_3(F_v)$ and ${\rm GL}_4(F_v)${\rm ,} resp{\rm .}\
Note that the case $D_6-3$ appears as the second
$L$-function in the case $E_7-1$ {\rm (}\/see Section {\rm 3).}

\endproclaim
 
{\it Proof}. For simplicity, we drop the subscript $v$.
For the case $D_n-2$ ($m=2$), we can use either Asgari's thesis ([As]) in
which many other results on Conjecture 7.1 of [Sh1] are proved,
or use Ramakrishnan's result [R1] on the local lift of
${\rm GL}_2(F)\times {\rm GL}_2(F)$ to ${\rm GL}_4(F)$ to reduce it to a Rankin-Selberg
$L$-function for ${\rm GL}_{n-2}(F)\times {\rm GL}_4(F)$, namely,
$$L(s,\pi_1\times\pi_2\times\sigma)=L(s,\pi_2\times
(\pi_1\boxtimes\sigma)).
$$

For the case $D_6-3$, we need to use [As]. In this case, the trouble
is when $\pi_{2v}$ is a Steinberg representation, for then factors
cancel between numerator and denominator of $\gamma$-factors.
If we use $\Spin_{12}$, the Levi subgroup is very complicated and it
is difficult to use multiplicativity of $\gamma$-functions (part 3 of
Theorem 3.5 of [Sh1]).
Asgari's idea is to use $G\Spin_{12}$, instead of
$\Spin_{12}$. The Levi subgroup is ${\rm GL}_3\times G\Spin_6$ and
multiplicativity of $\gamma$-functions (part 3 of Theorem 3.5 of
[Sh1]) becomes transparent.

The cases $E_6-1$ and $E_7-1$ are dealt with case by case analysis.
We first consider the case $E_6-1$ ($m=3$):

\demo{Case 1} If $\pi_1$ is not a discrete series, then by multiplicativity
of $\gamma$-functions (part 3 of Theorem 3.5 of [Sh1]),
$\gamma(s,\pi,r_1,\psi)$ is a product of $\gamma$-functions for
rank-one situations for ${\rm SL}_3\times
 {\rm SL}_3$. Hence $L(s,\pi,r_1)$ is a product of  $L$-functions for ${\rm SL}_3\times {\rm SL}_3$, and is therefore holomorphic for
${\rm Re}\, s>0$.
\enddemo

\demo{Case 2} If $\pi_1$ is a special representation, given as the
subrepresentation of $\Ind\, \mu|\ |^{\frac 12}\otimes\mu|\ |^{-\frac 12}$,
then by multiplicativity of $\gamma$-functions (part 3 of Theorem 3.5
of [Sh1]),
$$\gamma(s,\pi,r_1,\psi)=\gamma\left(s+\tfrac
12,\sigma_1,\psi\right)\gamma\left(s-\tfrac 12,\sigma_2,\psi\right),
$$
where $\sigma_1,\sigma_2$ are tempered representations of $F$-points of $\bM'$ for which the
derived group $\bM'_D$ is ${\rm SL}_3\times {\rm SL}_3$, and the $\gamma(s,\sigma_i,\psi)$'s are the Rankin-Selberg $\gamma$-factors for ${\rm
GL}_3\times {\rm GL}_3$. Note that $L(s,\sigma_i)$ is holomorphic for ${\rm Re}\, s>0$ and hence the only possible poles of
$L(s,\pi,r_1)$ are with ${\rm Re}\, s=\frac 12$, which is excluded.
\enddemo

\demo{Case 3} Representation $\pi_1$ is supercuspidal.
\enddemo

\demo{Case 3.1} If $\pi_2$ is supercuspidal, then $L(s,\pi,r_2)$ is
trivial unless $\sigma$ is also supercuspidal and we are reduced to
a case of Proposition 7.2 of [Sh1].
\enddemo

\demo{Case 3.2} Suppose $\pi_2$ is not a discrete series. Then by
multiplicativity of $\gamma$-functions (part 3 of Theorem 3.5 of
[Sh1]), $\gamma(s,\pi,r_1,\psi)$ is a product of
$\gamma$-functions for rank-one situations for either $D_5-2$ or ${\rm SL}_2\times
{\rm SL}_3$.
Hence $L(s,\pi,r_1)$ is a product of
$L$-functions for either $D_5-2$ or ${\rm SL}_2\times
{\rm SL}_3$, and it is therefore holomorphic for ${\rm Re}\, s>0$.
\enddemo

\demo{Case 3.3} Suppose $\pi_2$ is a special representation, given as the
subrepresentation of $\Ind\, \mu|\ |\otimes\mu\otimes \mu|\ |^{-1}$.
Then again by multiplicativity of $\gamma$-functions (part 3 of Theorem
3.5 of [Sh1]),
$$
\gamma(s,\pi,r_1,\psi)=\gamma(s+1,\sigma_1,\psi)\gamma(s,\sigma_2,
\psi)\gamma(s-1,\sigma_3,\psi),
$$
where the $\sigma_i$'s are tempered representations of $F$-points of $\bM'$ whose derived
group is ${\rm SL}_2\times {\rm SL}_3$, and the $\gamma(s,\sigma_i,\psi)$'s are the Rankin-Selberg $\gamma$-factors for ${\rm GL}_2\times {\rm
GL}_3$. Note that $L(s,\sigma_i)=1$ unless $\sigma$ is of the form
$\sigma=\Ind\, \tau\otimes\eta$, where $\tau$
is a supercuspidal representation of ${\rm GL}_2(F)$.
In this case, we use the argument in case 3.2 when $\pi_2$ is not in the discrete series.
\enddemo

Next we look at the case $E_7-1$ ($m=4$):

\demo{Case 1} If $\pi_1$ is not a discrete series, then by multiplicativity
of $\gamma$-functions (part 3 of Theorem 3.5 of [Sh1]),
$\gamma(s,\pi,r_1,\psi)$ is a product of $\gamma$-functions for
rank-one situations for ${\rm SL}_3\times {\rm SL}_4$. Hence $L(s,\pi,r_1)$ is a product of $L$-functions for 
${\rm SL}_3\times {\rm SL}_4$, and is consequently holomorphic for ${\rm Re}\, s>0$.
\enddemo

\demo{Case 2} Suppose $\pi_1$ is a special representation, given as the
subrepresentation of $\Ind\, \mu|\ |^{\frac 12}\otimes\mu|\ |^{-\frac 12}$;
then
$$\gamma(s,\pi,r_1,\psi)=\gamma\left(s+\tfrac
12,\sigma_1,\psi\right)\gamma\left(s-\tfrac 12,\sigma_2,\psi\right),
$$
where $\sigma_i$'s are tempered representations of $F$-points of $\bM'$ for which the
derived group $\bM'_D$ is ${\rm SL}_3\times {\rm SL}_4$, and the $\gamma(s,\sigma_i,\psi)$'s are the Rankin-Selberg $\gamma$-factors for ${\rm
GL}_3\times {\rm GL}_4$. Since $L(s,\sigma_i)$ is holomorphic for ${\rm Re}\, s>0$,
the only poles of $L(s,\pi,r_1)$ are at ${\rm Re}\, s=\frac 12$, which
are excluded.
\enddemo

\demo{Case 3} Representation $\pi_1$ is supercuspidal.
\enddemo

\demo{Case 3.1} Suppose $\pi_2$ is not a discrete series. Then by
multiplicativity of $\gamma$-functions (part 3 of Theorem 3.5 of
[Sh1]), $\gamma(s,\pi,r_1,\psi)$ is a product of
$\gamma$-functions for rank-one situations for $D_6-2$, or
${\rm SL}_2\times {\rm SL}_4$. Hence $L(s,\pi,r_1)$ is a product of
$L$-functions for either $D_6-2$ or ${\rm SL}_2\times
{\rm SL}_4$, and is therefore holomorphic for ${\rm Re}\, s>0$.
\enddemo

\demo{Case 3.2} Suppose $\pi_2$ is supercuspidal. If $\sigma$ is not a
discrete series, then by multiplicativity of $\gamma$-functions (part
3 of Theorem 3.5 of [Sh1]), $\gamma(s,\pi,r_1,\psi)$ is a product
of $\gamma$-functions for rank-one situations for either $E_6-1$,
$D_5-2$,
or ${\rm SL}_2\times {\rm SL}_3$. Hence $L(s,\pi,r_1)$ is a product
of $L$-functions for either cases $E_6-1$, $D_5-2$,
or ${\rm SL}_2\times {\rm SL}_3$, and is consequently holomorphic for ${\rm Re}\, s>0$ by the above result for the case $E_6-1$.

If $\sigma$ is supercuspidal, we are reduced to a case of Proposition
7.2 of [Sh1].

Suppose $\sigma$ is given as the subrepresentation of
$\Ind\, \rho|{\rm det}|^{\frac 12}\otimes\rho|{\rm det}|^{-\frac 12}$, where $\rho$ is a supercuspidal representation of ${\rm GL}_2(F)$.
Then
$$
\gamma(s,\pi,r_1,\psi)=\gamma\left(s+\tfrac
12,\sigma_1,\psi\right)\gamma\left(s-\tfrac 12,\sigma_2,\psi\right),
$$
where $\sigma_i$'s are tempered representations of $F$-points of $\bM'$ for which the
derived group $\bM'_D$ is ${\rm SL}_2\times {\rm SL}_2\times {\rm SL}_2$.
This is case $D_4-2$, and by the above result, $L(s,\sigma_i)$ is holomorphic for ${\rm Re}\, s>0$. Hence $L(s,\pi,r_1)$ can have a pole at most with ${\rm Re}\, s=\frac 12$,  which is excluded.

If $\sigma$ is given as the subrepresentation of $\Ind\, \mu|\ |^{\frac
32}\otimes\mu|\ |^{\frac 12}\otimes\mu|\ |^{-\frac 12}\otimes\mu|\
|^{-\frac 32}$, then all rank-one situations are that of ${\rm SL}_3\times {\rm SL}_2$,
which is non-self-conjugate. Hence $L(s,\pi,r_1)=1$.
\enddemo

\demo{Case 3.3} Suppose $\pi_2$ is a special representation, given as the
subrepresentation of $\Ind\, \mu|\ |\otimes\mu\otimes\mu|\ |^{-1}$. Then
$$
\gamma(s,\pi,r_1,\psi)=\gamma(s+1,\sigma_1,\psi)\gamma(s,\sigma_2,
\psi)\gamma(s-1,\sigma_3,\psi),
$$
where $\sigma_i$'s are tempered representations of $F$-points of  $\bM'$ for which the
derived group $\bM'_D$ is ${\rm SL}_2\times {\rm SL}_4$, and $\gamma(s,\sigma_i,\psi)$'s are the Rankin-Selberg $\gamma$-factors for ${\rm GL}_2\times {\rm GL}_4$. Note that $L(s,\sigma_i)$ is not trivial only when $\sigma$ is of the form $\Ind\, \tau\otimes\tau'$, where $\tau,\tau'$ are either supercuspidal representations of ${\rm GL}_2(F)$, or are in the discrete series, given each as the irreducible subrepresentation of $\Ind\, \rho|{\rm det}|^{\frac 12}\otimes\rho|{\rm det}|^{-\frac 12}$, where $\rho$ is a supercuspidal representation of ${\rm GL}_2(F)$.
If $\sigma$ is of the first form, then we are in the first part of Case 3.2. \pagebreak 
If $\sigma$ is in the discrete series, then the $L$-function is given by $(1-u q^{-\frac 12-s})^{-1}$, where $u$ is a complex number with absolute value 1. Hence $L(s,\sigma,r_1)$ can have a pole only at ${\rm Re}\, s=\frac 12$, which is excluded.
\hfill\qed
\enddemo

\nonumproclaim{Proposition  4.2} The normalized local intertwining operators
$$N(s\tilde\alpha,\pi_v,w_0)$$ are holomorphic and nonzero for ${\rm Re}\, s\geq
\frac 12$ and for all $v${\rm .}
\endproclaim

\demo{{P}roof} We proceed as in [Ki2, Prop.\ 3.4]. If $\pi_v$
is tempered, then the unnormalized operator is holomorphic and
nonzero for ${\rm Re}\, s>0$. By Lemma 4.1, $L(s,\pi_v,r_i)$ is holomorphic for ${\rm Re}\, s>0$. Hence the normalized operator $N(s\tilde\alpha,\pi_v,w_0)$ is holomorphic and nonzero for ${\rm Re}\, s>0$.

If $\pi_v$ is nontempered, we write $I(s,\pi_v)$ as in [Ki1, p.~841],
$$
I(s,\pi_v)=I(s\tilde\alpha+\Lambda_0,\pi_0)=\Ind_{\bold M_0(F_v)\bold N_0(F_v)
}^{\bold G(F_v)}\, \pi_0\otimes q^{\langle s\tilde\alpha+\Lambda_0, H_{P_0}
(\ )\rangle},
$$
where $\pi_0$ is a tempered representation of $\bold M_0(F_v)$ and
$\bold P_0=\bold M_0\bold N_0$ is another parabolic subgroup of $\bold G$.
Moreover, $(\Lambda_0,\pi_0)$ is the Langlands parameter of $\pi_v$.
We can identify
the normalized operator $N(s\tilde\alpha,\pi_v,w_0)$ with the normalized
operator $N(s\tilde\alpha+\Lambda_0,\pi_0, w_0)$, which is a
product of rank-one operators attached to tempered representations
(cf. [Z, Prop.~1]).

By direct observation, we see that all the rank-one operators are operators attached to tempered representations of a parabolic subgroup whose Levi subgroup has a derived group isomorphic to
${\rm SL}_k\times {\rm SL}_l$ inside a group whose derived group is ${\rm SL}_{k+l}$, except one case.
It is in the case $E_7-1$, when $\pi_1,\pi_2$ are tempered, and $\sigma$ is the nontempered representation, given by the quotient of $\Ind\, |\det|^r\rho\otimes |\det|^{-r}\rho$, where $\rho$ is a tempered representation of ${\rm GL}_2(F)$ and $0<r<\frac 12$.
The rank-one operator is that of the case $D_5-2$, attached to $\pi_1,\pi_2$ and $\rho$.
However in this case, $s\tilde\alpha+\Lambda_0$ is in the corresponding positive Weyl chamber for ${\rm Re}\, s\geq\frac 12$, and hence $N(s\tilde\alpha+\Lambda_0,\pi_0,w_0)$ is holomorphic for ${\rm Re}\, s\geq \frac 12$ [Ki1, Lemma 2.4].

The rank-one operators attached to tempered representations of ${\rm SL}_k\times {\rm SL}_l$ are then restrictions to
${\rm SL}_{k+l}$ of corresponding standard operators for ${\rm GL}_{k+l}$.
By [MW2, Prop.~I.10] one knows that these rank-one operators are
holomorphic for ${\rm Re}\, s>-1$.
Hence by identifying roots of $G$ with
respect to a parabolic subgroup, with those of $G$ with respect to
the maximal torus, it is enough to check $\langle
s\tilde\alpha+\Lambda_0,\beta^{\vee}\rangle>-1$ for all positive
roots $\beta$ if ${\rm Re}\, s\geq \frac 12$.
We do this case by case as follows: 
\enddemo

\demo{Case $D_5-2$} In the notation of Bourbaki [Bou],
$\tilde\alpha=e_1+e_2+e_3$;
$\Lambda_0=r_1e_1-r_1e_3+r_2(e_4-e_5)+r_3(e_4+e_5)$, where
$\frac 12>r_1,r_2,r_3\geq 0$.
Here $\pi_{1v}$ is tempered if $r_1=0$.
Hence
$$
s\tilde\alpha+\Lambda_0=(s+r_1)e_1+se_2+(s-r_1)e_3+(r_2+r_3)e_4+
(-r_2+r_3)e_5.
$$
We observe that the least value of
${\rm Re}(\langle s\tilde\alpha+\Lambda_0,\beta^{\vee}\rangle)$ is
${\rm Re}\, s-(r_1+r_2+r_3)$ when $\beta=e_3-e_4$. It is larger than $-1$  if ${\rm Re}\, s\geq \frac 12$. Consequently,
 $N(s\tilde\alpha+\Lambda_0,\pi_0,w_0)$ is holomorphic for ${\rm Re}\, s\geq \frac 12$. By Zhang's lemma (cf. [Ki2, Lemma 1.7] and [Z]), it is
nonzero as well.
\enddemo

\demo{Case $E_6-1$} In the notation of Bourbaki [Bou],
$\tilde\alpha=\varpi_4=2\alpha_1+3\alpha_2+4\alpha_3+6\alpha_4+4
\alpha_5+2\alpha_6$; $\Lambda_0=r_1\alpha_1+r_1\alpha_3+r_2\alpha_5+r_2\alpha_6+r_3\alpha_2$,
where $\frac 12>r_1,r_2,r_3\geq 0$. Hence
$$\align
s\tilde\alpha+\Lambda_0=&\ (2s+r_1)\alpha_1+ (3s+r_3)\alpha_2+(4s+r_1)\alpha_3
\cr
&\ +6s\alpha_4+(4s+r_2)\alpha_5+(2s+r_2)\alpha_6. \endalign
$$
We observe that the least value of
${\rm Re}(\langle s\tilde\alpha+\Lambda_0,\beta^{\vee}\rangle)$ is
${\rm Re}\, s-(r_1+r_2+r_3)$ when $\beta=\alpha_4$. It is larger than $-1$, if ${\rm Re}\, s\geq \frac 12$. Consequently, $N(s\tilde\alpha+\Lambda_0,\pi_0,w_0)$ is holomorphic for ${\rm Re}\, s\geq \frac 12$. By Zhang's lemma (cf. [Ki2, Lemma 1.7] and [Z]), it is nonzero as well.
\enddemo
 
\demo{Case $E_7-1$} In the notation of Bourbaki [Bou],
$\tilde\alpha=\varpi_4=4\alpha_1+6\alpha_2+8\alpha_3+12\alpha_4+9
\alpha_5+6\alpha_6+3\alpha_7$;
$\Lambda_0=r_1\alpha_1+r_1\alpha_3+r_2\alpha_5+(r_2+r_3)\alpha_6+r_2
\alpha_7+r_4\alpha_2$, where $\frac 12>r_1,r_4\geq 0$ and $\frac 12>r_2\geq
r_3\geq 0$.
Hence
$$ \align
s\tilde\alpha+\Lambda_0=&\ (4s+r_1)\alpha_1+ (6s+r_4)\alpha_2+(8s+r_1)\alpha_3+
12s\alpha_4\cr
&+\ (9s+r_2)\alpha_5+(6s+r_2+r_3)\alpha_6+(3s+r_2)\alpha_7. \endalign
$$
We observe that the least value of
${\rm Re}(\langle s\tilde\alpha+\Lambda_0,\beta^{\vee}\rangle)$ is
${\rm Re}\, s-(r_1+r_2+r_4)$ when $\beta=\alpha_4$. It is larger than $-1$ if ${\rm Re}\, s\geq \frac 12$. Consequently, 
$N(s\tilde\alpha+\Lambda_0,\pi_0, w_0)$ is holomorphic for ${\rm Re}\, s\geq \frac 12$. By Zhang's lemma (cf. [Ki2, Lemma 1.7]), it is nonzero as
well. \hfill\qed
\enddemo

\def\bB{{\bf B}}
\def\bT{{\bf T}}
\def\bU{{\bf U}}
\def\bA{{\bf A}}
\def\bG{{\bf G}}
\def\bP{{\bf P}}
\def\bM{{\bf M}}
\def\bN{{\bf N}}
\def\bZ{{\Bbb Z}}

\def\bR{{\Bbb R}}
\def\bC{{\Bbb C}}

\section{Functorial products for ${\rm GL}_2\times {\rm GL}_3$}

In this section we will prove our main theorem. Recall that $\pi_1=\hskip-6pt  \mathbold{\otimes}_v \pi_{1v}$ and $\pi_2=\hskip-6pt  \mathbold{\otimes}_v \pi_{2v}$ are cuspidal representations of
${\rm GL}_2(\bAB_F)$ and ${\rm GL}_3(\Bbb A_F)$, respectively. Moreover, for each $v$,
$\pi_{1v}\boxtimes\pi_{2v}$ is the irreducible admissible representation of ${\rm GL}_6(F_v)$ attached to $\pi_{1v}\otimes\pi_{2v}$ through the local Langlands correspondence. Let $\pi_1\boxtimes\pi_2=\hskip-6pt  \mathbold{\otimes}_v (\pi_{1v}\boxtimes\pi_{2v})$. It is an irreducible admissible representation of ${\rm GL}_6(\Bbb A_F)$.

\nonumproclaim{Theorem 5.1}
 The representation $\pi_1\boxtimes\pi_2$ of ${\rm GL}_6(\Bbb A_F)$ is automorphic{\rm .} It is irreducibly induced from cuspidal
representations{\rm ,} i.e{\rm .,} $\Pi=\Ind \, \sigma_1\otimes\cdots\otimes\sigma_k$, where $\sigma_i$\/{\rm '}\/s 
are cuspidal representations of ${\rm
GL}_{n_i}(\Bbb A_F)${\rm ,}
$n_i>1${\rm .}
The cases\break $k=3,\ n_1=n_2=n_3=2${\rm ,} or $k=2,\ n_1=2${\rm ,} $n_2=4${\rm ,}
 occur if and only if $\pi_2$ is a twist of ${\rm Ad}(\pi_1)$ by a gr{\rm \"{\it o}}ssencharacter{\rm .}
\endproclaim

We will first prove the last statement, assuming the earlier ones.
Suppose $n_i=2$ for some $i$ and consider the partial $L$-function
$$
L_S (s,\Pi\times\ \tilde\sigma_i)
$$
which has a pole at $s=1$ (cf.~[JS1] and [Sh4]), where $S$ is a finite set of places outside of which everything is unramified.
It equals
$$
L_S(s, (\pi_1\boxtimes \pi_2)\times\ \tilde\sigma_i)=L_S (s,\pi_2\times(\pi_1\boxtimes\tilde\sigma_i)).
$$

Since any quadratic base change of $\pi_2$ is cuspidal (Theorem 4.2 of [AC]), we may assume that $\pi_1\boxtimes\tilde\sigma_i$ is not an automorphic induction from a quadratic extension of $F$.
For then $L_S(s,\pi_2\times (\pi_1\boxtimes\tilde\sigma_i))$ would never have a pole.
Consequently, by the discussion in the proof of Part II of Lemma 3.1.1 of [R1], we may assume neither $\pi_1$, nor $\tilde\sigma_i$, is monomial.
By [R1, Th.~M, pg.~54], $\pi_1\boxtimes \tilde\sigma_i$ is cuspidal unless $\sigma_i\cong\pi_1\otimes\eta$ for some gr\"ossencharacter $\eta$
and only then. In this case, $\pi_1\boxtimes\tilde\sigma_i={\rm Ad}(\pi_1)\otimes\eta^{-1}\boxplus\eta^{-1}$.
Since $L_S(s,\pi_2\times ({\rm Ad}(\pi_1)\otimes\eta^{-1}))$ must now have a pole at $s=1$, we have $\pi_2\cong {\rm Ad}(\pi_1)\otimes\eta$.
Conversely, if $\pi_2\cong {\rm Ad}(\pi_1)\otimes\eta$, then
$$
L_S\left(s,\Pi \times \ (\tilde\pi_1\otimes \eta^{-1})\right)=L_S\left(s,\pi_2\times ({\rm Ad}(\pi_1)\otimes \eta^{-1})\right)
L_S\left(s,\pi_2\otimes\eta^{-1}\right)
$$
will have a pole at $s=1$.
Thus $\sigma_i=\pi_1\otimes\eta$ appears in the inducing data for $\Pi=\pi_1\boxtimes\pi_2$.

We now proceed to prove the main part of Theorem 5.1 and it is appropriate to remind the reader of the definitions of $L$-functions and root numbers
using our method ([Sh1]). We shall freely use definitions and results from\break [Sh1, 4,  7].
We start by applying Theorem 3.5 of [Sh1] to each of our four cases explained in Section 3.
This allows us to define a $\gamma$-factor $\gamma(s,\pi_{1v}\times \pi_{2v}\times\sigma_v,\psi_v)$ at each place $v$.
When $v=\infty$, we will define $L(s,\pi_{1v}\times \pi_{2v}\times\sigma_v)$ and $\var(s,\pi_{1v}\times \pi_{2v}\times\sigma_v,\psi_v)$ using parametrization ([La4] and [Sh7]).
They satisfy
$$ \align
&\gamma(s,\pi_{1v}\times\pi_{2v}\times\sigma_v,\psi_v)\tag5.1\cr
&\qquad =\var(s,\pi_{1v}\times\pi_{2v}\times\sigma_v,\psi_v)
L(1-s,\tilde\pi_{1v}\times\tilde\pi_{2v}\times\tilde\sigma_v)/L(s,\pi_{1v}\times\pi_{2v}\times\sigma_v).\endalign
$$  
Now, suppose $v < \infty$ and $\pi_{iv},i=1,2$, and $\sigma_v$ are all tempered.
As explained in [Sh1, \S 7], we then define $L(s,\pi_{1v}\times \pi_{2v}\times\sigma_v)$ as the inverse of a polynomial in $q_v^{-s}$ whose
constant term is 1 and which has the same zeros as $\gamma(s,\pi_{1v}\times\pi_{2v}\times\sigma_v,\psi_v)$. The root number
$\var(s,\pi_{1v}\times\pi_{2v}\times \sigma_v,\psi_v)$ is now defined using (5.1). It is now clear that if $\pi_{iv}$, $i=1,2$, and $\sigma_v$ are
all tempered, then the $\gamma$-function determines the root number and the $L$-function uniquely.

Defining the $L$-function for nontempered ones is more delicate.
We need to use Langlands classification and multiplicativity (Part 3 of Theorem 3.5 in [Sh1]).
Then $\gamma(s,\pi_{1v}\times\pi_{2v}\times\sigma_v,\psi_v)$ can be written as a product of $\gamma$-functions defined by the quasi-tempered data which gives the Langlands parameter for $\pi_{1v}\otimes\pi_{2v}\otimes\sigma_v$ as its unique Langlands subrepresentation.
(To apply multiplicativity, one must consider the representation $\pi_{1v}\otimes\pi_{2v}\otimes\sigma_v$ as a subrepresentation, in fact the unique one, of the induced representation with a parameter in the negative Weyl chamber.
But even if one uses the standard module, i.e., the one with a parameter in the positive Weyl chamber, then the two defining local coefficients differ
only by a local coefficient which is defined by a standard intertwining operator for the group $M$, rather than one for $G$. Consequently, the
corresponding local coefficient is independent of $s$ and can be used to normalize the defining operator, which can have a pole or a zero, leading to a
local coefficient, now defined by this normalized operator, equal to $1$.) More precisely, each of these $\gamma$-functions is defined by means of a
maximal Levi in a connected reductive group via Theorem 3.5 of [Sh1]. For each of them an $L$-function is defined by means of analytic continuation
of the $L$-function from the tempered data, as explained earlier, to the quasi-tempered ones. In our case, they are all triple product $L$-functions
coming from groups of lower rank. The $L$-function $L(s,\pi_{1v}\times\pi_{2v}\times\sigma_v)$ is then the product of these $L$-functions. This
is what we call multiplicativity for $L$-functions. This definition agrees completely with parametrization and in particular when $v=\infty$
(cf.~[La4]). The root numbers are now defined by means of (5.1).
We observe that, as discussed above, both the $L$-function and the root number are independent of whether the representation is considered as a
Langlands subrepresentation or a Langlands quotient. More generally, local coefficients are defined to depend only on the equivalence classes of
inducing data, by setting any local coefficient which can be defined by means of a intertwining operator between different realizations of inducing
data, which in fact does not depend on the complex parameter, equal to 1. Moreover, as discussed in Section 9 of [Sh1], the same definition can be
used to define the factors even for Langlands quotients which are not generic. We refer to Sections 7 and 9 of [Sh1] for more detail. The reader can
consult our proof of Lemma 4.1 in Section 4 to see how multiplicativity is applied and what the lower rank factors are.

\nonumproclaim{Proposition  5.2} Suppose either $\pi_{1v}${\rm ,}
 or $\pi_{2v}${\rm ,} is not supercuspidal{\rm .} Then the equalities {\rm (3.2.1)} and {\rm (3.2.2)} hold{\rm .} More precisely{\rm ,}
$$
L(s,\pi_{1v}\times \pi_{2v}\times \sigma_v)=L(s,(\pi_{1v}\boxtimes \pi_{2v})\times \sigma_v)\tag5.2.1
$$
and
$$
\varepsilon(s,\pi_{1v}\times\pi_{2v}\times\sigma_v,\psi_v)=\varepsilon(s,(\pi_{1v}\boxtimes \pi_{2v})\times \sigma_v,\psi_v)\tag5.2.2
$$
for every irreducible admissible generic representation $\sigma_v$ of ${\rm GL}_n(F_v)${\rm ,} $n=1,2,3,4${\rm .}
\endproclaim

\demo{{P}roof} By multiplicativity of $\gamma$-factors, it is enough to
prove the equality of corresponding $\gamma$-functions for supercuspidal representations $\sigma_v$. If either $\pi_{1v}$, or $\pi_{2v}$, is a component of a principal series, then by multiplicativity of $\gamma$ and $L$-factors, our assertion follows as in Proposition 3.6. So the only case left 
is when $\pi_{1v}$ is supercuspidal, and $\pi_{2v}=\Ind \, \rho\otimes\mu$, where $\rho$ is a supercuspidal representation of ${\rm GL}_2(F_v)$
and $\mu$ is a character of $F_v^{\times}$. By multiplicativity of $\gamma$-factors,
$$
\gamma(s,\pi_{1v}\times\pi_{2v}\times\sigma_v,\psi_v)=\gamma(s,\pi_{1v}\times \rho\times \sigma_v,\psi_v)\gamma(s,(\pi_{1v}\otimes\mu)\times\sigma_v,\psi_v).
$$
It is enough to prove that
$$\gamma(s,\pi_{1v}\times \rho\times \sigma_v,\psi_v)=\gamma(s,(\pi_{1v}\boxtimes\rho)\times \sigma_v,\psi_v),
$$
where $\pi_{1v}\boxtimes\rho$ is an irreducible representation of ${\rm GL}_4(F_v)$, given by the local Langlands correspondence. By
[Sh1, Prop.~5.1], we can find cuspidal representations $\pi_1=\hskip-6pt  \mathbold{\otimes}_w \pi_{1w}, \pi'=\hskip-6pt  \mathbold{\otimes}_w
\pi_w'$ of ${\rm GL}_2(\Bbb A_F)$ and a cuspidal representation $\sigma=\hskip-6pt  \mathbold{\otimes}_w \sigma_w$ of
${\rm GL}_n(\Bbb A_F)$ such that $\pi_{1w},\pi_w'$, and $\sigma_w$ are all unramified for $w\ne v$, $w<\infty$, and $\pi_v'=\rho$. Let $\pi_1\boxtimes\pi'$ be the functorial product, constructed by Ramakrishnan 
[R1].
Consider the case $D_{n+2}-2$ in [Sh2], and compare the functional equations for $L(s,\pi_1\times\pi'\times\sigma)$ and $L(s,(\pi_1\boxtimes\pi')\times\sigma)$;
$$\align
L(s,\pi_1\times\pi'\times\sigma)&=\varepsilon(s,\pi_1\times
\pi'\times\sigma)L(1-s,\tilde\pi_1\times\tilde\pi'\times\tilde\sigma),\\
L(s,(\pi_1\boxtimes\pi')\times\sigma)&=
\varepsilon(s,(\pi_1\boxtimes\pi')\times\sigma)L(1-s,(\tilde\pi_1\boxtimes\tilde\pi')\times\tilde\sigma).
\endalign
$$
Since $\pi_{1w}, \pi_{w}'$, and $\sigma_w$ are all unramified for $w\ne v$, $w<\infty$,
$$
\gamma(s,\pi_{1w}\times\pi_w'\times\sigma_w,\psi_w)=\gamma(s,(\pi_{1w}\boxtimes\pi_w')\times\sigma_w,\psi_w),
$$
for all $w\ne v$.
Hence we have
$$
\gamma(s,\pi_{1v}\times\pi_v'\times\sigma_v,\psi_v)=\gamma(s,(\pi_{1v}\boxtimes\pi_v')\times\sigma_v,\psi_v),
$$
since similar inequalities hold at archimedean places, the factors being those of Artin [Sh7]. 
\enddemo 

It is more difficult to prove the equalities (3.2.1) and (3.2.2) if both $\pi_{1v}$ and $\pi_{2v}$ are supercuspidal.
Let $T$ be the finite set of places $v|2$ such that
$\pi_{1v}$ is an extraordinary supercuspidal representation of ${\rm GL}_2(F_v)$, while $\pi_{2v}$ is
a supercuspidal representation of ${\rm GL}_3(F_v)$ attached to a character of a
nonnormal cubic extension $L$ of $F_v$.

\nonumproclaim{Proposition  5.3} Suppose $\pi_{1v}$ and $\pi_{2v}$ are both supercuspidal{\rm .}
Then\/{\rm :}

\itemitem{\rm a)} Equalities {\rm (5.2.1)} and {\rm (5.2.2)} hold for every irreducible admissible generic
 representation $\sigma_v$ of ${\rm GL}_n (F_v),\
n=1,2,3${\rm .}

\itemitem{\rm b)} Suppose $v\not\in T${\rm .}
Then {\rm (5.2.1)} and {\rm (5.2.2)} hold for every irreducible admissible generic representation $\sigma_v$ of ${\rm GL}_n(F_v),\ n=1,2,3,4${\rm .}
\endproclaim

\demo{Proof} By multiplicativity of $\gamma$-factors, it is enough to
prove the equality of corresponding $\gamma$-functions for supercuspidal representations $\sigma_v$.

Let $\pi_{1v},\pi_{2v},\sigma_v$ be supercuspidal representations of ${\rm GL}_2(F_v)$, ${\rm GL}_3(F_v)$, ${\rm GL}_n(F_v)$, resp., where
$n=2,3,4$. Let $\rho_{1v},\rho_{2v},\rho_{3v}$ be the corresponding Weil group representations. We shall prove the equivalent equality:
$$\gamma(s,\pi_{1v}\times\pi_{2v}\times\sigma_v,\psi_v)=\gamma(s,\rho_{1v}\otimes\rho_{2v}\otimes\rho_{3v},\psi_v).\tag5.3.1
$$
Here $\gamma(s,\pi_{1v}\times\pi_{2v}\times\sigma_v,\psi_v)$ is the $\gamma$-factor defined in [Sh1], while $\gamma(s,\rho_{1v}\otimes\rho_{2v}\otimes\rho_{3v},\psi_v)$ is the corresponding Artin $\gamma$-factor. Since
$$\gamma(s,(\pi_{1v}\boxtimes\pi_{2v})\times\sigma_v,\psi_v)=\gamma(s,(\rho_{1v}\otimes\rho_{2v})\otimes\rho_{3v},\psi_v),
$$
and $\pi_{1v}\boxtimes\pi_{2v}$ corresponds to $\rho_{1v}\otimes\rho_{2v}$ by
the local Langlands correspondence, Proposition 5.3 then follows,
namely,
$$\gamma(s,\pi_{1v}\times\pi_{2v}\times\sigma_v,\psi_v)=\gamma(s, (\pi_{1v}\boxtimes\pi_{2v})\times\sigma_v,\psi_v).
$$

\demo{Proof of {\rm (5.3.1)}} We first consider the case where $v|2$.
Then $v\nmid 3$ and therefore $\pi_{2v}$ is tame and by classification of tame supercuspidal representations there exists a cubic extension $L$ of $F_v$ such that $\pi_{2v}$ corresponds to Ind$(W_{F_v},W_L,\eta)$, where $\eta$ is a character of $L^*$.
To proceed, we embed $\pi_{1v},\sigma_v$ as local components of cuspidal representations
$\pi_1,\sigma$
with unramified finite components everywhere else.

Next, we choose a cubic extension $\bL/F$ such that $\bL_w=L,\ w|v$; we also choose a gr\"ossencharacter $\chi$ of $\bL$ satisfying
$\chi_w=\eta$. We   assume further that $\chi_u$ is unramified for every $u < \infty,\ u\not= w$.
Let $\pi_2$ be the cuspidal representation of ${\rm GL}_3(\bAB_F)$ corresponding to $\rho_2=\Ind (W_F,W_{\bL},\chi)$ [AC], [JP-SS2].
Observe that $\bL/F$ need not be normal.
\enddemo

Let $\pi_1^\bL$ and $\sigma^\bL$ be the base changes of $\pi_1$ and $\sigma$ to $\bL$.
Observe that when $v\in T$, we are assuming $n\leq 3$ so that the nonnormal base change of $\sigma$ to $\bL$ exists [JP-SS2].
We compare the functional equations for $L(s,\pi_1\times\pi_2\times\sigma)$ and $L(s,(\pi_1^\bL\otimes\chi)\times\sigma^\bL)$.
More precisely,
$$
L(s,\pi_1\times\pi_2\times\sigma)=\var(s,\pi_1\times\pi_2\times\sigma)
L(1-s,\tilde\pi_1\times\tilde\pi_2\times\tilde\sigma)
$$
and
$$
L\big(s,\big(\pi_1^\bL\otimes\chi\big)\times\sigma^\bL\big)=\var\big(s,\big(\pi_1^\bL\otimes\chi\big)
\times\sigma^\bL\big)
L\big(1-s,\big(\tilde\pi_1^\bL\otimes\chi^{-1}\big)\times\tilde\sigma^\bL\big).
$$

Let $L_u=\bL\otimes_F F_u$.
Then either $L_u/F_u$ is a cubic field extension, $L_u\simeq F_u\oplus L_w$ with $L_w/F_u$ a quadratic extension, or $L_u\simeq F_u\oplus F_u\oplus F_u$.
Thus if $w$ is a place of $\bL$ over $u$, then the local base changes to $L_w$ are either cubic, quadratic, or trivial.

To continue we need a general discussion as follows.

Let $K/F$ be a finite separable extension of local fields.
Fix a nontrivial additive character $\psi_F$ of $F$ and let $\psi_{K/F}=\psi_F\cdot Tr_{K/F}$.
Let $\lambda (K/F,\psi_F)$ be the corresponding Langlands $\lambda$-function (cf.~[La5], [AC]).
Let $\pi_1$ and $\pi_2$ be two irreducible admissible representations of ${\rm GL}_m(F)$ and ${\rm GL}_n(F)$, respectively.
Fix continuous representations $\rho_1$ and $\rho_2$ of the Deligne-Weil group $W_F$ attached to $\pi_1$ and $\pi_2$ by the local Langlands correspondence
 [HT], [He], respectively.
Assume $\rho_2=\text{ Ind}(W_F,W_K,\tau_2)$ with $\dim\tau_2=p$.
Let $\pi_1^K$ be the base change of $\pi_1$, i.e.~the representation attached to $\rho_1|W_K$.
Then
$$
\var(s,\rho_1\otimes\rho_2,\psi_F)=\lambda (K/F,\psi_F)^{mp}
\var(s,\rho_1|W_K\otimes\tau_2,\psi_{K/F}).
$$
Consequently
$$
\var(s,\pi_1\times\pi_2,\psi_F)=\lambda (K/F,\psi_F)^{mp}
\var(s,\pi_1^K\times\beta_2,\ \psi_{K/F})
$$
where $\beta_2$ is the representation of ${\rm GL}_p(K)$, $p=\dim\tau_2$, attached to $\tau_2$.
The identities for $L$-functions hold without any $\lambda$-functions.

In our setting $p=1,\tau_2=\eta$, and $\rho_{2v}=\Ind(W_{F_v},W_L,\eta)$.
Thus
$$
\var(s,(\rho_{1v}\otimes\rho_{3v})\otimes \rho_{2v},\psi_v)=\lambda (L/F_v,\psi_v)^{2n}\var(s,{\rm Res}_{W_L}
(\rho_{1v}\otimes\rho_{3v})\otimes\eta,\ \psi_{L/F_v})
$$
since $\dim(\rho_{1v}\otimes\rho_{3v})=2n$.
By the local Langlands correspondence,
$$
\align
\var(s,\rho_{1v}\otimes\rho_{2v}\otimes\rho_{3v},\psi_v)&=\lambda
(L/F_v,\psi_v)^{2n}
\var(s,(\pi_{1v}\boxtimes\sigma_v)^L\otimes\eta,\psi_{L/F_v})\\
&= \lambda (L/F_v,\psi_v)^{2n}\var(s,(\pi_{1v}^L\otimes\eta)
\times\sigma_v^L,\psi_{L/F_v}).
\endalign
$$

On the other hand for $v_1$, a place of $F$ different from $v$, we note that $\rho_2=\Ind(W_F,W_\bL,\chi)$, as a representation of $W_{v_1}=W_{F_{v_1}}$, is equal to
$$
\bigoplus_{w_1|v_1} I(\chi_{w_1}),
$$
where
$$
I(\chi_{w_1})=\Ind (W_{v_1},W_{w_1},\chi_{w_1}),
$$
with $W_{w_1}=W_{L_{w_1}}$.
Thus
$$
\align 
&\tag5.3.2\\
\noalign{\vskip-30pt}
\hskip.5in & \var(s,(\rho_{1v_1}\otimes\rho_{3v_1})\otimes\rho_{2v_1},\psi_{v_1})\cr
& \qquad=\, \prod_{w_1|v_1}
\var(s,(\rho_{1v_1}\otimes\rho_{3v_1})\otimes I(\chi_{w_1}),\psi_{v_1})\\
&\qquad=\, \prod_{w_1|v_1} \lambda
(L_{w_1}/F_{v_1},\psi_{v_1})^{2n}\var(s,(\pi_{1v_1}\boxtimes\sigma_{v_1})^{L_{w_1}}\otimes\chi_{w_1},\psi_{Lw_1/F_{v_1}})\\
&\qquad=\, \prod_{w_1|v_1} \lambda (L_{w_1}/F_{v_1},\psi_{v_1})^{2n}\var\big(s,\big(\pi_{1v_1}^{L_{w_1}}\otimes
\chi_{w_1}\big)\times\sigma_{v_1}^{L_{w_1}},\psi_{L_{w_1}/F_{v_1}}\big)
\endalign 
$$
by the local Langlands correspondence.
Similarly
$$
L(s,\rho_{1v_1}\otimes\rho_{2v_1}\otimes\rho_{3v_1})=\prod_{w_1|v_1} L\big(s,\big(\pi_{1v_1}^{L_{w_1}}\otimes\chi_{w_1}\big)
\times\sigma_{v_1}^{L_{w_1}}\big).\tag5.3.3
$$

To prove (5.3.1), it would then be enough to show that
$\var(s,\pi_{1v_1}\times\pi_{2v_1}\times\sigma_{v_1},\psi_{v_1})\text{ and }
L(s,\pi_{1v_1}\times\pi_{2v_1}\times\sigma_{v_1})$
are equal to the right-hand sides of (5.3.2) and (5.3.3), for every $v_1\not= v$, for then (5.3.1) would follow immediately by comparison of functional
equations for $L(s,\pi_1\times\pi_2\times\sigma)$ and $L(s,(\pi_1^\bL\otimes\chi)\times\sigma^\bL)$. But these last equalities are obvious since
$\pi_{iv_1},i=1,2$, and $\sigma_{v_1}$ are either unramified or archimedean.

Since when $v\notin T$, we can use cubic cyclic base change, which is available for all $n$, part b) follows as well.

If $v\nmid 2$, then $\pi_{1v}$ is dihedral, and the arguement goes as before except that this time we use a quadratic base change. \hfill\qed
\enddemo 

It remains to prove (5.3.1) when $\sigma_v$ is a supercuspidal representation of ${\rm GL}_4(F_v)$ for $v\in T$.
Namely, $v|2,\ \pi_{1v}$ is an extraordinary supercuspidal representation of ${\rm GL}_2(F_v)$, and $\pi_{2v}$ is a supercuspidal representation of ${\rm GL}_3(F_v)$ attached to a character of a nonnormal cubic extension $L$ of $F_v$.
The difficulty is that the theory of nonnormal cubic base change for ${\rm GL}_4$ is not available at present.
We need to proceed as follows.
We first construct a local lift $\Pi_v$ in the sense of Definition 3.3, and show that it differs from $\pi_{1v}\boxtimes\pi_{2v}$ by at most a quadratic character.
Then the appendix by Bushnell and Henniart [BH] proves that, in fact, $\Pi_v\simeq\pi_{1v}\boxtimes\pi_{2v}$.
We start with:

\nonumproclaim{Proposition  5.4} For each $v\in T${\rm ,}
 there exists a local lift $\Pi_v$ of $\pi_{1v}\otimes\pi_{2v}$ in the sense of Definition {\rm 3.3.}
\endproclaim

\demo{Proof}
We start by letting
$\pi_1\otimes\pi_2=\hskip-6pt  \mathbold{\otimes}_w(\pi_{1w}\otimes\pi_{2w})$ be a cuspidal representation of ${\rm GL}_2(\bAB_F)\times {\rm GL}_3(\bAB_F)$ such that $\pi_{1w}\otimes\pi_{2w}$ is unramified for all $w < \infty,\ w\not= v$ (cf.~[Sh1]).
Let $\Pi=\hskip-6pt  \mathbold{\otimes}_w \Pi_w$ be the weak lift of $\pi_1\otimes\pi_2$ with respect to the set $S=\{v\}$ obtained from Proposition 3.10.
We will show that $\Pi_v$ is a local lift of $\pi_{1v}\otimes\pi_{2v}$, i.e., that it satisfies (3.3.1) and (3.3.2).
To proceed, we must be extra careful.
In fact, we no longer know, {\it a priori}, that $\Pi_v$ is tempered, and therefore the equality of $\gamma$-functions does not directly imply that of $L$ and $\var$ factors which we used occasionally earlier, as it will require the knowledge of its Langlands parameter.
All that we know is that $\Pi_v$ is irreducible, unitary, and generic.
But fortunately this is enough.
In fact, $\Pi_v$ is of the form Ind $\tau_1|\det|^{s_1}\otimes\cdots\otimes\tau_\ell|\det|^{s_\ell}
\otimes\tau_{\ell+1}
\otimes\cdots\otimes\tau_{\ell+u}\otimes\tau_\ell|
\det|^{-s_\ell}\otimes\cdots\otimes
\tau_1|\det|^{-s_1}$, where the $\tau_i$'s are discrete series representations of smaller ${\rm GL}$'s, and $0 < s_l \leq\cdots\leq s_1 < {1\over 2}$
(cf.~[Tad]).

For $\sigma_v$ in the discrete series of ${\rm GL}_n(F_v),\ n=1,2,3,4$, the $L$-function $L(s,\Pi_v\times\sigma_v)$ is equal to
$$
\prod_{k=1}^\ell L(s-s_k,\tau_k\times\sigma_v)L(s+s_k,\tau_k\times\sigma_v)\prod^u_{j=1} L(s,\tau_{\ell+j}\times\sigma_v).
$$
By strict inequalities $0 < s_k < 1/2$ and the holomorphy of each $L(s,\tau_k\times\sigma_v)$ for ${\rm Re}(s) > 0$, it is easy to see that
as a function of $q_v^{-s}$, $L(s,\Pi_v\times\sigma_v)^{-1}$ has the same zeros as $\gamma(s,\Pi_v\times\sigma_v,\psi_v)$ and therefore the equalities
$$
L(s,\Pi_v\times\sigma_v)=L(s,\pi_{1v}\times\pi_{2v}\times\sigma_v)
$$
and
$$
\var(s,\Pi_v\times\sigma_v,\psi_v)=\var(s,\pi_{1v}\times\pi_{2v}\times\sigma_v,\psi_v),
$$
follow from
$$
\gamma(s,\Pi_v\times\sigma_v,\psi_v)=\gamma(s,\pi_{1v}\times\pi_{2v}\times\sigma_v,\psi_v),\tag5.4.1
$$
since $\pi_{iv}$ and $\sigma_v$ are tempered (cf.~Section 7 of [Sh1] as well as the beginning of Section 5 here).

It is therefore enough to prove (5.4.1) for $\sigma_v$ in the discrete series.
The case of an irreducible admissible generic $\sigma_v$ follows from this by multiplicativity.

To prove (5.4.1), we take a cuspidal representation $\sigma$ of ${\rm GL}_n(\bAB_F),\ n=1,2,3,4$, which has $\sigma_v$ as its $v^{\rm th}$ component [AC].
Although we know little about other components of $\sigma$, we still have
$$
\gamma(s,\Pi_w\times\sigma_w,\psi_w)=\gamma(s,\pi_{1w}\times\pi_{2w}\times\sigma_w,\psi_w)
$$
for every $w\not= v$; for $\pi_{iw}$ and $\Pi_w$ are now unramified, if $w\not= v,\ w < \infty$, and therefore $\gamma$-functions on both sides are 
the same products as those of Godement-Jacquet by multiplicativity (cf.~[Sh1]).
We recall that for $w=\infty$, the factors on both sides are those of Artin [Sh7], and are therefore automatically equal.
The equality (5.4.1) is now an immediate consequence of a comparison of the functional equations for $L(s,\Pi\times\sigma)$ and $L(s,\pi_1\times\pi_2\times\sigma)$.
We should point out that it would not be much harder to conclude the equality of root numbers and $L$-functions, if one establishes (5.4.1)
 only for supercuspidal~$\sigma_v$.
\enddemo

In the following proposition, let $K/k$ be a quadratic extension of local fields.
Let $\rho_1$ and $\rho_2$ be supercuspidal representations of ${\rm GL}_2(k)$ and ${\rm GL}_3(k)$, resp., and let $\Omega$ be the local lift of $\rho_1\otimes\rho_2$, constructed in Proposition 5.4.
Let $\rho_1^K$ and  $\rho_2^K$ be the base changes of $\rho_1$ and $\rho_2$, respectively, to $K$.
We shall assume that $\rho_1^K$ and $\rho_2^K$ are both supercuspidal, and let $\Omega_K$ be the local lift of $\rho_1^K\otimes\rho_2^K$.

\nonumproclaim{Proposition  5.5} Let $\Omega^K$ be the base change of $\Omega$ to
$K${\rm .} Then
$$\Omega_K=\Omega^K;
$$
i.e.{\rm ,} our local lift commutes with the base change{\rm .}
\endproclaim\pagebreak

\demo{Proof} By [Sh1, Prop.~5.1], we can take a number field $F$ with $F_v=k$, and $\pi_1=\hskip-6pt  \mathbold{\otimes}_w \pi_{1w},
\pi_2=\hskip-6pt  \mathbold{\otimes}_w \pi_{2w}$, cuspidal representations of ${\rm GL}_2(\Bbb A_F), {\rm GL}_3(\Bbb A_F)$, resp. such that
$\pi_{1v}=\rho_1,\pi_{2v}=\rho_2$ and $\pi_{1w},\pi_{2w}$ are unramified for $w\ne v$,\break $w<\infty$. Let $E/F$ be a quadratic extension such
that
$E_v=K$ ($v$ is inert). Let $\Sigma_v$ be a supercuspidal representation of ${\rm GL}_n(K)$, $n=1,2,3,4$, and\break $\Sigma=\hskip-6pt  \mathbold{\otimes}_u
\Sigma_u$, a cuspidal representation of
${\rm GL}_n(\Bbb A_E)$
such that $\Sigma_u$ is unramified for $u\ne v$, $u<\infty$.
By the definition of the local lift, we have
$$\gamma(s,\Sigma_v\times\Omega_K,\psi_v)=\gamma\left(s,\Sigma_v\times\rho_1^K\times\rho_2^K,\psi_v\right).
$$
We only need to prove that
$$\gamma(s,\Sigma_v\times\Omega^K,\psi_v)=\gamma\left(s,\Sigma_v\times\rho_1^K\times\rho_2^K,\psi_v\right).
$$

Since $\rho_1^K$ and $\rho_2^K$ are supercuspidal, $\pi_1^E$ and $\pi_2^E$ are cuspidal representations.
Let $\Pi$ be the weak lift constructed in the proof of Proposition 5.4
such that $\Pi_v=\Omega$.
Compare functional equations for $L(s,\Sigma\times\Pi^E)$ and
$L(s,\Sigma\times\pi_1^E\times\pi_2^E)$, and note that for $u\ne v$,
$$L(s,\Sigma_u\times (\Pi^E)_u)=L(s,\Sigma_u\times (\pi_1^E)_u\times (\pi_2^E)_u).
$$
Our equality now follows.
\enddemo

We continue with the same notation as before Proposition 5.4.

\nonumproclaim{Proposition  5.6} Suppose $v\in T${\rm .} Fix a normal closure $K$ of $L${\rm .}
Let $E/F_v$ be the unique quadratic extension of $F_v$ inside $K${\rm .}  Then
$$
L\big(s,\pi_{1v}^E\times \pi_{2v}^E\times\sigma\big)=L\big(s,\big(\pi_{1v}\boxtimes \pi_{2v}\big)^E\times\sigma\big)
$$
and
$$
\varepsilon\big(s,\pi_{1v}^E\times\pi_{2v}^E\times\sigma,\psi_{E/F_v}\big)=
\varepsilon\big(s,\big(\pi_{1v}\boxtimes\pi_{2v}\big)^E\times\sigma,\psi_{E/F_v}\big),
$$
for every irreducible admissible generic representation $\sigma$ of ${\rm GL}_n(E)${\rm ,}
$n=1,2,3,4${\rm .}  Here for each representation $\tau$ of ${\rm GL}_m(F_v)${\rm ,} $\tau^E$
denotes its base change to ${\rm GL}_m(E)${\rm ,} $m=2,3${\rm ,} and $\psi_{E/F_v}=\psi_v\cdot Tr_{E/F_v}${\rm .}
\endproclaim

\demo{Proof} Observe that $\pi_{iv}^E$, $i=1,2$, is still supercuspidal since
$v\in T$
(and moreover $(\pi_{1v}\boxtimes\pi_{2v})^E=\pi_{1v}^E\boxtimes\pi_{2v}^E)$.  Consequently, $\pi_{2v}^E$ is attached to $\text{Ind}(W_E, W_K,\eta)$, where $\eta$ is a
character of $W_K$ and $K/E$ is a cubic cyclic extension.  We can now apply
Proposition 5.3.b) to the pair $(\pi_{1v}^E,\pi_{2v}^E)$ to conclude the proof of the
proposition.  
\enddemo

\nonumproclaim{Corollary 5.7} Suppose $v\in T${\rm .}
 Then there exists a local lift $\Pi_v$ of $\pi_{1v}\otimes\pi_{2v}${\rm ,} such that $\Pi_v\simeq (\pi_{1v}\boxtimes\pi_{2v})\otimes\eta${\rm ,}
 where
$\eta^2=1${\rm .}
\endproclaim\pagebreak

\demo{Proof} By Propositions 5.4, 5.5, and 5.6, applied to $\pi_{1v}^E$ and   $\pi_{2v}^E$, and using the equality
$(\pi_{1v}\boxtimes \pi_{2v})^E=\pi_{1v}^E\boxtimes \pi^E_{2v}$,
one concludes that
$$
(\Pi_v)^E=(\pi_{1v}\boxtimes\pi_{2v})^E,
$$
by appealing to the local converse theorem proved in [CP-S1, \S 7], and [Ch].
Our result follows with $\eta$ whose kernel contains
$N_{E/F_v} (E^*)$. 
\enddemo

\nonumproclaim{Proposition  5.8} For $v\in T${\rm ,}
 $\Pi_v\simeq \pi_{1v}\boxtimes\pi_{2v}${\rm .} In particular{\rm ,} equalities {\rm (5.2.1)} and {\rm (5.2.2)}
 hold for $n=1,2,3,4$ at every place $v$ of $F${\rm .}
\endproclaim

\demo{Proof} We only need to apply Corollary 5.7 and Proposition 5.3.a) to the main theorem of the appendix [BH].
The proposition follows.
\enddemo

\demo{Proof of Theorem {\rm 5.1}} Let $\Pi=\hskip-6pt  \mathbold{\otimes}_v \Pi_v$, where $\Pi_v=\pi_{1v}\boxtimes\pi_{2v}$.
It is an irreducible admissible representation of ${\rm GL}_6(\Bbb A_F)$. Pick two finite places $v_1,v_2$, where $\pi_{jv_1},
\pi_{jv_2}$ are unramified for $j=1,2$.
Let $S_i=\{v_i\}$, $i=1,2$.
We apply the converse theorem twice
to $\Pi=\hskip-6pt  \mathbold{\otimes}_v \Pi_v$ with $S_1$ and $S_2$. We find two
automorphic representations $\Pi_1,\Pi_2$ of ${\rm GL}_6(\Bbb A_F)$ such
that $\Pi_{1v}\simeq\Pi_v$ for $v\ne v_1$, and $\Pi_{2v}\simeq\Pi_v$ for
$v\ne v_2$.
Hence $\Pi_{1v}\simeq \Pi_{2v}$ for all $v\ne v_1,v_2$.
By Proposition 3.10, $\Pi_1,\Pi_2$ are of the form
$\sigma_1\boxplus\cdots\boxplus\sigma_k$, where $\sigma_i$'s are
(unitary) cuspidal representations of ${\rm GL}$.
By the  classification theorem [JS1], $\Pi_1\simeq \Pi_2$, in particular,
$\Pi_{1v_i}\simeq \Pi_{2v_i}\simeq\Pi_{v_i}$ for $i=1,2$.
Thus $\Pi$ is automorphic.
\enddemo

\section{Functorial symmetric cubes for ${\rm GL}_2$}

Let ${\rm Sym}^m\colon {\rm GL}_2(\Bbb C)\longrightarrow {\rm GL}_{m+1}(\Bbb C)$ be the map given by
the $m^{\rm th}$ symmetric power
representation of ${\rm GL}_2(\bC)$ on the space of symmetric tensors of rank~$m$.
Let $\pi=\hskip-6pt  \mathbold{\otimes}_v\pi_v$ be a cuspidal representation of
${\rm GL}_2(\Bbb A_F)$ with central character $\omega_{\pi}$. By the local
Langlands correspondence, ${\rm Sym}^m(\pi_v)$ is well-defined for all $v$.
More precisely, it is the representation of ${\rm GL}_{m+1}(F_v)$ attached to ${\rm Sym}^m(\rho_v)$, where $\rho_v$ is the two-dimensional
representation of the  Deligne-Weil group attached to $\pi_v$. Hence Langlands' functoriality is equivalent to the assertion that
${\rm Sym}^m(\pi)=\hskip-6pt  \mathbold{\otimes}_v {\rm Sym}^m(\pi_v)$ is an automorphic representation
of ${\rm GL}_{m+1}(\Bbb A_F)$.
It is convenient to introduce
$A^m(\pi)={\rm Sym}^m(\pi)\otimes\omega_{\pi}^{-1}$ (called ${\rm Ad}^m(\pi)$ in [Sh3]).
If $m=2$, $A^2(\pi)={\rm Ad}(\pi)$.
Gelbart and Jacquet [GeJ] showed that ${\rm Ad}(\pi)$ is an automorphic representation of ${\rm GL}_3(\bAB_F)$, which is cuspidal unless $\pi$ is
monomial; i.e., $\pi\simeq\pi\otimes\eta$, where $\eta\not= 1$ is a gr\"ossencharacter of $F$. We prove:

\nonumproclaim{Theorem 6.1} The representation ${\rm Sym}^3(\pi)$ is an automorphic representation of ${\rm GL}_4(\Bbb A_F)${\rm .}
It is cuspidal{\rm ,} unless either $\pi$ is a monomial representation{\rm ,}
 or ${\rm Ad}(\pi)\simeq {\rm Ad}(\pi)\otimes\eta${\rm ,} for a nontrivial gr{\rm \"{\it o}}ssencharacter $\eta${\rm .}
Equivalently{\rm ,} ${\rm Sym}^3(\pi)$ is cuspidal{\rm ,} unless $\pi$ is either dihedral or
it has a cubic cyclic base change which is dihedral{\rm ,} i.e.{\rm ,} $\pi$ is of tetrahedral type{\rm .}
In particular{\rm ,}
 if $F=\bQ$ and $\pi$ is the automorphic cuspidal representation attached to a nondihedral holomorphic form of weight $\geq 2${\rm ,} then ${\rm
Sym}^3(\pi)$ is cuspidal{\rm .}
\endproclaim

We first consider:

\demo{{\rm 6.1.} $\pi$ is a monomial cuspidal representation}
That is, $\pi\otimes\eta\simeq\pi$ for a nontrivial gr\"ossencharacter $\eta$. Then $\eta^2=1$ and $\eta$ determines a quadratic extension $E/F$.
According to [LL], there is a gr\"ossencharacter $\chi$ of $E$ such that $\pi=\pi(\chi)$, where $\pi(\chi)$ is the automorphic representation whose
local factor at $v$ is the one attached to the representation of the local Weil group induced from $\chi_v$. Let $\chi'$ be the conjugate of $\chi$ by
the action of the nontrivial element of the Galois group. Then ${\rm Ad} (\pi)$ is given by
$${\rm Ad}(\pi)=\pi\big(\chi{\chi'}^{-1}\big)\boxplus\eta.$$
There are two cases:

Case 1. $\chi{\chi'}^{-1}$ factors through the norm; i.e., $\chi{\chi'}^{-1}=\mu\circ N_{E/F}$ for a gr\"ossen\-character $\mu$ of $F$. Then
$\pi(\chi{\chi'}^{-1})$ is not cuspidal. In fact, $\pi(\chi{\chi'}^{-1})=\mu\boxplus\mu\eta$. In this case,
$$A^3(\pi)=\pi\big(\chi{\chi'}^{-1}\big)\boxtimes\pi=\left(\mu\otimes\pi)\boxplus (\mu\eta\otimes\pi\right).
$$

Case 2. $\chi{\chi'}^{-1}$ does not factor through the norm. In this case, $\pi(\chi{\chi'}^{-1})$ is a cuspidal representation. Then

$$A^3(\pi)=\pi(\chi{\chi'}^{-1})\boxtimes\pi=\pi(\chi^2{\chi'}^{-1})\boxplus \pi.
$$
Here we used the fact that $\pi(\chi)_E=\chi\boxplus \chi'$ ([R1, Prop.~2.3.1]) and furthermore that $\pi'\boxtimes\pi=I_F^E (\pi'_E\otimes\chi)$, if
$\pi=\pi(\chi)$ ([R1,\S 3.1]). The index $E$ signifies the base change to $E$.
Observe that we are now using subscripts to denote the base change rather than superscripts   used in Section 5.
A superscript seemed to be a more appropriate notation for that section.
\enddemo

 6.2. $\pi$ {\it is not monomial}. Then ${\rm Ad}(\pi)$ is a cuspidal representation of ${\rm GL}_3(\Bbb A_F)$.
We first prove:

\nonumproclaim{Lemma 6.2} Let $\sigma$ be a cuspidal representation of
${\rm GL}_2(\Bbb A_F)${\rm .}
Then the triple $L$\/{\rm -}\/function
$L_S(s,{\rm Ad}(\pi)\times \pi\times\sigma)$
 has a pole at $s=1$ if and only if $\sigma\simeq\pi\otimes\chi$ and ${\rm Ad}(\pi)\simeq {\rm Ad}(\pi)\otimes (\omega_{\pi}\chi)$
for some gr{\rm \"{\it o}}ssencharacter $\chi${\rm .}
Here $S$ is a finite set of places for which $v\not\in S$ implies that both
$\pi_v$ and $\sigma_v$ are unramified{\rm .}
\endproclaim\pagebreak

\demo{Proof} Suppose $L_S(s,{\rm Ad}(\pi)\times \pi\times\sigma)$ has a pole at $s=1$. Consider $\pi\boxtimes\sigma$. It is an automorphic representation of ${\rm GL}_4(\Bbb A_F)$.
As argued at the beginning of the proof of Theorem 5.1, we may assume that $\sigma$ is not monomial.
By [R1, Th.\ M, p.~54] $\pi\boxtimes\sigma$ is then cuspidal, unless $\sigma\simeq\pi\otimes\chi$ for some gr\"ossencharacter $\chi$. If
$\pi\boxtimes\sigma$ is cuspidal, then $L_S(s,{\rm Ad}(\pi)\times \pi\times\sigma)$ is entire. Hence $\sigma\simeq\pi\otimes\chi$ for some
gr\"ossencharacter $\chi$. Consider the following $L$-function identity:
$$\multline
L_S(s,{\rm Ad}(\pi)\times\pi\times (\pi\otimes\chi))\\ =L_S(s,{\rm Ad}(\pi)\times ({\rm Ad}(\pi)\otimes (\omega_{\pi}\chi)))L_S(s,{\rm
Ad}(\pi)\otimes (\omega_{\pi}\chi)).\endmultline
$$
Since $L_S(s,{\rm Ad}(\pi)\otimes (\omega_{\pi}\chi))$ has no zero at $s=1$, $L_S(s,{\rm Ad}(\pi)\times ({\rm Ad}(\pi)\otimes (\omega_{\pi}\chi)))$ has a pole at $s=1$. Hence
${\rm Ad}(\pi)\simeq {\rm Ad}(\pi)\otimes (\omega_{\pi}\chi)$ since ${\rm Ad}(\pi)$ is self-contragredient.
The converse is clear from the above identity.  
\enddemo

We consider the functorial product $\pi\boxtimes {\rm Ad}(\pi)$ as in Theorem 5.1.
By Lemma 6.2 and the classification theorem [JS1],
$$
\pi\boxtimes {\rm Ad}(\pi)=\tau\boxplus\pi,
$$
where $\tau$ is an automorphic representation of ${\rm GL}_4(\Bbb A_F)$.
Since $\pi_v\boxtimes {\rm Ad}(\pi_v)=A^3(\pi_v)\boxplus \pi_v$,
we conclude $\tau_v\simeq A^3(\pi_v)$ for all $v$. Hence we have:

\nonumproclaim{Proposition  6.3} 
The representation $A^3(\pi)$ is an automorphic representation of ${\rm GL}_4(\Bbb A_F)${\rm .}
 It is not cuspidal if and only if there exists a nontrivial
gr{\rm \"{\it o}}ssencharacter $\eta$ such that ${\rm Ad}(\pi)\simeq {\rm Ad}(\pi)\otimes\eta${\rm .} In this case
$$
A^3(\pi)=(\pi\otimes\eta)\boxplus (\pi\otimes\eta^2).
$$
\endproclaim

\demo{Proof} We only need to prove the last assertion.
Clearly $$L_S(s,(\pi\boxtimes {\rm Ad}(\pi))\times\tilde\pi)$$ has a pole at $s=1$.
Thus $\pi\boxtimes {\rm Ad}(\pi)=\pi\boxplus\tau$, where $\tau$ is an automorphic representation of ${\rm GL}_4(\bAB_F)$.
If $\tau$ is not cuspidal, then $\tau=\sigma_1\boxplus\sigma_2$, where $\sigma_i,\ i=1,2$, are cuspidal representations of ${\rm GL}_2(\bAB_F)$.
Then $L_S(s,(\pi\boxtimes {\rm Ad}(\pi))\times\tilde\sigma_i)$ must have a pole at $s=1$.
We can now proceed as in either Lemma 6.2 or as in the proof of the last statement of Theorem 5.1, to conclude that $\sigma_i=\pi\otimes\eta^i$, where ${\rm Ad}(\pi)\cong {\rm Ad}(\pi)\otimes\eta,\ \eta\not= 1$.
\enddemo

\nonumproclaim{{C}orollary 6.4} The representation ${\rm Sym}^3(\pi)$ is an automorphic representation of ${\rm GL}_4(\bAB_F)${\rm .}
It is not cuspidal if and only if there exists a nontrivial gr{\rm \"{\it o}}ssencharacter 
$\eta$ such that ${\rm Sym}^2(\pi)\cong$ ${\rm Sym}^2(\pi)\otimes\eta${\rm .}
In this case
$$
{\rm Sym}^3(\pi)=(\pi\otimes\eta\omega_\pi)\boxplus (\pi\otimes\eta^2\omega_\pi).
$$
\endproclaim

To complete the proof of Theorem 6.1, we need

\nonumproclaim{Lemma 6.5} Suppose $\pi$ is not monomial{\rm ,} but ${\rm Ad}(\pi)$ is{\rm ,} i.e.{\rm ,} there exists a nontrivial 
gr{\rm \"{\it o}}ssencharacter $\eta${\rm ,} necessarily cubic{\rm ,} such that ${\rm Ad}(\pi)\otimes\eta\cong$ ${\rm Ad}(\pi)${\rm .}
Then $\pi$ is of tetrahedral type{\rm ,} i.e.{\rm ,} there exists a Galois representation $\sigma$ of tetrahedral type such that $\pi=\pi(\sigma)${\rm .}
\endproclaim

\demo{Proof} Let $E/F$ be the cubic cyclic extension defined by $\eta$.
Then as observed in [Sh8], the $E/F$-base change $\pi_E$ of $\pi$ is monomial.
Let $\sigma_E$ be the two-dimensional dihedral representation of $W_E$
attached to $\pi_E$.
Since $\sigma_E$ is invariant under the action of Gal$(E/F)$, it extends to a two-dimensional representation $\sigma$ of $W_F$ which is now of
tetrahedral type (cf.~[La2], [Ge]). Let $\pi'$ be the cuspidal representation of ${\rm GL}_2(\bAB_F)$ attached to $\sigma$, i.e.~$\pi'=\pi(\sigma)$.
Observe that $\pi'_{E}\cong\pi_E$ as they both correspond to $\sigma_E$ and therefore $\pi'\cong\pi\otimes\eta^a$ for some $a=0,1,2$.
But the lift $\sigma$ is unique only up to twisting by a power of $\eta$, and therefore changing the choice of $\sigma$ if necessary, 
we have $\pi\cong\pi'=\pi(\sigma)$.

The last assertion follows from the fact that holomorphic forms of weight $\geq 2$ can never be of tetrahedral type. 
\enddemo

\demo{{R}emark {\rm 6.6}} For the proof of the functoriality of $\pi\boxtimes {\rm Ad}(\pi)$, and hence ${\rm Sym}^3(\pi)$, we do not need the
appendix [BH]. The appendix is needed only for the general case of functoriality for ${\rm GL}_2\times {\rm GL}_3$.
The reason is the following.
By Proposition 5.4, we can construct a strong lift $\Pi=\hskip-6pt  \mathbold{\otimes}_v \Pi_v$ of $\pi\otimes {\rm Ad}(\pi)$.
Then we still have $\Pi=\tau\boxplus\pi$.
We only need to prove that $\tau_v\simeq A^3(\pi_v)$ for all $v$.
For that, it is enough to prove $\gamma(s,\sigma_v\times\tau_v,\psi_v)=\gamma(s,\sigma_v\times A^3(\pi_v),\psi_v)$ for every irreducible generic representation $\sigma_v$ of ${\rm GL}_n(F_v)$, $n=1,2,3$.
Let $\varphi_v$ be the corresponding Weil group representation attached to $\pi_v$.
Then $\varphi_v\otimes {\rm Ad}(\varphi_v)\simeq A^3(\varphi_v)\oplus\varphi_v$.
Now the equality of $\gamma$-factors follows immediately from (5.3.1) and (5.1.4).
We should remark that since we are able to twist up to $n=3$, the local converse theorem proved in [Ch] and [CP-S1] is not necessary in this case.
We should point out that in proving Theorem 5.1, we could have confined ourselves to more conventional converse theorems, both local and global, if we were also to use case $E_8-1$ of [Sh2], which would allow us to twist by ${\rm GL}_5$ as well.
In view of [CP-S1], we did not pursue this.
We should remind the reader that in our present approach, we have used the local converse theorem of [Ch] and [CP-S1] only in the proof of Corollary
5.7.
\enddemo

\section{New estimates towards the Ramanujan and Selberg conjectures}

An immediate consequence of the existence of ${\rm Sym}^3(\pi)$ is a new estimate on Hecke eigenvalues of a Maass form.
More precisely, let $\pi$ be a cuspidal representation of ${\rm GL}_2(\bAB_F)$.
Write $\pi=\hskip-6pt  \mathbold{\otimes}_v\pi_v$.
Assume $\pi_v$ is spherical, i.e.,
$$
\pi_v=\Ind(|\ \ |_v^{s_{1v}},\ |\ \ |_v^{s_{2v}}),
$$
$s_{iv}\in\bC$.
When $v < \infty$, we set $\alpha_{1v}=|\varpi_v|_v^{s_{1v}}$ and $\alpha_{2v}=|\varpi_v|_v^{s_{2v}}$, so that diag~$(\alpha_{1v},\alpha_{2v})$
represents the semisimple conjugacy class in ${\rm GL}_2(\bC)$ attached to $\pi_v$. Next, suppose $F=\bQ$.
Let $\pi$ be attached to a Maass form $f$ with respect to a congruence subgroup $\Gamma$.
Denote by $\lambda_1=\lambda_1(\Gamma)$ the smallest positive eigenvalue of Laplace operator $\Delta=-y^2\big({\partial^2\over  \partial x^2}+{
\partial^2\over  \partial y^2}\big)$ on $L^2(\Gamma\backslash\frakh)$, where $\frakh$ denotes the upper half plane. It depends on $f$.
More precisely,
$$
\Delta f={1\over 4}(1-s^2) f,
$$
with $\lambda_1(\Gamma)={1\over 4}(1-s^2)$, where $s=2 {\rm Re}(s_{1\infty})=-2 {\rm Re}(s_{2\infty})$.
Then, one expects

\vglue5pt {\elevensc Selberg's Conjecture}.  $\lambda_1(\Gamma)\geq 1/4${\rm .}

\vglue5pt {\elevensc Ramanujan-Petersson's Conjecture}. $|\alpha_{1v}|=|\alpha_{2v}|=1$.
\vglue5pt

We now prove

\nonumproclaim{Theorem 7.1} Let $\pi$ be a cuspidal representation of ${\rm GL}_2(\Bbb A_F)${\rm .}
 Let $\pi_v$ be a
local {\rm (}\/finite or infinite\/{\rm )} spherical component of $\pi${\rm ,} i.e.~$\pi_v=\Ind  (|\ |_v^{s_{1v}}, |\ |_v^{s_{2v}})${\rm .}
 Then $|{\rm Re}(s_{jv})|\leq
\frac 5{34},\ j=1,2${\rm .}
\endproclaim

\demo{Proof} Consider ${\rm Sym}^3(\pi)$. If it is not cuspidal, then $\pi$, being of dihedral or tetrahedral type, satisfies Ramanujan's conjecture.
Hence we may assume that it is cuspidal. We apply [LRS1] to ${\rm Sym}^3(\pi)$.
It states that if $\Pi=\hskip-6pt  \mathbold{\otimes}_v\Pi_v$ is a cuspidal representation of ${\rm GL}_n(\bAB_F)$ and if $\Pi_v$ is the spherical constituent of $Ind^{{\rm GL}_n(F_v)}_{B_n(F_v)}\ds\otimes^n_{j=1} |\ \ |_v^{t_{jv}},\ t_{jv}\in\bC,\ j=1,\ldots,n$, then ${\rm Re}(t_{jv})\leq \ds{1\over 2}-\ds{1\over n^2+1}$.
In our case $n=4$ and 
$$
3{\rm Re}(s_{iv})\leq \frac 12-\frac 1{4^2+1}.
$$
Our result follows. 
\enddemo

{\it Remark} {\rm 7.2}. Theorem 7.1 establishes a new estimate towards   Selberg's conjecture.
It is $\lambda_1={1\over 4}(1-s^2)\geq 66/289\cong 0.22837$.
We refer to Luo-Rudnick-Sarnak [LRS2] for the estimate $\lambda_1\geq 0.21$ obtained earlier.
When $v < \infty$, Theorem 7.1 can be written as
$$
|\alpha_{jv}| \leq q_v^{5/34}.
$$
The earlier exponents of $1/5$ and $5/28$ for arbitrary $F$ and $F=\bQ$ are due to [Sh2] and [LRS1], and [BDHI], respectively.

\section{Applications to analytic number theory}

This section has been suggested to us by Peter Sarnak.
We would like to thank him for his suggestion and helpful advice.

In this section let $F=\bQ$.
The bounds towards the Ramanujan conjectures
$$
|\alpha_{jp}|\leq p^{5\over 34}\tag8.1
$$
and especially their archimedean counterparts
$$
|{\rm Re}(s_{j\infty})|\leq {5\over 34}\tag8.2
$$
have numerous applications, see for example [I] and [IS].
They lead to improvements of exponents in a number of places, namely those where the question of small eigenvalues of the Laplacian on $\Gamma_0(N)\backslash\frak h$ enters (see [Se]).
One such example is that of cancellation in Kloosterman sums (see [Sa]).

There are also applications where the fact that (8.1) and (8.2) are sharper than ${1\over 6}$ has more fundamental consequences.
Let $\Gamma$ be a congruence subgroup of ${\rm SL}_2(\bZ)$.
Let $f$ be either a holomorphic cusp form of weight $k$, or a Maass cusp form for $\Gamma\backslash \frak h$.
In the latter case, $f$ is an eigenfunction of the Laplacian;
$$
y^2 \bigg({\partial^2\over\partial x^2}+{\partial^2\over\partial y^2}\bigg) f(z)+\bigg({1\over 4}-t^2\bigg) f(z)=0,\tag8.3
$$
with $t$ imaginary or $0 < t < {1\over 2}$ and $z=x+iy,\ y > 0$.
More precisely, we take $t=s_{1\infty}$.
In both cases, $f$ has a Fourier expansion
$$
f(z)=\sum_{n\not= 0} a_f (n)|n|^{k-1\over 2} W (nz),\tag8.4
$$
where $a(1)=1$.
Here if $f$ is a holomorphic cusp form, $W(z)=e^{2\pi iz}, n > 0$ and $k$ is the weight.
If $f$ is a Maass form, $k=0$ and
$$
W(z)=\sqrt{y} K_t (2\pi y) e^{2\pi ix},
$$
where $K_t$ is a Bessel function and $W(z)=W(\overline z)$, if $z$ is in the lower half plane.
With this normalization, the Ramanujan conjectures for $f$ (at the finite places) are equivalent to
$$
a_f(n)=O_\varepsilon (|n|^\varepsilon),\ \ \varepsilon > 0.\tag8.5
$$

The problem of cancellations in sums of shifted coefficients is as follows.
Fix $h\in\bZ,\ h\not= 0$, and for $X$ large, set
$$
D_{f,h}(X)\colon = \sum_{|n|\leq X} a_f (n) a_f (n+h).\tag8.6
$$
This sum has been studied extensively ([Go], [M], to name a few).
(We may include here the case of shifted divisor function sums.)
It is known [Go] that if the exceptional eigenvalues of the Laplacian on $\Gamma\backslash \frak h$ are denoted by $t_j$,\break $0 < t_j < {1\over
2}$, then as $X\to\infty$,
$$
\align
D_{f,h}(X)&=\sum_{0\leq t_j < {1\over 2}} b_{f,j} X^{{1\over 2}+t_j}+O(X^{2\over 3})\tag8.7\\
&=\sum_{{1\over 6} < t_j < {1\over 2}} b_{f,j} X^{{1\over 2}+t_j}+O(X^{2\over 3})\tag8.8
\endalign
$$
for suitable constants $b_{f,j}$ depending on $f$ and $t_j$.

When $f$ is holomorphic, equation (8.7) is proved in [Go].
(See Part i) of Theorem 1 of [Go]; with our normalization of Fourier coefficients in (8.4), $O(X^{k-1/3})$ in [Go] becomes $O(X^{2/3})$.)
One expects similar arguments to apply, and (8.7) must be valid, even if $f$ is a Maass form. But, so far as we know, no reference to that is available.

The error term of $X^{2\over 3}$ above is present because of real analysis issues involved with the form of the sharp cutoff in the sum defining $D_{f,h}(X)$ (i.e., the sum gives weight 1 for $n\leq X$ and weight 0 for $n > X$).
Insomuch as the main results, (8.1) and (8.2) above,  yield $|t_j|\leq {1\over 6}$, we have:

\nonumproclaim{Proposition  8.1} Let $\Gamma$ and $f$ be as above{\rm .}
Suppose $f$ is a holomorphic cusp form{\rm .} Fix $h\not= 0${\rm .} Then as $X\to\infty${\rm ,}
$$
D_{f,h}(X)=O(X^{2\over 3}).
$$
\endproclaim

Another application is to that of the  hyperbolic circle problem [I, p.~190].
Let $z,w\in \frak h$. The hyperbolic distance is
$$
\rho(z,w)=\log{ |z-\overline w|+|z-w|\over |z-\overline w|-|z-w| }.
$$
Let
$$
2\cosh \rho(z,w)=e^{\rho(z,w)}+e^{-\rho(z,w)}=2+4u(z,w),\quad
u(z,w)={|z-w|^2\over 4 {\rm Im}\, z {\rm Im}\, w}.
$$
For $X\geq 2$, set
$$
P(X)=\# \{\gamma\in\Gamma\ |\ 4u(\gamma z,w)+2\leq X\}.
$$

\nonumproclaim{Proposition  8.2} Let $|F|$ be the volume of the fundamental domain of $\Gamma\backslash\frak h${\rm .}
Then
$$
P(X)=\pi |F|^{-1} X+O(X^{2\over 3}).
$$
\endproclaim

This follows from [I, Th.\ 12.1], and again from the fact that there are no eigenvalues between ${1\over 6}$ and ${1\over 2}$. To be precise, note that
in the notation of [I], the eigenvalue $\lambda_j=s_j (1-s_j)={1\over 4}-t^2_j=2/9$ corresponds to $t_j=1/6$ and $s_j=2/3$, and therefore in the
equation (12.6) of [I] we may disregard the interval $2/3 \leq s_j < 1$.

Thus, the point is that with our present understanding of the (sharp cutoff) sums (8.6), such results yield as sharp a result as one 
would get if one assumes the full Ramanujan conjectures $({\rm Re}(t_j)=0)$.

\section{Siegel cusp forms of weight three}

In this section, using the existence of symmetric cubes which we proved
in Section 6 and a conjecture of Arthur [A2], we prove the existence of
illusive Siegel cusp forms of weight $3$.
We  thank Joseph Shalika for suggesting the problem and start with the following unpublished result:

\proclaimtitle{Jacquet, Piatetski-Shapiro, and Shalika}
\specialnumber{9.1}
\proclaim{Theorem}  
Let $\sigma=\hskip-6pt  \mathbold{\otimes}_v\sigma_v$ be a cuspidal automorphic representation of ${\rm GL}_4 (\bAB_F)$
 for which there exists a gr{\rm \"{\it o}}ssencharacter $\chi$ and a finite set of places $S$ for which every $\sigma_v$
 with $v\not\in S$ is unramified{\rm ,} such
that $L_S(s,\sigma,\Lambda^2\otimes\chi^{-1})$ has a pole at $s=1${\rm .}
 Then there exists a globally generic cuspidal automorphic representations
$\tau$ of ${\rm GSp}_4(\bAB_F)$ with central character $\chi$ such that $\sigma$ is the functorial lift of $\tau$ under the embedding ${\rm
GSp}_4(\bC)\hookrightarrow {\rm GL}_4(\bC)${\rm .}
\endproclaim

Let $\pi=\hskip-6pt  \mathbold{\otimes}_v\pi_v$ be a cuspidal representation of ${\rm GL}_2(\bAB_F)$ whose central character is $\omega_\pi$.
Assume $\pi$ is neither of dihedral, nor of tetrahedral type.
Then ${\rm Sym}^3(\pi)$ is a cuspidal representation of ${\rm GL}_4(\bAB_F)$.
It is easy to see that if $S$ is a finite set of places outside of which every $\pi_v$ is unramified, $L_S(s, {\rm Sym}^3(\pi),\ \Lambda^2\otimes \omega_\pi^{-3})$ has a pole at $s=1$.
Consider the identity $L_S(s, {\rm Ad}(\pi)\times {\rm Ad}(\pi))=L_S(s, {\rm Ad}(\pi)) L_S (s, A^3(\pi), \Lambda^2\otimes \omega_\pi^{-1})$.
Since ${\rm Ad}(\pi)$ is self-contragredient, the left-hand side has a pole at $s=1$.
But $L_S(s, {\rm Ad}(\pi))$ has no zeros at $s=1$ and consequently $L_S(s, A^3(\pi), \Lambda^2\otimes\omega_\pi^{-1})$ has a pole at $s=1$.
Now observe that $L_S(s, A^3 (\pi), \Lambda^2\otimes\omega_\pi^{-1})=L_S(s, {\rm Sym}^3 (\pi), \Lambda^2\otimes\omega_\pi^{-3})$.
Applying Theorem 9.1, we may consider ${\rm Sym}^3(\pi)$ as a cuspidal representation of ${\rm GSp}_4(\bAB_F)$.

Now, assume $F=\bQ$, and let $\pi$ correspond to a non-CM type holomorphic cusp form of weight 2. Then $\pi_\infty$ is a (holomorphic) discrete series parametrized by 
a two-dimensional representation $\varphi_\infty$ of the Weil group $W_\bR$. We assume $\omega_\pi$ is trivial. Then $\varphi_\infty|\bC^*$ is
given by
$$
z\mapsto \pmatrix z/|z|&0\\ 0&\overline z/|z|\endpmatrix.
$$
It is now clear that the $L$-packet of ${\rm Sym}^3(\pi_\infty)$ is parametrized by the homomorphism $\psi_\infty\colon W_\bR\to {\rm
GSp}_4(\bC)$ for which $\psi_\infty|\bC^*$ is given by
$$
z\mapsto \text{ diag}(z^3/|z|^3, z/|z|,\ \overline z/|z|, \overline z^3/|z|^3).
$$

Now, let $L_\bQ$ be the conjectural Langlands group whose two-dimensional representations parametrize automorphic forms on ${\rm
GL}_2(\bAB_\bQ)$. Let $\varphi\colon L_\bQ\to {\rm GL}_2(\bC)$ be the two-dimensional representation of $L_\bQ$ which parametrizes $\pi$. Then
$\psi={\rm Sym}^3(\varphi)={\rm Sym}^3\cdot \varphi$ factors through ${\rm GSp}_4(\bC)$ and $\psi|W_\bR=\psi_\infty$. Since $\pi$ is
nondihedral, $\varphi$ is surjective, and therefore the image of ${\rm Sym}^3(\varphi)$ in ${\rm GSp}_4(\bC)$ is the same as that of ${\rm Sym}^3$.
Consequently, by Schur's lemma, the centralizer of ${\rm Im}({\rm Sym}^3(\varphi))$ in ${\rm GSp}_4(\bC)$ consists of only scalars $\bC^*$. Thus,
since the centralizer is connected, the component group $S_\varphi$ of Arthur [A2]  is trivial. Now, his multiplicity formula in [A2], implies that
every member of the\break $L$-packet of ${\rm Sym}^3(\pi)$, i.e., the packet parametrized by $\psi$, is automorphic. In particular, we can change the
component
$\pi_\infty$ in its $L$-packet from the generic to the holomorphic one [V], so as to pick up a  weight $3$ Siegel modular cusp form. Thus we have
proved:

\nonumproclaim{Theorem 9.2} Let $\pi$ be a cuspidal representation of ${\rm GL}_2(\bAB_\bQ)$ attached to a non\/{\rm -CM} holomorphic form of
weight {\rm 2.} Assume the validity of Arthur\/{\rm '}\/s multiplicity formula for ${\rm GSp}_4(\bAB_\bQ)${\rm .}
Then every member of the $L$-packet of ${\rm Sym}^3(\pi)$ is automorphic.
In particular{\rm ,} Siegel modular cusp forms of weight $3$ exist{\rm .}
\endproclaim

\vglue-6pt
{\it Remark} {\rm 9.3}. While the group $L_\bQ$ is out of reach at present, one expects to prove Arthur's multiplicity formulas using his trace
formula. One can then also expect to define the group $S_\varphi$, using the trace formula and necessary local results (Langlands, Shelstad), without
any knowledge of $L_\bQ$ and $\varphi$. For ${\rm GSp}_4$, the stable trace formula is in very good shape, since all the fundamental lemmas are
already established. On the other hand, to be able to use Arthur's trace formula so as to prove Theorem 9.2, one actually needs to show that the
transfer of $\sym^3(\pi)$ to ${\rm GSp}_4(\bAB_\bQ)$ belongs to a stable $L$-packet. In general, this will require a comparison of the regular trace
formula for ${\rm GSp}_4(\bAB_F)$ with the most general type of twisted trace formula for ${\rm GL}_4(\bAB_F)$, i.e., the one which picks up $\Pi$
satisfying $\Pi(^tg^{-1})\cong \Pi(g)\omega(\det g)$ for a gr\"ossencharacter $\omega$. (As was pointed out to us by Ramakrishnan, this may be
accomplished by considering ${\rm GL}_4\times {\rm GL}_1$.) This does not seem to be as in good shape as the stable trace formula for ${\rm
GSp}_4(\bAB_F)$   alluded to before.
 
\vglue-6pt
\section{Applications to global Langlands correspondence\\ and Artin's conjecture}
\vglue-6pt

In view of recent consequences of Taylor's program [Tay] on Artin's conjecture (Buzzard, Dickinson, Shepherd-Barron, Taylor), and the work of
Langlands [La2] and Tunnell [Tu], with few exceptions, to every two-dimensional continuous representation $\sigma$ of $W_F$, the Weil group of
$\overline F/F$, one can attach an automorphic representation of ${\rm GL}_2(\bAB_F)$, preserving root numbers and\break $L$-functions. We
remark that, if $\sigma$ is of icosahedral type, we need to assume $F=\bQ$ and $\sigma$ is odd [Tay].

Given a two-dimensional irreducible continuous representation $\sigma$ of $W_F$, let $\pi(\sigma)$ be the corresponding cuspidal representation
of
${\rm GL}_2(\bAB_F)$. On the other hand, ${\rm Sym}^3(\sigma)$ defines a continuous four-dimensional representation of $W_F$.
By Theorem 6.1, ${\rm Sym}^3(\pi(\sigma))$ is an automorphic representation of\break ${\rm GL}_4(\bAB_F)$, and the map ${\rm
Sym}^3(\sigma)\mapsto {\rm Sym}^3(\pi(\sigma))$ preserves root numbers and $L$-functions. In fact, if $\tau$ is another one- or two-dimensional
irreducible continuous representation of $W_F$ for which $\pi(\tau)$ exists, then this correspondence preserves the factors for pairs ${\rm
Sym}^3(\sigma)\otimes\tau$ and ${\rm Sym}^3(\pi(\sigma))\times\pi(\tau)$. We record this as:

\nonumproclaim{Theorem 10.1} Let $\sigma$ be a continuous irreducible two\/{\rm -}\/dimensional representation of $W_F$ for which $\pi(\sigma)$ exists
{\rm (}\/see discussion above\/{\rm ).} Then{\rm ,} to the continuous four\/{\rm -}\/dimensional representation ${\rm Sym}^3(\sigma)$ of $W_F${\rm ,}
 there is attached an automorphic representation of ${\rm GL}_4(\bAB_F)${\rm ,} the symmetric cube ${\rm Sym}^3(\pi(\sigma))$
 of $\pi(\sigma)${\rm ,} which
preserves root numbers and $L$\/{\rm -}\/functions{\rm .} In particular{\rm ,}
 if $\sigma$ is not of dihedral type{\rm ,} then the Artin $L$\/{\rm -}\/function
$L(s,{\rm Sym}^3(\sigma)\otimes\tau)$ is entire for every one\/{\rm -}\/dimensional representation $\tau$ of $W_F${\rm .}
Similarly{\rm ,} if $\sigma$ is  of tetrahedral type{\rm ,}
 then the same statement is also true for any two\/{\rm -}\/dimensional irreducible continuous representation $\tau$ of $W_F${\rm ,}
 provided that $\pi(\sigma)$
exists{\rm .} Moreover{\rm ,} if $\sigma$ is a representation of icosahedral type for which $\pi(\sigma)$ exists{\rm ,}
 then $L(s,{\rm Sym}^3(\sigma))$ is a degree\/{\rm -}\/four irreducible primitive Artin $L$\/{\rm -}\/function which is entire{\rm .}
\endproclaim

It is easy to see, by the material in Section 6, that even when $\sigma$ is of octahedral type, ${\rm Sym}^3(\sigma)$, although irreducible, is not
primitive. In fact, there exists a quadratic extension $E/F$ such that $\sigma_E=\sigma|W_E$ is of tetrahedral type, to which there is attached a
cubic character $\eta$. Then by Corollary 6.4,
$$
{\rm Sym}^3(\sigma)|W_E={\rm Sym}^3(\sigma_E)=(\sigma_E\otimes\eta \det(\sigma_E))\oplus \left(\sigma_E\otimes\eta^2\det(\sigma_E)\right).
$$
This immediately implies that ${\rm Sym}^3(\sigma)$ is induced.
But those attached to icosahedral type representations will always be primitive since $A_5$ has no proper subgroup of index less than 5.

Next, let $\sigma$ and $\tau$ be two continuous two-dimensional irreducible representations of $W_F$ such that $\pi(\sigma)$ and $\pi(\tau)$
exist. Then by Theorem 5.1, the representation $\sigma\otimes {\rm Sym}^2(\tau)$ corresponds to $\pi(\sigma)\boxtimes {\rm Sym}^2(\pi(\tau))$,
i.e., to the six-dimensional representation $\sigma\otimes {\rm Sym}^2(\tau)$ of $W_F$, there is attached an automorphic representation of ${\rm
GL}_6(\bAB_F)$, preserving root numbers and $L$-functions.

Suppose $\eta$ is another two-dimensional continuous representation of $W_F$ which is not of dihedral type. We will assume that $\sigma$ and
$\eta$, have nonconjugate images in $P {\rm GL}_2(\bC)$. Then by [R1], $\pi(\sigma)\boxtimes \pi(\eta)$, is a cuspidal representation of ${\rm
GL}_4(\bAB_F)$.  Thus
$$
L(s,\pi(\sigma)\times {\rm Sym}^2 (\pi (\tau))\times \pi(\eta))\tag10.1
$$
is entire. Next assume $\sigma$ and $\eta$ have conjugate images in $P {\rm GL}_2(\bC)$. One can check quickly that by [JS1], (10.1) has a pole if
${\rm Sym}^2(\pi(\tilde\tau))$ is a twist of ${\rm Sym}^2(\pi(\sigma))$, or equivalently ${\rm Ad}(\pi(\tau))$ is a twist of ${\rm Ad}(\pi(\sigma))$.
We therefore see that (10.1) is entire if ${\rm Ad}(\tau)$ and ${\rm Ad}(\sigma)$ have nonconjugate images in ${\rm GL}_3(\bC)$. We state this as:

\nonumproclaim{Theorem 10.2} {\rm a)}
 Let $\sigma$ and $\tau$ be two continuous two\/{\rm -}\/dimensional representations of $W_F$
 for which $\pi(\sigma)$ and $\pi(\tau)$ exist{\rm .} Then to the six\/{\rm -}\/dimensional representation $\sigma\otimes {\rm Sym}^2(\tau)${\rm ,}
 there is attached an automorphic representation of ${\rm GL}_6(\bAB_F)${\rm ,} which
preserves root number and $L$\/{\rm -}\/functions for pairs{\rm .}

{\rm b)} Let $\eta$ be another such representation and assume $\sigma$ and $\eta$ have nonconjugate images in $P {\rm GL}_2(\bC)${\rm .}
Then $L(s,\sigma\otimes {\rm Sym}^2(\tau)\otimes\eta)$ is entire{\rm .}

{\rm c)} Assume $\sigma$ and $\eta$ have conjugate images in $P {\rm GL}_2(\bC)$ and neither $\sigma${\rm ,} nor~$\tau${\rm ,}
 nor $\eta$ is dihedral.
Moreover assume ${\rm Ad}(\sigma)$ and ${\rm Ad}(\tau)$ have nonconjugate images in $P {\rm GL}_3(\bC)${\rm .}
Then again $L(s,\sigma\otimes {\rm Sym}^2 (\tau)\otimes\eta)$ is entire{\rm .}
\endproclaim

{\it Remark}. We note that Theorem 10.2 is only new if either $\sigma$ or $\tau$ is of icosahedral type.
In fact, Theorem 10.2 does not give anything new beyond [R3] in the solvable cases.
\input Shahidi.ref

\font\scr=eusm10 
\let\le\leqslant
\let\ge\geqslant
\define\dual#1{\overset{{}_\vee} \to {#1}}
\define\ado{\text{\rm Ad}^0}

\let\ups\upsilon
\let\ve\varepsilon

\define\GL#1#2{\text{\rm GL}_{#1}({#2})}

\vglue-36pt
\title{Appendix:\\
On certain dyadic representations}
\twoauthors{Colin J.\ Bushnell}{Guy Henniart}
\institutions{King's College,  London, United Kingdom\\
{\eightpoint {\it E-mail address\/}: bushnell\@mth.kcl.ac.uk}
\\
\vglue6pt
Universit\'e de Paris-Sud, Orsay, France\\
{\eightpoint {\it E-mail address\/}:  Guy.Henniart\@math.u-psud.fr}}
\acknowledgements{This work was supported in part by the European network TMR ``Arithmetical Algebraic Geometry''. The first-named author also
thanks l'Universit\'e de Paris-Sud and l'Institut Henri Poincar\'e for hospitality during the conception of this work.
}

In this appendix, $F$ denotes a dyadic local field, that is, $F$ is a finite extension field of the 2-adic rational field $\Bbb Q_2$. We fix an algebraic closure $\overline F$ of $F$, and write 
$\hbox{\scr W}_F$ for the Weil group of $\overline F/F$. If $E/F$ is a finite extension (always assumed to be contained in $\overline F$) we denote by 
$\hbox{\scr W}_E$ the Weil group of $\overline F/E$. Also, if $\tau$ is a finite-dimensional continuous semisimple representation of 
$\hbox{\scr W}_E$, we write
$\Ind_{E/F}(\tau)$ for the representation of $\hbox{\scr W}_F$ induced by $\tau$.
 
Let $\psi$ be a nontrivial additive character of $F$. Let $n_1$, $n_2$ be positive integers and let $\pi_1$, $\pi_2$ be irreducible smooth representations of $\GL{n_1}F$, $\GL{n_2}F$ respectively. We form
$$
\gamma(\pi_1\times \pi_2,s,\psi) = \ve(\pi_1\times \pi_2,s,\psi)\, \frac{L(\dual\pi_1 \times \dual \pi_2,1{-}s)}{L(\pi_1\times\pi_2,s)}\,,
$$
as in \cite{6}, \cite{8}, where $s$ is a complex variable.
 
We now state our main result. This makes no pretence at generality: it is designed to serve only the paper \cite{7} to which this is an appendix. We need the following data:
\vglue4pt
\item{(1)} {\it An irreducible continuous representation $\rho_1$ of $\hbox{\scr W}_F${\rm ,} which is\/ {\rm primitive} and of dimension $2${\rm ,} and
\item{\rm (2)}
an irreducible continuous representation $\rho_2$ of\/ $\hbox{\scr W}_F${\rm ,} of dimension $3$ and satisfying\/{\rm :}\/
\itemitem{\rm (a)} There are a noncyclic cubic extension $L/F$ and a quasicharacter $\theta$ of\/ $\hbox{\scr W}_L$
 such that $\rho_2 \cong \Ind_{L/F}(\theta)${\rm ,} and
\itemitem{\rm (b)} $\rho_2\not\cong \chi\otimes \rho_2$ for any unramified quasicharacter $\chi \neq 1$ of\/ $\hbox{\scr W}_F${\rm .}}
 \vglue4pt 

We put $\rho= \rho_1\otimes \rho_2$, and we let $\pi$ be the irreducible smooth representation of $\GL6F$ corresponding to $\rho$ via the Langlands correspondence of \cite{4}, \cite{5}.
\nonumproclaim{Main Theorem}\hskip-4pt
Let $\pi'$ be an irreducible smooth representation of $\GL6F$ satisfying the following conditions\/{\rm :}
\vglue3pt
\item{\rm (1)}
Let $E/F$ be an unramified quadratic extension{\rm ,}
 and let $\pi_E$, $\pi'_E$ be the representations of $\GL6E$ obtained by base change from $\pi${\rm ,} $\pi'$ respectively{\rm .} Then
$$
\pi_E \cong \pi'_E.
$$
\item{\rm (2)}
$
\gamma(\pi\times\sigma,s,\psi) = \gamma(\pi'\times\sigma,s,\psi)
$ 
for every irreducible supercuspidal representation $\sigma$ of $\GL mF${\rm ,} $m=1,2,3${\rm .}
\endproclaim

Now
$$
\pi \cong \pi'.
$$

We start the paper with a discussion of the structure of $\rho$ in Section~1.
Then, in Section~2, we give a consequence of the conductor formula of \cite{2}
applicable to the present situation. We can then prove the main theorem
in Section~3.
\advance \sectioncount -10 
\section{Galois representations}

In this section, we are only concerned with the representations $\rho_1$, $\rho_2$ introduced above. By definition, the representation $\rho_1$ is not of the form $\Ind_{F'/F}(\chi)$, for any quadratic extension $F'/F$ and any quasicharacter $\chi$ of 
$\hbox{\scr W}_{F'}$. Also, $\rho_1$ remains irreducible on restriction to $\hbox{\scr W}_K$, for any finite, tamely ramified
extension $K/F$.
 
We need some notation  attached to the representation $\rho_2 = \Ind_{L/F}(\theta)$. Let $K/F$ be the normal closure of $L/F$, and let $E/F$ be the
maximal unramified sub-extension of $K/F$. Then $[E:F]=2$, $K=LE$, and ${\rm Gal}(K/F)$ is isomorphic to the symmetric group $S_3$. We let $g
\in {\rm Gal}(K/F)$ have order 3, and we let $h$ be the nontrivial involution which fixes $L$.
 
We write $\ado$ for the adjoint action of $\GL2{\Bbb C}$ on the Lie algebra of $\roman{{\rm SL}}_2(\Bbb C)$. Thus if $\tau$ is a two-dimensional representation of, for example, a finite group $G$, we have
$$
\tau\otimes\dual\tau = 1_G \oplus \ado(\tau),
$$
where $1_G$ denotes the trivial representation of $G$. We note that the representation $\ado(\tau)$ is reducible if and only if there is a nontrivial one-dimensional representation $\alpha$ of $G$ such that $\alpha\otimes\tau \cong\tau$. In that case, $\alpha$ is a component of $\ado(\tau)$.
 
We can now state and prove the main result of this section.
\nonumproclaim{Theorem 1}
{\rm (1)}
The representation $\rho = \rho_1\otimes \rho_2$ is reducible if and only if there is a quasicharacter $\chi$ of $\hbox{\scr W}_F$ such that
$$
\rho_2 \cong \chi\otimes \ado(\rho_1).
\tag 1.1
$$
\item{\rm (2)}
When condition\/ {\rm (1.1)} is satisfied{\rm ,}  
$$
\rho \cong (\chi\otimes \rho_1)\oplus (\chi\otimes\rho_1\otimes\lambda),
$$
where $\lambda$ denotes the unique irreducible two\/{\rm -}\/dimensional representation of\/ ${\rm Gal} (K/F)${\rm .}
 
\endproclaim

\demo{Proof}
If $\tau$ is a representation of $\hbox{\scr W}_F$ and $F'/F$ is a finite extension, we write $\tau_{F'}$
 for the restriction of $\tau$ to $\hbox{\scr W}_{F'}$. With this notation, 
$$
\rho=\rho_1\otimes \rho_2 = \Ind_{L/F}(\theta\otimes \rho_{1L}).
$$
Let $\langle\,,\,\rangle$ denote the standard inner product of semisimple representations. The representation $\rho$ is thus reducible if and only if
$$
\langle\rho,\rho\rangle = \langle \rho\otimes\dual\rho,1\rangle \ge2\,;
$$
that is, $\rho$ is reducible if and only if $\rho\otimes\dual\rho$ contains the trivial representation with multiplicity at least two.
 
With this in mind, the Mackey restriction-induction formula gives
$$
\rho\otimes \dual\rho = \Ind_{L/F}\left[\theta\otimes \rho_{1L} \otimes \theta^{-1} \otimes \dual\rho_{1L}\right]
\oplus \Ind_{K/F}\left[\theta_K \otimes \rho_{1K} \otimes \theta_K^{-g} \otimes \dual \rho_{1K}\right].
$$
The representation $\rho_{1L}$ is irreducible. It follows that the first term in this expression contains the trivial character with multiplicity one. Thus $\rho\otimes \dual\rho$ is reducible if and only if the second term contains the trivial representation. But, since $\rho_{2E} = \Ind_{K/E}(\theta_K)$ is irreducible, 
$\theta_K^g \neq \theta_K$, and so it follows that:
\nonumproclaim{Lemma 1}
The representation $\rho$ is reducible if and only if $\theta_K^{1-g} \otimes \ado(\rho_{1K})$ 
contains the trivial representation of $\hbox{\scr W}_K${\rm .}
\endproclaim
The representation $\rho_{1K}$ is irreducible, since $\rho_1$ is primitive and $K/F$ is tame. By \cite{9} therefore, there are just two, mutually exclusive possibilities: either $\rho_{1K}$ is primitive, or else it is {\it triply imprimitive.} This means there is a Galois extension $K'/K$, with ${\rm Gal}(K'/K) \cong \Bbb Z/2\Bbb Z \times \Bbb Z/2\Bbb Z$, such that $\rho_{1K}$ is induced from each of the three quadratic sub-extensions of $K'/K$.
\nonumproclaim{Lemma 2}
{\rm (1)}
If $\rho_{1K}$ is primitive{\rm ,} then $\ado(\rho_{1K})$ is irreducible{\rm .}
\item{(\rm 2)}
Suppose that $\rho_{1K}$ is triply imprimitive{\rm ,}
 attached to a quartic extension $K'/K$ as above{\rm .} Let $\eta$ be a nontrivial character of\/ ${\rm Gal}(K'/K)${\rm .} Then
$$
\ado(\rho_{1K}) = \eta\oplus \eta^g \oplus \eta^{g^2}.
$$
The characters $\eta${\rm ,} $\eta^g$ and $\eta^{g^2}$ are distinct{\rm .}
\endproclaim

\demo{Proof}
The first assertion, and the fact that $\eta$ is a component of $\ado(\rho_{1K})$ in part (2), both follow from earlier remarks. Indeed, these remarks show that, in part (2), $\ado(\rho_{1K})$ is the direct sum of the nontrivial characters of ${\rm Gal}(K'/K)$. The element $g$ certainly acts on this set of characters: we have to show that this action is nontrivial. If the action were trivial then, since $K/E$ is tame, each of these characters would be the restriction of a character of $\hbox{\scr W}_E$. It would follow that $\ado(\rho_{1E})$ is a direct sum of abelian characters, as well as being the restriction to $\hbox{\scr W}_E$ of the irreducible $3$-dimensional  representation $\ado(\rho_1)$. Since $\hbox{\scr W}_E$ is a normal subgroup of $\hbox{\scr W}_F$ of index $2$, this is
 impossible.  
\enddemo

Let us deal with the theorem in the case where $\rho_{1K}$ is primitive. By Lemmas 1 and 2, the representation $\rho$ is irreducible. On the other hand, $\rho_{2K}$ is certainly reducible, so (1.1) can never hold. Part (2) is vacuous in this case, so this proves the theorem for $\rho_{1K}$ primitive.
 
From now on, therefore, we assume that {\it $\rho_{1K}$ is triply imprimitive.}
\nonumproclaim{Lemma 3}
Suppose that $\rho_{1K}$ is triply imprimitive and that $\rho$ is reducible{\rm .}
 There exists a nontrivial character $\eta$ of\/ ${\rm Gal}(K'/K)$ such that
$$
\align
\eta\otimes \rho_{1K} &\cong \rho_{1K}, \\
\eta^h &= \eta, \\
(\theta_K/\eta)^g &= \theta_K/\eta.
\endalign
$$
\endproclaim

\demo{Proof}
We have $\rho_{1K} \cong \chi\otimes \rho_{1K}$ for every character
$\chi$ of the noncyclic $4$-group ${\rm Gal} (K'/K)$, so the first property holds for any choice of $\eta$. Since $\ado(\rho_{1K}) = \eta \oplus \eta^g \oplus \eta^{g^2}$, Lemma 1 implies that $\theta_K^{g-1}$ is one of $\eta$, $\eta^g$, $\eta^{g^2}$. In particular, $\theta_K^{g-1}$ has order 2. By definition, $\theta_K$ is invariant under $h$, so
that
$$
(\theta_K^{(g-1)g})^h = \theta_K^{g-g^2} = (\theta_K^{g-g^2})^{-1} = \theta_K^{(g-1)g} .
$$
Thus, $\theta_K^{(g-1)g}$ is also invariant under $h$. Changing notation if necessary, we can assume $\theta_K^{(g-1)g} = \eta$, which is therefore
fixed by $h$. Since $\eta.\eta^g. \eta^{g^2} =1$, we also have $\eta = \eta^{(g-1)g}$.
 It follows that $\theta_K/\eta$ is fixed by $g$ as required.  
\enddemo

We assume that $\rho_{1K}$ is reducible. Since $\theta_K/\eta$ (as in Lemma 3) is fixed by $g$, there exists a quasicharacter $\chi$ of $\hbox{\scr W}_E$ such that $\theta_K = \eta\cdot\chi_K$. Since $\theta_K$ and $\eta$ are fixed by $h$, so is $\chi_K$. Thus $\chi^{h-1}$ is a character of ${\rm Gal}(K/E)$. Moreover, $\chi$ is only determined modulo the character group of ${\rm Gal}(K/E)$, on which the element $h$ acts nontrivially. We can therefore choose $\chi$ to satisfy $\chi^h=\chi$. Thus $\chi = \chi'_E$, for some character $\chi'$ of $\hbox{\scr W}_F$. Likewise, $\eta = \eta'_K$, for some character $\eta'$ of $\hbox{\scr W}_L$. 
We have $\theta_K = \eta'_K\chi'_K$, so
$$
\theta = \eta'\chi'_L \quad \text{or}\quad \theta = \omega_{K/L}\, \eta'\chi'_L,
$$
where $\omega_{K/L}$ is the nontrivial character of ${\rm Gal}(K/L)$. This, however, is the restriction to $\hbox{\scr W}_L$ of the nontrivial character $\omega_{E/F}$ of ${\rm Gal}(E/F)$, so we can absorb this into $\chi'$ and assume $\theta = \eta'\chi'_L$.
 
Thus, when $\rho$ is reducible, 
$$
\rho_2 = \chi' \otimes \Ind_{L/F}(\eta'),
$$
for some quasicharacter $\chi'$ of $\hbox{\scr W}_F$. On the other hand, $\ado(\rho_{1L})$ contains either $\eta'$ or $\omega_{K/L}\eta'$. It follows that either $\ado(\rho_1)$ or $\omega_{E/F} \otimes \ado(\rho_1)$ is equivalent to $\Ind_{L/F}(\eta')$, so we conclude that, when $\rho$ is reducible, $\rho_2$ is a twist of $\ado(\rho_1)$.
 
We now treat the converse statement in part (1) of the theorem. We assume that $\rho_2 = \chi \otimes \ado(\rho_1)$ (and that $\rho_{1K}$ is triply imprimitive). Since, by definition, $\rho_1\otimes \dual \rho_1 = 1 \oplus
\ado(\rho_1)$, 
$$
\langle \rho_1,\rho_1\otimes \ado(\rho_1)\rangle = \langle \rho_1\otimes \dual\rho_1,\ado(\rho_1)\rangle \ge 1,
$$
and so $\rho_1\otimes \ado(\rho_1)$ contains $\rho_1$. Write $\ado(\rho_1) = \Ind_{L/F}(\eta')$, for some quasicharacter $\eta'$ of $\hbox{\scr W}_L$ such that $\eta'_K = \eta$. Then
$$
\rho_1 \otimes \ado(\rho_1) = \Ind_{L/F}(\eta'\otimes \rho_{1L}).
$$
We have $(\eta'\otimes \rho_{1L})_K = \eta\otimes \rho_{1K} = \rho_{1K}$, so that $\eta'\otimes \rho_{1L}$ is either $\rho_{1L}$ or
$\omega_{K/L}\rho_{1L}$. In the first case, we get $\rho_1\otimes \ado(\rho_1) = \rho_1 \oplus (\rho_1\otimes\lambda)$, and in the second,
$\omega_{E/F}\rho_1 \oplus (\rho_1\otimes\lambda)$. (In particular, $\rho$ is reducible, which completes the proof of (1).) Since $\rho_1$ is a
component of $\rho_1\otimes \ado(\rho_1)$, it is the first case which must hold.  This gives the formula required for part (2), and we have
completed
 the proof of Theorem 1.  \hfill\qed
\enddemo

\nonumproclaim{Variant}
Assume the hypotheses and notation of Theorem {\rm 1,} except that $L/F$  is now cyclic and totally ramified{\rm .} Then\/{\rm :}
\vglue2pt
\item{\rm (1)}
The representation $\rho$ is reducible if and only if there is a quasicharacter $\chi$ of $\hbox{\scr W}_F$ such that
$$
\rho_2 \cong \chi\otimes \ado(\rho_1).
\tag 1.2
$$
\item{\rm (2)}
When {\rm (1.2)} is satisfied{\rm ,} we have
$$
\rho \cong \chi\otimes \rho_1 \otimes \Ind_{L/F}\,1_L,
$$
where $1_L$ denotes the trivial character of $\hbox{\scr W}_L${\rm .}
 
\endproclaim

\demo{Proof}
This is parallel to that of the theorem, setting $E=F$, $K=L$, and $h=1$. 
\enddemo

We revert to the hypotheses of Theorem 1. It will be useful to have some further control of the representation $\rho=\rho_1\otimes \rho_2$. Recall that an irreducible representation $\sigma$ of $\hbox{\scr W}_F$ is
 {\it totally ramified\/} if $\sigma\otimes \chi \not \cong\sigma$ for any unramified quasicharacter $\chi\neq 1$ of $\hbox{\scr W}_F$.
Equivalently, $\sigma_{F'}$ is irreducible (and totally ramified) for any finite unramified extension $F'/F$. For example, the representations
$\rho_1$, $\rho_2$ are totally ramified.
\nonumproclaim{Proposition  1}
In the situation of Theorem {\rm 1,} suppose that $\rho$ is irreducible{\rm .} Then $\rho$ is totally ramified{\rm .}
\endproclaim

\demo{Proof}
Suppose otherwise; then there is an unramified quasicharacter $\alpha\neq 1$ of $\hbox{\scr W}_F$ such that $\alpha \otimes \rho \cong \rho$. The
quasicharacter $\alpha$ has finite order dividing $6$, so we assume it has order $d=2$ or $3$. We may view $\alpha$ as a character of ${\rm
Gal}(M/F)$, where $M/F$ is unramified of degree $d$.

\nonumproclaim{Lemma 4}
The representation $\rho_{1M}$ is primitive{\rm .}
\endproclaim

\demo{Proof}
Suppose that $\rho_{1M}$ is triply imprimitive, and take first the case $d=3$. Thus there is  a quadratic extension $M'/M$ and a quasicharacter $\xi$
of
$\hbox{\scr W}_{M'}$ such that $\rho_{1M} = \Ind_{M'/M}(\xi)$. Next, $\rho_{2M'}$ is irreducible, so $\rho_{M'}$ is the direct sum of two irreducible
3-dimensional representations. On the other hand, the relation $\chi\otimes\rho \cong \rho$ implies that $\rho$ is induced from $\hbox{\scr W}_M$,
and that $\rho_M$ is the direct sum of three irreducible 2-dimensional representations. This case can therefore not arise.
 
Take next the case $d=2$. The representation $\ado(\rho_1)$ is irreducible of dimension 3, so that  $\ado(\rho_1)_M = \ado(\rho_{1M})$ is
irreducible. On the other hand, if $\rho_{1M}$ is triply imprimitive, the representation $\ado(\rho_{1M})$ is the direct sum of three quasicharacters.
Again, this case cannot arise.  
\enddemo

Now we prove the proposition. Suppose first that $M/F$ is of degree 3. The hypotheses of the theorem apply to the representations $\rho_{iM}$, and
so we get a quasicharacter $\chi$ of $\hbox{\scr W}_M$ such that $\rho_{2M} \cong \chi\otimes \ado(\rho_{1M})$. If $y$ generates ${\rm Gal}(M/F)$,
we have $\chi^{y-1}\otimes \rho_{2M} \cong \rho_{2M}$. The representation $\rho_{2M}$ is totally ramified and induced from the noncyclic cubic
extension $LM/M$, so that  $\chi^y=\chi$. Thus there is a quasicharacter $\chi'$ of $\hbox{\scr W}_F$ with $\chi'_M=\chi$ and $\rho_2 \cong
\chi'\otimes
\ado(\rho_1)$. This is impossible, by Theorem 1.
 
This reduces to the case where $M/F$ is quadratic. We can apply
 the variant of the theorem to $\rho_{1M}$ and $\rho_{2M}$ to get a character $\chi$ of $\hbox{\scr W}_M$ such that $\rho_{2M} \cong \chi \otimes
\ado(\rho_{1M})$. If $h$ is a generator of ${\rm Gal}(M/F)$, the character $\chi^{h-1}$ fixes $\rho_{2M}$, and so is a character of ${\rm Gal}(LM/M)$.
Since $\chi$ is only determined modulo the character group of ${\rm Gal}(LM/M)$, we can in fact assume that $\chi$ is fixed by $h$ and finish the
proof as in the last case.\hfill\qed  
\enddemo
\section{Conductors of pairs}

In this section only, $F$ denotes an arbitrary non-Archimedean local field, with finite residue field of $q$ elements. Let $N$ denote a positive integer; we write $G_N= \GL NF$ and also $A_N = \roman M_N(F)$, the algebra of $N\times N$-matrices over $F$.
 
We shall use the classification theory of \cite{3}. In particular, we need the notion of a {\it simple stratum\/} in $A_N$, as in  [3, 1.5].
Attached to the simple stratum $[\frak A,n,0,\beta]$, we have the compact open subgroups $H^r(\beta,\frak A)$, $r\ge1$, of $G_N$. Next, we fix a
continuous character $\psi_0$ of $F$ nontrivial on the discrete valuation ring $\frak o$ in $F$, but trivial on the maximal ideal $\frak p$ of $\frak
o$. Such a choice of $\psi_0$ gives rise to a set $\hbox{\scr C}(\frak A,0,\beta)=\hbox{\scr C}(\frak A,\beta)$ of characters of $H^1(\beta,\frak A)$, called {\it
simple characters.} (The terminology and notation are as in [3, Ch.\,3].)
 
Let $\pi$ be an irreducible supercuspidal representation of $G_N$. If $\chi$ is a quasicharacter of $F^\times$, we write $\chi\cdot\pi$ for the representation $g\mapsto \chi(\det g)\pi(g)$ of $G_N$. We say that $\pi$ is 
{\it totally ramified\/} if $\chi\cdot\pi \not\cong \pi$ for any unramified quasicharacter $\chi\neq 1$ of $F^\times$.

\nonumproclaim{Proposition  2}
Let $\pi$ be an irreducible supercuspidal representation of $G_N${\rm ,} $N\ge2$. The following are equivalent\/{\rm :}
\vglue2pt
\item{\rm (1)}
$\pi$ is totally ramified\/{\rm ;}\/
\item{\rm (2)}
There is a simple stratum $[\frak A,n,0,\beta]$ in $A_N$ and a simple character $\theta\in \hbox{\scr C}(\frak A,\beta)$ such that\/{\rm :}
\itemitem{\rm (a)} $\theta$ occurs in $\pi${\rm ,} and
\itemitem{\rm (b)} the field extension $F[\beta]/F$ is totally ramified of degree $N${\rm .}
 
\endproclaim

\demo{Proof}
This follows immediately from  [3, (6.2.5), (8.4.1)].  
\enddemo

The integer $n=n(\pi)$ is an invariant of the totally ramified representation~$\pi$.

\nonumproclaim{Lemma 5}
Let $\pi$ be a totally ramified{\rm ,} irreducible{\rm ,} supercuspidal representation of $G_N${\rm ,} $N\ge2${\rm .} The following are equivalent\/{\rm
:}
\vglue2pt
\item{\rm (1)}
There is a quasicharacter $\chi$ of $F^\times$ such that $n(\chi\cdot\pi)<n(\pi)${\rm ;}
\item{\rm (2)}
$n(\pi) \equiv 0 \pmod N${\rm .}
 
\endproclaim

\demo{Proof}
Take a simple stratum $[\frak A,n,0,\beta]$ and a simple character $\theta$ occurring in $\pi$, as in Proposition 2. The integer $n=n(\pi)$ is then
$-\ups_{F[\beta]}(\beta)$, where $\ups$ denotes the normalized additive valuation.
 
We choose a simple stratum $[\frak A,n,n{-}1,\alpha]$ in $A_N$ which is equivalent to $[\frak A,n,n{-}1,\beta]$. The algebra $F[\alpha]$ is a field. The ramification index $e(F[\alpha] \mid F)$ and residue class degree $f(F[\alpha]\mid F)$ divide the corresponding invariants of $F[\beta]/F$. In particular, $F[\alpha]/F$ is totally ramified.
 
Since $[\frak A,n,n{-}1,\alpha]$ is simple, the integer
$$
\ups_{F[\alpha]}(\alpha) = -e(F[\beta]\mid F)^{-1}e(F[\alpha]\mid F)\,n = -N^{-1}e(F[\alpha]\mid F)\,n
$$
is relatively prime to $e(F[\alpha]\mid F) = [F[\alpha]:F]$. If $N$ divides $n$ therefore, $e(F[\alpha]\mid F)$ must divide $\ups_{F[\alpha]}(\alpha)$, and this forces $F[\alpha] = F$.
 
By definition  [3, 3.2], the restriction of the character $\theta$ to the congruence unit group $U^n(\frak A)$ of the order $\frak A$ takes the
form
$$
x\longmapsto \psi_0(\roman{tr}(\alpha(x{-}1)),\quad x\in U^n(\frak A),
$$
where tr denotes the matrix trace $A_N\to F$. Since $\alpha\in F$, there is a quasicharacter $\chi$ of $F^\times$ such that $\chi(x) = \psi_0(\alpha(x{-}1))$, for $x\in U^{n/N}_F$. The representation $\pi_1=\chi^{-1}\cdot \pi$ then satisfies $n(\pi_1) <n$.
 
Thus $(2)\,\Rightarrow (1)$; the converse follows similarly.  
\enddemo

Let $N$, $N'$ be positive integers, and let $\pi$, $\pi'$ be irreducible smooth representations of $G_N$, $G_{N'}$ respectively. We form the local constant $\ve(\pi\times \pi',s,\psi)$ of the pair $(\pi,\pi')$, as in 
\cite{6} or \cite{8}. Here, $s$ is a complex variable and $\psi$ is a nontrivial character of the additive group of $F$. The local constant takes the
form
$$
\ve(\pi\times\pi',s,\psi) = q^{(\frac12 -s)f(\pi\times\pi',\psi)}\,\ve(\pi\times\pi',\tfrac 12,\psi),
$$
for a certain integer $f(\pi\times\pi',\psi)$. Let $c(\psi)$ denote the conductor of $\psi$; thus $c(\psi)$ is the least integer $k$ such that $\psi$ is trivial on $\frak p^k$. There is then an integer $a(\pi\times\pi')$, independent of $\psi$, such that
$$
f(\pi\times\pi',\psi) = a(\pi\times\pi') - NN'c(\psi).
$$
We can now give the main result of the section.

\nonumproclaim{Theorem 2}
Let $\pi$ be an irreducible{\rm ,} totally ramified supercuspidal representation of $G_6${\rm .}
 There are a positive integer $m$ dividing $3${\rm ,} and an irreducible{\rm ,} totally ramified supercuspidal representation
 $\pi'$ of $G_m$ such that
$$
a(\pi\times\pi') \not\equiv 0 \pmod 2.
$$
\endproclaim

\demo{Proof}
We choose a simple stratum $[\frak A,n,0,\beta]$ and a simple character $\theta \in \hbox{\scr C}(\frak A,\beta)$ which occurs in $\pi$. If $\chi$ is a quasicharacter of
$F^\times$, 
$$
a(\chi\cdot\pi \times \pi') = a(\pi\times\chi\cdot\pi').
$$
Therefore, by Lemma 5, we can replace $\pi$ by $\chi\cdot\pi$, for a suitable $\chi$, to ensure that $n \not\equiv 0 \pmod 6$. That done, if we choose a simple stratum $[\frak A,n,n{-}1,\gamma]$ equivalent to $[\frak A,n,n{-}1,\beta]$, we have $[F[\gamma]:F] >1$. The quantity
$$
-\ups_{F[\gamma]}(\gamma) = e(F[\gamma]\mid F)\,n/6
$$
is relatively prime to $e(F[\gamma]\mid F) = [F[\gamma]:F]$ so, if $[F[\gamma]:F]$ is even, the integer $n$ is odd and there is nothing to prove. We therefore assume $[F[\gamma]:F] = 3$.
 
Let $r$ denote the integer $-k_0(\beta,\frak A)$ (notation of [3, 1.4.5]). By the assumption on $F[\gamma]$, we have $n>r\ge 1$. 
We choose a simple stratum $[\frak A,n,r,\gamma']$ equivalent to $[\frak A,n,r,\beta]$. The field extension $F[\gamma']/F$ is again of degree 3
(since it divides $6$ properly, and is divisible by $[F[\gamma]:F]=3$). We could have taken $\gamma= \gamma'$ above, and so we now simplify
notation by
  doing so.

\nonumproclaim{Lemma 6}
The integer $r$ is odd{\rm .}
\endproclaim

\demo{Proof}
Let $B$ denote the centralizer of $\gamma$ in $A=A_6$, and put $\frak B = \frak A\cap B$. Let $s_\gamma:A \to B$ be a tame corestriction
 ([3, 1.3]) relative to $F[\gamma]/F$, and write $\beta = \gamma +c$. Then $[\frak B,r,r{-}1,s_\gamma(c)]$ is equivalent (in $B$) to a simple
stratum $[\frak B, r,r{-}1,\alpha]$ (\/{\it ibid\.} (2.4.1)). Now,
$$
\align
e(F[\gamma,\alpha]\mid F) &= e(F[\beta]\mid F) = 6, \\
f(F[\gamma,\alpha]\mid F) &= f(F[\beta]\mid F) = 1,
\endalign
$$
and so the field extension $F[\alpha,\gamma]/F[\gamma]$ is totally ramified of degree $2$. However, by the definition of simple stratum, the integer
$r=-\ups_{F[\alpha,\gamma]}(\alpha)$ is relatively prime to $e(F[\alpha,\gamma]\mid F[\gamma])$, and so $r$ is odd, as required.  
\enddemo

We now start the construction of the desired representation $\pi'$. 
First, $\theta\in \hbox{\scr C}(\frak A,\beta)$ is a character of the group $H^1(\beta,\frak A)$. We can also form $H^1(\gamma,\frak A)$, and we
have $H^{r+1}(\gamma,\frak A)= H^{r+1}(\beta,\frak A)$ [3, (3.1.9)]. There is a simple character $\theta_0\in \hbox{\scr C}(\frak A,\gamma)$
such that
$$
\theta_0\mid H^{r+1}(\beta,\frak A) = \theta\mid H^{r+1}(\beta,\frak A)
$$
(\/{\it ibid\.} (3.3.21)). We now embed the field $F[\gamma]$ in $A_3 = \roman M_3(F)$; we then get a simple stratum $[\frak A',n',0,\gamma]$ in
$A_3$, for a uniquely determined principal order $\frak A'$, with $n'=n/2$. We have a canonical bijection between $\hbox{\scr C}(\frak A,\gamma)$
and $\hbox{\scr C}(\frak A',\gamma)$ [3, (3.6.14)]: let $\theta_1 \in \hbox{\scr C}(\frak A',\gamma)$ correspond to $\theta_0$.
 
Let $\pi_1$ be an irreducible representation of $G_3$ which contains $\theta_1$. Then $\pi_1$ is supercuspidal and totally ramified.
 In the language of [2, \S6], the pair $(\pi,\pi_1)$ has a best common approximation of the form $([\varLambda,n,0,\gamma],r, \vartheta)$, with
$r$ and $n$ as above. We apply the formula of [2, 6.5(iii)]  to get
$$
a(\pi\times\dual\pi_1) = 18\left(1+\frac{\frak c(\gamma)}9 +\frac r{18}\right) ,
$$
for a certain integer $\frak c(\gamma)$, and this yields
$$
a(\pi\times\dual\pi_1) \equiv r \equiv 1 \pmod 2.
$$
The representation $\pi'=\dual\pi_1$ thus has the required property.\hfill\qed  
\enddemo

\section{Proof of the Main Theorem}

We now prove the Main Theorem stated in the introduction. We therefore revert to the notation  there and  in \S1. We treat first the case where
$\rho$ is irreducible. Thus $\pi$ is supercuspidal. The representation $\rho_E$ is irreducible by Proposition 1, and $\pi_E$ corresponds to $\rho_E$
under the Langlands correspondence. The relation $\pi_E \cong \pi'_E$ implies that $\pi'_E$ is supercuspidal, and hence that $\pi'$ is supercuspidal.
 
Since $\pi$ is supercuspidal, we have $\gamma(\pi\times\sigma,s,\psi) = \ve(\pi\times\sigma,s,\psi)$, for any irreducible smooth representation $\sigma$ of $\GL mF$, $m\le 5$.
That is similar  for $\pi'$, and the relation $\ve(\pi\times\sigma,s,\psi) = \ve(\pi' \times \sigma,s,\psi)$ implies $a(\pi\times \sigma) =
a(\pi'\times
\sigma)$. The relation $\pi_E \cong \pi'_E$ further implies that either $\pi \cong \pi'$ or $\pi \cong \omega_{E/F}\cdot \pi'$, where
$\omega_{E/F}$ is the character of $F^\times$ corresponding to the nontrivial character of ${\rm Gal}(E/F)$ [1, I Prop\. 6.7]. By Theorem 2, we
can find an irreducible supercuspidal representation $\sigma$ of $\GL mF$, $m=1$ or $3$, such that $a(\pi\times \sigma)= a(\pi'\times\sigma)$ is
odd. Then
$$
\ve(\omega_{E/F}\pi'\times \sigma,s,\psi) = (-1)^{a(\pi'\times\sigma)} \ve(\pi'\times \sigma,s, \psi) = - \ve(\pi\times\sigma,s,\psi).
$$
It is therefore only the case $\pi\cong\pi'$ which arises, as desired.
 
We now assume that $\rho$ is reducible. Thus $\rho_2 = \chi\otimes\ado(\rho_1)$, for some quasicharacter $\chi$ of $\hbox{\scr W}_F$ (Theorem 1). Clearly, it is enough to treat the case $\chi=1$, so that $\rho = \rho_1\oplus (\rho_1\otimes\lambda)$, 
as in Theorem 1(2). The representation $\lambda$ is of the form $\Ind_{E/F}(\phi)$, where $\phi$ is a nontrivial character of ${\rm Gal}(K/E)$. It
follows easily that $\rho_1\otimes\lambda$ is irreducible, and satisfies
$$
\rho_1\otimes\lambda\otimes\omega_{E/F} \cong \rho_1\otimes\lambda.
$$
On the other hand, $\rho_1$ is irreducible and totally ramified.
 
Let $\pi_1$ be the irreducible supercuspidal representation of $\GL2F$ corresponding to $\rho_1$, and $\varPi$ that of $\GL4F$ corresponding to $\rho_1\otimes\lambda$. Thus $\pi_1$ is totally ramified, while $\omega_{E/F}\cdot\varPi \cong \varPi$ (where $\omega_{E/F}$ now denotes the character of $F^\times$ corresponding to the nontrivial character of ${\rm Gal}(E/F)$). We have
$$
\gamma(\pi\times\dual\pi_1,s,\psi) = \frac{1-q^{-s}}{1-q^{s-1}}\,\,\ve(\pi_1\times \dual\pi_1,s,\psi)\,\ve(\varPi\times\dual\pi_1,s,\psi).
\tag 3.1
$$
The representation $\varPi_E$ corresponds, via the Langlands correspondence, to $(\rho_1\otimes \lambda)_E = (\phi\otimes \rho_{1E}) \oplus (\phi^2\otimes\rho_{1E})$. Thus
$$
\pi_E \cong \pi'_E \cong \pi_{1E}\,\boxplus\, \phi\cdot \pi_{1E} \,\boxplus\, \phi^2\cdot \pi_{1E},
$$
where $\phi$ now denotes the character of $E^\times$ corresponding to some nontrivial character of ${\rm Gal}(K/E)$. This implies
$$
\pi' \cong \omega_{E/F}^a\pi_1 \,\boxplus\, \varPi,
$$
for some integer $a$. We now have
$$
\gamma(\pi'\times\dual\pi_1,s,\psi) = \frac{1-(-1)^aq^{-s}}{1-(-1)^aq^{s-1}}\,\, \ve(\omega_{E/F}^a \pi_1\times\dual\pi_1,s,\psi)\,\ve(\varPi\times\dual\pi_1,s,\psi).
$$
Comparing the location of the poles of this expression with those
of (3.1), we get $(-1)^a=1$ and $\omega_{E/F}^a = 1$, as desired. \hfill\qed  

\vglue-12pt

\AuthorRefNames[10]
\references

[1]
\name{J.\ Arthur} and \name{L.\ Clozel},
{\it Simple algebras{\rm ,} base change{\rm ,} and the advanced theory of the trace 
formula},
{\it Ann.\ of Math. Studies}  {\bf 120}, Princeton Univ.\  Press,
Princeton, NJ, 1989.
 
[2]
\name{C.\ J.\ Bushnell, G.\ Henniart,} and \name{P.\ C.\ Kutzko},
 Local Rankin-Selberg convolutions for $\roman{GL}_n$: Explicit conductor formula,
{\it J.A.M.S.}  {\bf 11} (1998), 703--730.
 
[3]
\name{C.\ J.\ Bushnell} and \name{P.\ C.\ Kutzko},
{ The admissible dual of ${\rm GL}(N)$ via compact open subgroups,
{\it Ann.\ of Math.\ Studies} {\bf 129}, Princeton 
Univ.\  Press, Princeton, NJ, 1993.
 
[4]
\name{M.\ Harris} and \name{R.\ Taylor},
 The geometry and cohomology of some simple Shimura varieties,
{\it Ann.\ of Math. Studies} {\bf 151},   Princeton 
Univ.\  Press, Princeton, NJ, 2001.
 
[5]
 \name{G.\ Henniart},
 Une preuve simple des conjectures locales de Langlands pour $\roman{GL}_n$ sur un corps $p$-adique,
{\it Invent.\ Math.\/}  {\bf 139} (2000), 439--455.

[6]
 \name{H.\ Jacquet, I.\ I.\ Piatetskii-Shapiro}, and \name{J.\ A.\ Shalika},
 Rankin-Selberg convolutions,
{\it Amer. J. Math.\/}  {\bf 105} (1983), 367--464.

[7] 
\name{H. H.\ Kim} and \name{F.\ Shahidi},
 Functorial products for $\roman{GL}_2 \times \roman{GL}_3$ and the 
symmetric cube for $\roman{GL}_2$,
{\it Ann.\ of  Math.\/}  {\bf 155} (2002), 837--883.

[8] 
\name{F.\ Shahidi},
 Fourier transforms of intertwining operators and Plancherel measures for $\roman{GL}(n)$, {\it Amer.\ J. Math.\/} {\bf 106}  (1984), 67--111.
 
[9]
\name{A.\ Weil},
 Exercices dyadiques,
{\it Invent.\ Math.\/}  {\bf 27} (1974), 1--22.
 
\endreferences

\centerline{(Received October 2, 2000)}
\enddocument

%% file: Shahidi.ref
 \AuthorRefNames [CKP-SS]

\references
[A1] \name{J.\ Arthur}, Endoscopic $L$-functions and a combinatorial
identity, {\it Canadian J. Math\/}.\ {\bf 51} (1999), 1135--1148.

[A2] \bibline, Unipotent automorphic representations:
conjectures,
in {\it Repr\'esentations II. Groupes $p$-adiques et R\'eels\/},
{\it  Ast\'erisque} {\bf 171--172} (1989),. 13--71.

[AC] \name{J. Arthur} and \name{L. Clozel},
{\it  Simple Algebras, Base Change, and the Advanced Theory of the 
Trace Formula\/}, {\it  Ann.\  of Math.\ Studies}  {\bf 120}, Princeton Univ.\
Press, Princeton, NJ, 1989.

[As] \name{M.\ Asgari},
 Local $L$-functions for split spinor groups, {\it Canadian J.\ Math\/}., to appear.

 [Bor] \name{A.~Borel},  Automorphic $L$-functions, in  {\it Automorphic Forms and 
Automorphic Representations\/}, {\it Proc.~Sympos.~Pure Math\/}.\ {\bf 33} (1979), 27--61,
Amer.~Math.~Soc.,  Providence, RI.

[Bou] \name{N.\ Bourbaki},
{\it  Groupes et Alg\`ebres de Lie\/}, {\rm Chap.\ 4--6},
Paris, Hermann, 1968.

[BK] \name{A.~Braverman} and \name{D.~Kazhdan}, $\gamma$-functions of representations and lifting, in
{\it Visions in Mathematics Towards\/} 2000, Birkh\"auser, 2000.

[BDHI] \name{D. Bump, W. Duke, J.\ Hoffstein}, and \name{H.\ Iwaniec},
An estimate for the Hecke eigenvalues of Maass forms,
{\it IMRN} {\bf 4} (1992), 75--81.

[BH] \name{C.~J.~Bushnell} and \name{G.~Henniart},  On certain dyadic
representations (Appendix to this paper), {\it Ann. of Math\/}.\ {\bf 155} (2002), 883--893.

[BHK] \name{C.\ J.\ Bushnell, G.\ Henniart}, and \name{P.\ C.\ Kutzko},
 Local Rankin-Selberg convolutions for $\text{GL}_n$: Explicit conductor
formula, {\it J.\ Amer.\ Math.\ Soc\/}.\ {\bf 11} (1998), 703--730.

[CS] \name{W.\ Casselman} and \name{F.\ Shahidi},
 On irreducibility of standard modules for generic representations,
{\it  Ann.\ Sci.\  {\rm \'{\it E}}cole Norm.\ Sup\/}.\ {\bf 31} (1998),  561--589.

[Ch] \name{J.-P.~Jeff Chen},  Local factors, Central characters, and Representations of general linear 
groups over non-archimedean fields,  Doctoral dissertation, Yale University (1996).

[CP-S1] \name{J.\ W.\ Cogdell} and \name{I.\ I.\ Piatetski-Shapiro},
Converse theorems for $\text{GL}_n$ II, {\it  J.\ Reine Angew.\ Math\/}.\ {\bf 507} (1999),
165--188.

[CP-S2]  \bibline,
Unitarity and functoriality,
{\it Geom.\ Funct.\ Anal\/}.\ {\bf 5} (1995), 164--173.

[CP-S3] \bibline,
Converse theorems for $\text{GL}_n$ and their applications to lifting,
{\it Proc.\ Internat.\ Conf.\ on Cohomology of Arithmetic Groups\/}, {\it $L$-functions, and 
Automorphic Forms\/}, Tata Institute of Fundamental Research  (Mumbai 1998), Editor: T. N. Venkataramana, TIFR-Narosa Publ.,
AMS dist. (2001), 1--34.

[CKP-SS] \name{J.\ W.~Cogdell, H.~Kim, I.\ I.~Piatetski-Shapiro}, and \name{F.~Shahidi}, 
On lifting from classical groups to $\text{GL}_N$, {\it Publ.~Math.~IHES\/} 
{\bf 93} (2001), 5--30.

[DI] \name{J.\ M.\ Deshouillers} and \name{H.\ Iwaniec},
Kloosterman sums and Fourier coefficients of cusp forms,
{\it  Invent.\ Math\/}.\ {\bf 70} (1982), 219--288.

[Ge] \name{S. Gelbart},
Three lectures on the modularity of $\overline\rho\sb {E,3}$ and the Langlands reciprocity
conjecture,
in {\it  Modular Forms and Fermat\/{\rm '}\/s Last Theorem\/},
 155--207,
Springer-Verlag, New York, 1997.

[GeJ] \name{S.~Gelbart} and \name{H.~Jacquet},
A relation between automorphic representations of $\text{GL}(2)$ and
$\text{GL}(3)$,
{\it  Ann.~Sci.~{\rm \'{\it E}}cole Norm.~Sup\/}.\ {\bf 11} (1978), 471--552.

[GeS] \name{S.~Gelbart}  and \name{F.~Shahidi},
 Boundedness of automorphic ${L}$-functions in vertical strips,
{\it J.A.M.S}.\ {\bf 14} (2001),  79--107.

[Go] \name{A.\ Good},
 On various means involving the Fourier coefficients of cusp
forms,
{\it Math. Z\/}.\  {\bf 183} (1983),  95--129.

[HT]  \name{M.\ Harris} and \name{R.\ Taylor},
 On the geometry and cohomology of some simple Shimura varieties
{\it Ann.\  of Math.~Studies\/} {\bf 151},  Princeton Univ.\  Press,
 Princeton, NJ, 2001.

[He1] \name{G.\ Henniart},
 Une preuve simple des conjectures de Langlands pour $\text{GL}(n)$ sur un corps $p$-adique,
{\it Invent.~Math\/}.\ {\bf 139} (2000),  439--455.

[He2]  \bibline,  La conjecture de Langlands locale pour
$\text{GL}(3)$, {\it  M{\rm \'{\it e}}moires de la Soc.\  Math.\  de France\/}
{\bf 11/12} (1984),  1--186.

[I] \name{H.\ Iwaniec},
{\it  Introduction to the Spectral Theory of Automorphic Forms},
Rev.\  Mat.\  Iberoamericana, Madrid,  1995.

[IS]
 \name{H.\ Iwaniec} and \name{P.\ Sarnak},
 Perspectives on the analytic theory of $L$-functions,
 {\it Geom.\ Funct.\ Anal\/}.\ {\bf 2000} (Special Issue), 705--741.

[JP-SS1] \name{H.~Jacquet, I.~Piatetski-Shapiro}, and \name{J.~Shalika},
  Rankin-Selberg convolutions,
{\it  Amer.\ J.\ Math\/}.\  {\bf 105} (1983),  367--464.

[JP-SS2] \bibline,
Rel\`evement cubique non normal
{\it C.\ R.\ Acad.\ Sci.\ Paris S\'er.\ I Math\/}.\ {\bf 292} (1981), 
567--571.

[JP-SS3]  \bibline,
 Automorphic forms on $\text{GL}(3)$ I, II,
 {\it Ann.\ of Math\/}.\ {\bf 109} (1979),  169--212, 213--258.

[JS1] \name{H.\ Jacquet} and \name{J.\ Shalika},
 On Euler products and the classification of automorphic forms II, {\it  Amer.\ J.\ Math\/}.\
 {\bf 103} (1981),  777--815.

[JS2] \bibline,
 A lemma on highly ramified $\epsilon$-factors,
{\it  Math.\ Ann\/}.\ {\bf 271} (1985),  319--332.

[Ki1]  \name{H.\ Kim},
 Langlands-Shahidi method and poles of automorphic $L$-functions:\ 
 Application to exterior square $L$-functions,
{\it   Canadian J.\ Math\/}.\ {\bf 51} (1999),  835--849.

[Ki2]  \bibline,
 Langlands-Shahidi method and poles of automorphic $L$-functions II,
{\it  Israel J.\ Math\/}.\ {\bf 117} (2000), 261--284.

[Ki3]  \bibline,
 Correction to: Langlands-Shahidi method and poles of automorphic $L$-functions II,
{\it Israel J.\ Math\/}.\  {\bf 118} (2000), 379.

[Ki4] \bibline,
 Functoriality for the exterior square of $\text{GL}_4$ and symmetric fourth of $\text{GL}_2$,
 preprint, 2000.

[Ki5] \name{H.\ Kim},
 On local $L$-functions and normalized intertwining operators,
 preprint, 2000.

[KS1]  \name{H.\ Kim} and \name{F.\ Shahidi},
Symmetric cube $L$-functions for $\text{GL}_2$ are entire,
{\it Ann. of Math\/}.\ {\bf 150} (1999), 645--662.

[KS2] \bibline,
 Functorial products for $\text{GL}_2\times \text{GL}_3$ and functorial 
symmetric cube for $\text{GL}_2$, {\it C.\ R.\ Acad.\ Sci.\ Paris\/} {\bf 331} (2000),
 599--604.

[KS3] \bibline,
Cuspidality of symmetric powers with applications, 
{\it Duke Math.\ J\/}., to appear.

[Ku] \name{P.~Kutzko},  The Langlands conjecture for $\text{GL}_2$ of a
local field, {\it Ann.\ of Math\/}.\ {\bf 112} (1980),  381--412.

[LL] \name{J.-P.\ Labesse} and \name{R.\ P.\ Langlands},
$L$-indistinguishability for $\text{SL}(2)$, {\it  Canadian J. Math\/}.\ {\bf 31} (1979),
726--785.

[La1]  \name{R.~Langlands},
{\it  On the Functional Equations Satisfied by Eisenstein Series\/}, 
{\it Lecture Notes in Math\/}.\ {\bf 544},  Springer-Verlag, New York, 1976.

[La2] \bibline, Base change for $\text{GL}(2)$,
{\it Ann.\ of Math.~Studies\/} {\bf 96}, Princeton Univ.\ Press,
Princeton, NJ,   
 1980.

[La3] \bibline,  Automorphic representations, Shimura varieties, and
motives, {\it  Ein\break m\"archen\/}, {\it Proc.\ Sympos.\ Pure Math\/}.\ {\bf 33} (1979), 205--246,
A.M.S., 
Providence, RI.

[La4] \bibline,  On the classification of irreducible
representations of real algebraic groups,  in {\it Representation
Theory and Harmonic Analysis on Semisimple Lie Groups\/}  (P. J.~Sally, 
Jr.~and D.\ A.~Vogan, eds.),  {\it Math.\  Surveys and Monographs\/}, 
A.\ M.\ S. {\bf 31} (1989), 101--170.

[La5]  \bibline,
 {\it On Artin's $L$-functions\/},
 Rice Univ.\  Studies {\bf 56}, Houston, TX, 1970.

[La6] \bibline,
{\it  Euler Products\/},
Yale Univ.\  Press, New Haven, CT, 1971.

[LRS1]  \name{W.\ Luo, Z.\ Rudnick}, and \name{P.\ Sarnak},
On the generalized Ramanujan conjecture for $\text{GL}(n)$,
{\it  Proc.\  Sympos.\  Pure Math\/}.\ {\bf 66} (1999), 301--310.

[LRS2] \bibline,
 On Selberg's eigenvalue conjecture,
{\it Geom.\ Funct.\ Anal\/}.\ {\bf 5} (1995),  387--401.

[MW1]  \name{C.\ M\oe glin} and \name{J.-L.\ Waldspurger},
{\it  Spectral Decomposition and Eisenstein Series, une Paraphrase de
l'Ecriture}, {\it Cambridge Tracts in Math\/}.\ {\bf 113}, 
 Cambridge Univ.\  Press, Cambridge,  1995.

[MW2] \bibline, 
Le spectre r\'esiduel de $\text{GL}(n)$,
 {\it Ann.~Sci.\ E\'cole~Norm.~Sup\/}.\ {\bf 22} (1989),  605--674.

[M] \name{Y.\ Motohashi},
The binary additive divisor problem,
{\it  Ann.\ Sci.\ \'Ecole Norm.\ Sup\/}.\ {\bf 27} (1994), 529--572.

[R1] \name{D.~Ramakrishnan},  Modularity of the
Rankin-Selberg $L$-series, and multiplicity one for $\text{SL}(2)$,
{\it Ann.~of Math\/}.\ {\bf 152} (2000),  45--111.

[R2] \bibline, 
 On the coefficients of cusp forms, {\it Math.\ Research Letters\/} {\bf 4} (1997), 295--307.

[R3] \bibline, 
Modularity of solvable Artin representations of $\text{GO}(4)$-type,
{\it IMRN} (2002), No.\ 1, 1--54.

[Sa]  \name{P.\ Sarnak},
 Kloosterman sums, quadratic forms and modular forms,
{\it  Nieuw\break Archief voor Wiskunde\/} {\bf 1} (2000),  385--389.

[Se]  \name{A.\ Selberg},
 On the estimation of Fourier coefficients of modular forms, {\it Proc.\ Sympos.\ Pure Math\/}.\
 {\bf 8} (1965),  1--15.

[Sh1]  \name{F.\ Shahidi},
 A proof of Langlands conjecture on Plancherel measures; Complementary series for $p$-adic groups,
 {\it  Ann.\  of Math\/}.\ {\bf 132} (1990),  273--330.

[Sh2] \bibline, On the Ramanujan conjecture and finiteness of
poles for certain $L$-functions, {\it Ann.\ of Math\/}.\ {\bf 127}
(1988), 547--584.

[Sh3] \bibline,  Third symmetric power $L$-functions for $\text{GL}(2)$,
{\it  Compositio Math\/}.\ {\bf 70} 
(1989),  245--273.

[Sh4] \name{F.\ Shahidi}, On certain $L$-functions, {\it Amer.\ J.\ Math\/}.\ {\bf 103}
(1981), 297--355.

[Sh5] \bibline,  Fourier transforms of intertwining operators and Plancherel measures for 
$\text{GL}(n)$,
{\it  Amer.\ J.\  Math\/}.\ {\bf 106} 
(1984), 67--111.

[Sh6] \bibline, Twists of a general class of $L$-functions by
highly ramified characters, {\it   Canadian Math.\ Bull\/}.\ {\bf 43} (2000),  380--384.

[Sh7] \bibline, 
 Local coefficients as Artin factors for real groups,
{\it  Duke Math.\ J\/}.\ {\bf 52} (1985),  973--1007.

[Sh8] \bibline,  Symmetric power $L$-functions for
$\text{GL}(2)$, in {\it  Elliptic Curves and Related Topics\/} 
(H.~Kisilevsky and
M.\ R.~Murty, eds.),  {\it CRM Proc.~and Lectures Notes in Math\/}.\ {\bf 4} (1994),  159--182,
 Montr\'eal.

[Tad]  \name{M.~Tadic},
 Classification of unitary representations in irreducible representations of general linear group 
 (non-Archimedean case), {\it  Ann.~Sci.~{\rm \'{\it E}}cole Norm.~Sup\/}.\ {\bf 19} (1986),  335--382.

[Tay] \name{R. Taylor}, Icosahedral Galois representations
(in volume dedicated to Olga Taussky-Todd), {\it Pacific J.\ Math\/}.\ 
{\bf 1997} (Special Issue), 337--347.

[Tu] \name{J.\ Tunnell}, Artin's conjecture for representation of
octahedral type, {\it Bull.\ A.M.S\/}.\ {\bf 5} (1981),  173--175.

[V]  \name{D.\ Vogan},  Gelfand-Kirilov dimension for Harish-Chandra
modules, {\it Invent.\ Math\/}.\ {\bf 48} (1978),  75--98.

[W] \name{A.\ Wiles}, Modular elliptic curves and Fermat's last theorem,
{\it  Ann.\  of Math\/}.\ {\bf 142} (1995),  443--551.

[Z] \name{Y.\ Zhang},   The holomorphy and nonvanishing of normalized
local intertwining operators, {\it  Pacific J.\ Math\/}.\ {\bf 180} (1997),  385--398.

\endreferences